\documentclass[preprint,3p]{elsarticle}
\usepackage[toc,page,title,titletoc,header]{appendix}
\usepackage{stmaryrd}
\SetSymbolFont{stmry}{bold}{U}{stmry}{m}{n}
\usepackage{amsmath}
\usepackage{amssymb}
\usepackage{hyperref}
\hypersetup{
    colorlinks,
    allcolors={green!50!black},
    urlcolor={red!50!black}
}
\usepackage{url}
\usepackage{algorithm}
\usepackage{amsthm}
\usepackage{tabularx}
\usepackage{tikz} 
\usepackage{tikz-cd}
\usepackage{epsfig}
\usepackage{indentfirst}
\usepackage{cases}
\biboptions{sort&compress}
\usepackage{graphicx}
\usepackage{epstopdf}
\usepackage[nameinlink]{cleveref}
\usepackage{color}
\usepackage{float}
\usepackage{algorithmic}
\usepackage{bm}
\usepackage{multirow}
\usepackage{colonequals}
\usepackage{subcaption}

\newcommand{\be}{\begin{equation}}
\newcommand{\ee}{\end{equation}}
\newcommand{\ba}{\begin{array}}
\newcommand{\ea}{\end{array}}
\newcommand{\bea}{\begin{eqnarray}}
\newcommand{\eea}{\end{eqnarray}}
\newcommand{\beas}{\begin{eqnarray*}}
\newcommand{\eeas}{\end{eqnarray*}}

\newtheorem{thm}{Theorem}[section]
\newtheorem{lemma}{Lemma}[section]

\newtheorem{remark}{Remark}[section]

\numberwithin{equation}{section}
\renewcommand{\theequation}{\arabic{section}.\arabic{equation}}

\crefname{lemma}{lemma}{lemmas}
\Crefname{lemma}{Lemma}{Lemmas}
\Crefname{thm}{Theorem}{Theorems}
\crefname{remark}{remark}{remarks}
\Crefname{remark}{Remark}{Remarks}

\def\R{{\mathbb R}}
\def\K{{\mathbb K}}
\def\V{{\mathbb V}}
\def\W{{\mathbb W}}
\def\P{{\mathbb P}}

\newcommand{\p}{\partial}

\newcommand{\cE}{\mathcal E}

\def\id{\mathrm{id}}

\renewcommand{\l}{\left}
\renewcommand{\r}{\right}

\def\cJ{{\mathcal J}}

\renewcommand{\d}{\mathrm{d}}

\newcommand\bX{{\boldsymbol{X}}}
\newcommand\bx{{\boldsymbol{x}}}
\newcommand\bu{{\boldsymbol{u}}}

\newcommand\bA{{\boldsymbol{A}}}
\newcommand\bB{{\boldsymbol{B}}}

\newcommand\bv{{\boldsymbol{v}}}
\newcommand\bY{{\boldsymbol{Y}}}

\newcommand\cH{\mathcal{H}}
\newcommand\bn{{\boldsymbol{n}}}
\def\bF{{\bf F}}

\newcommand{\bxHat}{\widehat{\bx}}

\def\qx{{\widehat{\xi}}}
\def\qbx{{\widehat{\boldsymbol{\xi}}}}
\def\qw{{\widehat{\omega}}}

\newcommand{\includegraphicswithlegend}[5][1.0]{%
	\setlength{\unitlength}{#2}%
  \begin{picture}(1,1)%
    \put(0,0){\includegraphics[width=#1\unitlength]{#3}}%
    \put(1.05,0.27){\includegraphics[width=.03\unitlength]{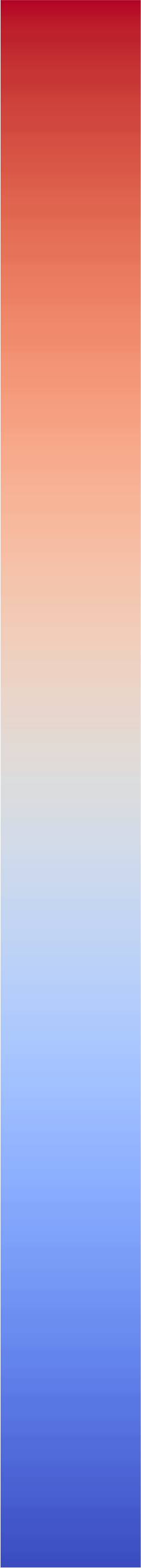}}%
    \put(0.95,0.65){\scriptsize{#4}}%
    \put(0.95,0.17){\scriptsize{#5}}%
  \end{picture}%
}

\renewcommand{\theequation}{\arabic{section}.\arabic{equation}}

\begin{document}

\begin{frontmatter}

\title{Isoparametric finite element methods for \\ mean curvature flow and surface diffusion}

\author[1]{Harald Garcke}
\ead{harald.garcke@ur.de}
\address[1]{Fakult{\"a}t f\"ur Mathematik, Universit{\"a}t Regensburg, 93040 Regensburg, Germany}

\author[2]{Robert N\"urnberg}
\ead{robert.nurnberg@unitn.it}
\address[2]{Dipartimento di Mathematica, Universit\`a di Trento,
38123 Trento, Italy}

\author[3]{Simon Praetorius}
\ead{simon.praetorius@tu-dresden.de}
\address[3]{Institute of Scientific Computing, Technische Universit\"at Dresden, 01062 Dresden, Germany}

\author[4]{Ganghui Zhang}
\ead{zhangganghui2025@gmail.com}
\address[4]{Yau Mathematical Sciences Center, Tsinghua University, Beijing, 100084, China}


\begin{abstract}
We propose higher-order isoparametric finite element approximations for mean curvature flow and surface diffusion. The methods are natural extensions of the piecewise linear finite element methods introduced by Barrett, Garcke, and N\"urnberg (BGN) in a series of papers in 2007 and 2008. The proposed schemes exhibit unconditional energy stability and inherit the favorable mesh quality of the original BGN methods. Moreover, in the case of surface diffusion we present structure-preserving higher-order isoparametric finite element methods. In addition to being unconditionally stable, these also conserve the enclosed volume. Extensive numerical results demonstrate the higher-order spatial accuracy, the unconditional energy stability, the volume preservation for surface diffusion, and the good mesh quality.

\end{abstract}




\end{frontmatter}
\section{Introduction}

Numerical methods for geometric evolution problems have been extensively studied over the last few decades, with a particular focus on gradient flows such as mean curvature, surface diffusion and Willmore flow. One of the most important numerical methods is the parametric finite element method, first proposed by Dziuk \cite{Dziuk1991} using the linear finite element method to solve mean curvature flow in three dimensions. Since then, the parametric finite element method has been widely developed for a range of geometric evolution problems, as discussed in the survey articles \cite{BGN20,Deckelnick-Dziuk-Elliott}.

The majority of contributions to the parametric finite element method for geometric evolution problems have focused on the use of piecewise linear finite elements. Early pioneers for the use of higher-order parametric finite elements have been Schmidt \cite{Schmidt93,Schmidt96} and Heine \cite{Heine03PhD,Heine04}. 
In recent years, higher-order finite element methods have also been studied for elliptic problems on a stationary surface \cite{Demlow}, parabolic partial differential equations on evolving surfaces \cite{Kovacs}, as well as partial differential equations involving coupled bulk-surface systems \cite{Elliott2013,Elliott}. Significant progress has also been made in the numerical analysis of higher-order numerical schemes for geometric evolution problems. Building upon  Huisken's evolution equations, Kov{\'a}cs, Lubich, and Li introduced evolving finite element methods under a matrix formulation for mean curvature flow \cite{Kovacs-Li-Lubich2019} and proved its convergence for higher-order finite elements. Their idea was then later extended to other geometric flows, such as the Willmore flow \cite{Kovacs-Li-Lubich2021}, forced mean curvature flow \cite{Kovacs-Li-Lubich2020}, and generalized mean curvature flow \cite{Binz}. For Dziuk's original scheme of mean curvature flow, Li \cite{Li1,Li2} has established convergence results for curves and surfaces with higher-order finite elements. To address the issue of mesh quality, Hu and Li \cite{Hu-Li} proposed a convergent higher-order method with the introduction of an artificial tangential velocity. More recently, Bai and Li considered a new approach for Dziuk's scheme and proved convergence based on the projected distance \cite{Bai2023}. A similar idea for the convergence of a stabilized BGN scheme was presented in \cite{Bai2024, Bai-Hu-Li}. It is noteworthy that the aforementioned numerical methods leverage higher-order techniques to ensure stability estimates through the use of inverse inequalities. However, as far as the authors know, there is still a lack of proof for the unconditional energy stability and structure-preserving properties in the context of higher-order isoparametric finite elements for mean curvature flow and surface diffusion.

The primary goal of this paper is to construct unconditionally energy-stable higher-order isoparametric finite element methods for mean curvature flow and surface diffusion, applicable to both curve and surface evolution problems. This work extends the results established for linear finite elements by Barrett, Garcke, and N\"urnberg (BGN) in \cite{BGN07A,BGN07B,BGN08A} to higher-order elements. The proposed method is built upon the isoparametric surface finite elements framework presented in \cite{Elliott,Heine04,Demlow}. Higher-order polynomials piecewisely parametrize the approximating curve or surface over the flat reference grid, and the simulation can be implemented by determining the movement of the Lagrange points \cite{Elliott,Praetorius2022}. A key observation is that the difference between the curve/surface elements is controlled by a quantity determined by the local flow map. Specifically, in \Cref{lemma:curve,lemma:surface} we prove the following pointwise estimates for the curve evolution case and surface evolution case, respectively:
\begin{align*}
		\big[\nabla_s\bX\cdot \nabla_s(\bX-\id)\big]\circ \bY \, |\p_{\widehat{\bm{\tau}}}\bY|&\ge |\p_{\widehat{\bm{\tau}}}\widetilde{\bY}|-|\p_{\widehat{\bm{\tau}}}\bY|,\ \ \quad \text{for curves,}\\
	\big[\nabla_{\sigma} \bX: \nabla_{\sigma}(\bX-\id)\big]\circ \bY \, |\cJ(\bY)| &\ge  |\cJ(\widetilde{\bY})| - |\cJ(\bY)|,\quad \text{for surfaces,}
\end{align*}
where  $\nabla_s$ and $\nabla_{\sigma}$ is the curve/surface gradient, $\bY$ and $\widetilde{\bY}$ are the parametrization of the $C^1$-pieces $\sigma$ and $\widetilde{\sigma}$ over the reference grid, $|\p_{\widehat{\bm{\tau}}}\bY|$ and $|\cJ(\bY)|$ are the curve/surface element over the reference element, and $\bX=\widetilde{\bY}\circ\bY^{-1}\colon\sigma\to \widetilde{\sigma}$ is the local transformation map.

These estimates enable us to integrate over the reference element and derive the corresponding energy stability results. In contrast to the linear finite element case \cite{BGN07A,BGN07B,BGN08A}, the normal vector and the curve/surface integration element over the approximating curve/surface are not piecewise constant, nor even polynomial functions. Therefore, the use of numerical integration and appropriate quadrature rules is necessary for the implementation of these higher-order methods. Another key advantage of the proven pointwise estimates is that they lead to numerical energy stability under any quadrature rule with positive weights. This demonstrates the practical robustness of the proposed higher-order methods, as the energy stability is guaranteed regardless of the specific quadrature rules employed.

It is well-known that surface diffusion exhibits two important geometric structures: (i) the perimeter/area of the curve/surface is decreasing, and (ii) the area/volume of the region enclosed by the curve/surface is preserved. Much effort has been devoted to designing numerical methods to preserve both key geometric structures. 
The original BGN method \cite{BGN07A,BGN08A} satisfies an unconditional energy stability result on the fully discrete level, as well as stability and volume conservation for a semidiscrete approximation.
One of the pioneering works towards fully discrete structure-preserving methods is the study by Jiang and Li \cite{Jiang21}. Later, Bao and Zhao \cite{Bao-Zhao} proposed a weakly implicit method that approximates the normal vector using information from the adjacent time step. This idea has been widely adopted and extended to other scenarios, including axisymmetric geometric evolution equations \cite{Bao-Garcke-Nurnberg-Zhao2022}, anisotropic surface diffusion \cite{Bao-Jiang-Li,Bao-Li2023,Bao-Li2024}, two-phase Navier--Stokes flow \cite{Garcke-Nurnberg-Zhao2023}, or the Mullins--Sekerka evolution \cite{Nurnberg}.

The next goal of this paper is to design a structure-preserving isoparametric finite element method of surface diffusion, applicable to both curve and surface evolution problems. To achieve this,  we employ the idea introduced by Bao and Zhao \cite{Bao-Zhao}, which uses the approximation of the normal vector based on the current and next-time step information. A key observation in our work is that under the parametrization of the reference element, both the integrands of the area/volume functional and the time derivative of the area/volume of the curve/surface are polynomial functions. Therefore, we can obtain the exact preservation of the area/volume property by using sufficiently accurate numerical integration and quadrature rules.

Finally, we would like to remark that while this paper focuses on high-order spatial discretizations, the time discretization remains of low order. Recent developments of high-order, energy-decaying time discretization methods for geometric flows include an averaged vector-field collocation method \cite{Duan2021}, a Lagrange multiplier approach \cite{Garcke2024}, and a Crank--Nicolson type predictor-corrector scheme \cite{csftime}.

The remainder of the paper is organized as follows: Section 2 focuses on the curve evolution problem. We first review the preliminaries of isoparametric finite element methods for curves and then prove the unconditional energy stability of the proposed schemes. Furthermore, we introduce a structure-preserving numerical method for the surface diffusion of curves.  Section 3 extends the constructions and proofs presented in Section 2 to the surface evolution case. The development of the isoparametric finite element methods and the corresponding stability and structure-preservation analyses are conducted for the surface case. Section 4 presents extensive numerical results to illustrate the accuracy and robustness of the proposed higher-order isoparametric finite element methods for both curve and surface evolution problems. Finally, in Section 5, we draw some conclusions.

\section{Curve evolution}

In this section, we consider curvature flow and surface diffusion for a family of planar curves $\Gamma(t)$, determined by
\begin{align*}
	\mathcal{V}=\begin{cases}
		\kappa,\qquad &\text{curvature flow of curves},\\
		-\partial_{ss}\kappa,\qquad\ &\text{surface diffusion of curves},
	\end{cases}
\end{align*}
where $\mathcal{V}$ is the normal velocity along the curve $\Gamma$, $\kappa$ is the curvature and $s$ is the arc-length.

\subsection{Isoparametric finite element approximation}

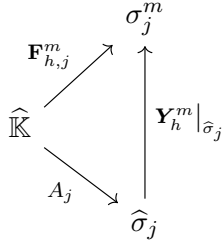
\begin{figure}[ht]
\[
    \begin{tikzcd}
        & \sigma^m_j\\
        \widehat{\K}\arrow[ru,rightarrow,"\bF_{h,j}^m"]\arrow[rd,rightarrow,"A_{j}"'] &\\
        & \widehat{\sigma}_j\arrow[uu,rightarrow,"\l.\bY^m_h\r|_{\widehat{\sigma}_j}"']
    \end{tikzcd}
    \]
    \caption{\label{fig:curve-mapping-diagram}Mappings between the reference element $\widehat{\K}$, an element of the reference grid $\widehat{\sigma}_j$, and a curved element in the parametrized grid $\sigma^m_j$.}
\end{figure}

We define higher-order analogues of the piecewise linear finite elements studied in \cite{BGN07A}. We follow the framework in \cite[Section 6]{Elliott} and first define the reference finite element triple $(\widehat{\K},\widehat{\V},\widehat{\Sigma})$. Here $\widehat{\K}=[0,1]\subset \R$ is the reference 1-simplex, $\widehat{\V}=\P^{\ell}(\widehat{\K})$ is a finite-dimensional space consisting of polynomial functions of degree $\le \ell$ over $\widehat{\K}$, $\widehat{\Sigma}\subset \widehat{\K}$ is the set of  Lagrange nodes as the dual space of $\widehat{\V}$, which determines $\widehat{\V}$ in the sense that if $\widehat{\varphi}(\widehat{x})=0$ for all $\widehat{x}\in \widehat{\Sigma}$, then $\widehat{\varphi}=0\in \widehat{\V}$.  Define a reference grid $\widehat{\Gamma}_h$ consisting of flat segments $\widehat{\sigma}_j$ as the images of the reference simplex by affine linear maps $A_j$
\[
\widehat{\Gamma}_h=\bigcup_{j=1}^N \widehat{\sigma}_j,\qquad  \widehat{\sigma}_j=A_j(\widehat{\K}),\qquad h=\max_{j=1}^N\mathrm{diam}(\widehat{\sigma}_j),
\]
where $\mathrm{diam}(\widehat{\sigma}_j)$ is the length of the segment $\widehat{\sigma}_j$. To approximate the evolution of closed curves  $\Gamma(t_m)$, $t_m\in [0,T]$, by higher order polynomials, following \cite{Demlow,Heine04}, we look for the parametric representation $\Gamma^m_h=\bY^m_h(\widehat{\Gamma}_h)$ over the reference grid $\widehat{\Gamma}_h$. Equivalently, we search for the piecewise $\P^\ell$ functions $\bF_{h,j}^m \colonequals \l.\bY^m_h\r|_{\widehat{\sigma}_j}\circ A_j$ over the reference simplex, and decompose $\Gamma^m_h=\bigcup_{j=1}^N \sigma^m_j$ into piecewise curved line segments, see \Cref{fig:curve-mapping-diagram}.

For simplicity, we fix $\tau>0$ as a uniform time step size and set $t_m=m\tau$, $m\ge 0$. Suppose for $\Gamma(t_m)$ we already have an approximating curve $\Gamma^m_h=\bY^m_h(\widehat{\Gamma}_h)$ with a piecewise $\P^\ell$ function $\l.\bY^m_h\r|_{\widehat{\sigma}_j}\in \P^{\ell}(\widehat{\sigma}_j)$, we define the surface finite element space $\V_\ell^m$ over $\Gamma^m_h$ as
\begin{align*}
    \V_\ell^m
    &\colonequals \l\{v_h\in H^1(\Gamma_h^m)\,:\, v_h\circ  \l.\bY_h^m\r|_{\widehat{\sigma}_j}\in \P^{\ell}(\widehat{\sigma}_j), \; j=1,\ldots,N\r\}.
\end{align*}
We consider the following formulations for two geometric evolution equations of curves \cite{BGN07A,BGN07B}:  mean curvature flow (MCF) and  surface diffusion (SD) of curves
    \begin{align*}
        \p_t\bY\circ\bY^{-1} \cdot \bn&=\begin{cases}
            \kappa,\quad &\text{MCF in 2d},\\
            -\partial_{ss}\kappa,\quad &\text{SD in 2d},
        \end{cases}\\
        \kappa\bn&=\partial_{ss}\id,
    \end{align*}
where $\bY(t):\widehat{\Gamma}\rightarrow \Gamma(t)$ is the parametrization of $\Gamma(t)$, $\widehat{\Gamma}$ is a fixed reference curve with the same topology as $\Gamma(t)$, $\kappa$ is the curvature, $s$ represents the arc-length and $\partial_s$ is the differentiation with respect to the arc-length parametrization. Here, the coupled equations are defined on the evolving curve $\Gamma(t)$ and  $\bn$ denotes the unit outer normal to $\Gamma(t)$. We then consider the following higher-order in space method:

\smallskip

\noindent \textbf{Isoparametric finite element method with exact integration:} Find the solution (called the local flow map) $\bX_h^{m+1}\in [\V_\ell^m]^2$ and $\kappa^{m+1}_h\in \V_\ell^m$ that satisfy the following equations
\begin{subequations}\label{eq:2d}
    \begin{align}
        \l(\frac{\bX^{m+1}_h-\id}{\tau},\varphi_h\bn^{m}_h \r)_{\Gamma^m_h}&=\begin{cases}
        \l(\kappa^{m+1}_h,\varphi_h\r)_{\Gamma^m_h},\quad &\text{MCF in 2d},\\
        \l(\nabla_{s}\kappa^{m+1}_h,\nabla_{s}\varphi_h\r)_{\Gamma^m_h},\quad &\text{SD in 2d,}
        \end{cases} \label{eq:2d1} \\
        \l(\kappa^{m+1}_h\bn^m_h,\bm{\eta}_h \r)_{\Gamma^m_h} &+\l(\nabla_s\bX^{m+1}_h,\nabla_s \bm{\eta}_h \r)_{\Gamma^m_h}=0,\label{eq:2d2}
    \end{align}
\end{subequations}
for any $\varphi_h\in \V_\ell^m$ and $\bm{\eta}_h\in [\V_\ell^m]^2$. Here, the unit normal $\bn^m_h$ and the discrete differential operator $\nabla_s$ over $\Gamma^m_h$ are all defined piecewise as
\begin{equation*} 
    \bn^m_h \circ \bY_h^m=\frac{(\p_{\widehat{\bm{\tau}}} \bY_h^m)^\perp}{|\p_{\widehat{\bm{\tau}}}\bY_h^m|} \ \text{on}\ \widehat\Gamma_h,
    \quad\text{and}\quad \nabla_sf=\p_{\bm{\tau}^m}f\, \bm{\tau}^m \ \text{on}\ \Gamma^m_h,
\end{equation*}
 where $\widehat{\bm{\tau}}$ is the unit tangent vector of $\widehat{\Gamma}_h$, $\p_{\widehat{\bm{\tau}}}$ is the directional derivative with respect to $\widehat{\bm{\tau}}$, and $\bm{\tau}^m$ is the unit tangent vector of $\Gamma^m_h$. Moreover, the inner product $(\cdot,\cdot)_{\Gamma^m_h}$ is the $L^2$--inner product over $\Gamma^m_h$, defined as
\begin{equation}\label{exact integration}
(u,v)_{\Gamma^m_h}
	 =\int_{\Gamma^{m}_h} u\cdot v\ \d s =  \int_{\widehat{\Gamma}_h} (u\circ \bY^m_h )\cdot (v\circ \bY^m_h)\, |\p_{\widehat{\bm{\tau}}}\bY^m_h|\, \d s.
\end{equation}
Now, suppose we have the solution of \eqref{eq:2d}. Then we define $\bY^{m+1}_h \colonequals \bX^{m+1}_h\circ \bY^{m}_h$, and $\Gamma^{m+1}_h \colonequals \bY^{m+1}_h(\widehat{\Gamma}_h)$ is the desired higher-order isoparametric approximation of $\Gamma(t_{m+1})=\bY(t_{m+1})(\widehat{\Gamma})$.

In practice, performing the integration in \eqref{exact integration} exactly for $\ell\ge 2$ is not easily possible, since the integrand is not a polynomial. Therefore, we now investigate using numerical integration with the help of suitable quadrature rules. Let $\{\qx_i\}_{i=1}^Q$ be the distinct quadrature points and let $\{\qw_{i}\}_{i=1}^Q$ be the corresponding quadrature weights. Then we define the following numerical integration over the reference element
\begin{equation}\label{quadrature rule 1d}
    \int^h_{\widehat{\K}}f(\widehat{x})\ \d \widehat{x} \colonequals  \sum_{i=1}^{Q}f(\qx_i)\ \qw_{i}.
\end{equation}
Similarly, we can introduce the numerical integration over each of the flat elements $\widehat{\sigma}_j$ of the reference grid $\widehat{\Gamma}_h$ as follows,
\begin{equation*}
    \int^h_{\widehat{\sigma}_j}f(x)\ \d s \colonequals \sum_{i=1}^Q (f\circ A_j) (\qx_i)\, |\nabla A_j| (\qx_i) \ \qw_{i}.
\end{equation*}

We say the numerical integration \eqref{quadrature rule 1d} is exact to degree $p$, if for any $f\in \P^{p}(\widehat{\K})$ it holds that $\int^h_{\widehat{\K}}f(\widehat{x})\ \d\widehat{x}=\int_{\widehat{\K}}f(\widehat{x})\ \d\widehat{x}$. By quadrature over the reference element, we introduce the inner product $(u,v)_{\Gamma^m_h}^h$ as 
\begin{align}
    (u,v)^h_{\Gamma^m_h}
    & \colonequals  \int_{\widehat{\Gamma}_h}^h (u\circ \bY^m_h )\cdot (v\circ \bY^m_h)\, |\p_{\widehat{\bm{\tau}}}\bY^m_h|\, \d s\notag\\
    &\colonequals\sum_{j=1}^N\int^h_{\widehat{\sigma}_j}(u\circ \bY^m_h )\cdot (v\circ \bY^m_h)|\p_{\widehat{\bm{\tau}}}\bY^m_h|\ \d s \notag\\
    &=\sum_{j=1}^N\int^h_{\widehat{\K}} (u\circ\bF^m_{h,j})\cdot  (v\circ\bF^m_{h,j})\  |\nabla \bF^m_{h,j}| \ \d \widehat{x} \notag \\
    &=\sum_{j=1}^N\sum_{i=1}^Q (u\circ\bF^m_{h,j})(\qx_i)\cdot  (v\circ\bF^m_{h,j})(\qx_i)\  |\nabla \bF^m_{h,j}(\qx_i)|\  \qw_i, \label{Numerical integration}
\end{align}
where the surface element $|\nabla \bF^m_{h,j}(\qx_i)|$ is defined pointwise as the norm of the vector $\nabla \bF^m_{h,j}\in \R^{2}$, see also the analogous definitions in \cite{Bai2023,Bai2024}.
In the above, we have
used the definition $\bF_{h,j}^m = \bY^m_h\circ A_j$ and the chain rule to get $|\nabla \bF^m_{h,j}|=|\p_{\widehat{\bm{\tau}}}\bY^m_h||\nabla A_j|$.

In practice, classical Newton--Cotes, Gauss or Gauss--Lobatto quadrature rules can be used. However, to obtain unconditionally stable schemes with numerical integration,  the quadrature weights need to be positive. Observe that \eqref{quadrature rule 1d} with $Q=2$ and with the quadrature points $\qx_1=0,\, \qx_2=1$ and weights $\qw_1=\qw_2=\frac{1}{2}$ corresponds to the lumped mass inner product considered in \cite{BGN07A}. With the numerical integration in \eqref{Numerical integration} at hand, we then consider the following higher-order method in practice:

\smallskip

\noindent \textbf{Isoparametric finite element method with numerical integration:} Find the solution (called the local flow map) $\bX_h^{m+1}\in [\V^m_\ell]^2$ and the approximative curvature $\kappa^{m+1}_h\in \V^m_\ell$ that satisfy the following equations
\begin{subequations}\label{eq:2dnum}
    \begin{align}
        \l(\frac{\bX^{m+1}_h-\id}{\tau},\varphi_h\bn^{m}_h \r)^h_{\Gamma^m_h}&=\begin{cases}
            \l(\kappa^{m+1}_h,\varphi_h\r)^h_{\Gamma^m_h},\quad &\text{MCF in 2d},\\
            \l(\nabla_{s}\kappa^{m+1}_h,\nabla_{s}\varphi_h\r)^h_{\Gamma^m_h},\quad &\text{SD in 2d,}
        \end{cases} \label{eq:2d1num} \\
        \l(\kappa^{m+1}_h\bn^m_h,\bm{\eta}_h \r)^h_{\Gamma^m_h} &+\l(\nabla_s\bX^{m+1}_h,\nabla_s \bm{\eta}_h \r)_{\Gamma^m_h}^h=0,\label{eq:2d2num}
    \end{align}
\end{subequations}
for any $\varphi_h\in \V_\ell^m$ and $\bm{\eta}_h\in [\V_\ell^m]^2$.
Combined with the numerical integration \eqref{Numerical integration}, we only need to evaluate the normal vector and surface gradient operator at the quadrature points.

\subsection{Well-posedness and unconditional energy stability}

In this subsection, we establish the well-posedness of \eqref{eq:2d} and \eqref{eq:2dnum} under mild conditions, as well as the unconditional energy stability for both schemes.
\begin{thm}[Well-posedness for \eqref{eq:2d}]\label{Well-posed:2d:exact}
    Assume that the discrete curve $\Gamma^m_h$ satisfies the following mild conditions \begin{enumerate}
    \item[$(\mathcal{A}1)$] The curve elements are nondegenerate, i.e., for each $\widehat{\sigma}_j$, $j=1,\ldots,N$, it holds
    \begin{equation*}
	|\p_{\widehat{\bm{\tau}}} \bY^m_h|> 0,\quad  \text{on} \quad \widehat{\sigma}_j.
    \end{equation*}
    \item[$(\mathcal{A}2)$] Nonparallel assumption, that is
    \begin{equation*}\label{Well-posedness:assumption:A2}
        \mathrm{dim}\ \mathrm{span} \l\{\l(\bn^m_h\varphi_h,1 \r)_{\Gamma^m_h}| \, \varphi_h\in \V_\ell^m \r\}= 2.
    \end{equation*}
    \end{enumerate}
    Then the system \eqref{eq:2d} has a unique solution.
\end{thm}

\begin{proof}
	We first consider the linear algebraic systems \eqref{eq:2d}. It suffices to show that the following homogeneous equation  only has the  zero solution 
in 
$[\V^m_\ell]^2\times \V^m_\ell$. We hence require for any $\varphi_h\in \V^m_\ell$ and $\bm{\eta}_h\in [\V^m_\ell]^2$
\begin{subequations}
	\begin{align}
		\l(\frac{\bX_h}{\tau},\varphi_h\bn^{m}_h \r)_{\Gamma^m_h}&=\begin{cases}
			\l(\kappa_h, \varphi_h\r)_{\Gamma^m_h},\quad &\text{MCF},\\
			\l(\nabla_s\kappa_h, \nabla_s\varphi_h\r)_{\Gamma^m_h},\quad &\text{SD},
		\end{cases}\label{eq:hom1} \\
			\l(\kappa_h\bn^m_h,\bm{\eta}_h \r)_{\Gamma^m_h} +&\l(\nabla_s\bX_h,\nabla_s \bm{\eta}_h \r)_{\Gamma^m_h}=0.\label{eq:hom2}
	\end{align}
\end{subequations}
Indeed, taking $\varphi_h=\kappa_h$ in \eqref{eq:hom1} and $\bm{\eta}_h =\bX_h$ in \eqref{eq:hom2}, we obtain
\begin{equation*}
	\l(\nabla_s\bX_h,\nabla_s\bX_h \r)_{\Gamma^m_h}=\begin{cases}
		-\tau\l(\kappa_h,\kappa_h\r)_{\Gamma^m_h},\quad &\text{MCF},\\
		-\tau\l(\nabla_s\kappa_h,\nabla_s\kappa_h\r)_{\Gamma^m_h},\quad &\text{SD}.
	\end{cases}
\end{equation*}
Using the non-degeneracy assumption $(\mathcal{A}1)$, this leads to $\kappa_h=0$ and $\bX_h=\bX^c$ a constant vector in the 
curvature flow case and $\kappa_h=\kappa^c$ and $\bX_h=\bX^c$ a constant vector in the surface diffusion case. Now substituting $\kappa_h$ and $\bX_h$ back into \eqref{eq:hom1} and \eqref{eq:hom2}, we arrive at
\[
	0
		=\l(\bX_h,\varphi_h\bn^{m}_h \r)_{\Gamma^m_h}=\bX^c\cdot \l(1,\varphi_h\bn^{m}_h \r)_{\Gamma^m_h},
\]
which implies that $\bX_h=0$ by assumption $(\mathcal{A}2)$. It remains to show that $\kappa^c=0$ in the surface diffusion case.
To this end, we choose a $\varphi_{h,0} \in \V^m_\ell$ such that $\bm v_0 \colonequals \l(\bn^m_h\varphi_{h,0},1 \r)_{\Gamma^m_h}$ is nonzero,
which is possible thanks to assumption $(\mathcal{A}2)$. Then we choose $\bm\eta_h=\kappa^c \bm v_0 \varphi_{h,0}$ in \eqref{eq:hom2},
and obtain
\[
0 = \kappa^c \l(\bn^m_h,\bm{\eta}_h \r)_{\Gamma^m_h} 
= (\kappa^c)^2 \l(\bn^m_h, \bm v_0 \varphi_{h,0} \r)_{\Gamma^m_h} 
= (\kappa^c)^2 \bm v_0 \cdot \l(\bn^m_h, \varphi_{h,0} \r)_{\Gamma^m_h}
= (\kappa^c)^2 |\bm v_0|^2.
\]
Hence $\kappa^c=0$, which completes the proof.
\end{proof}

A similar argument can be applied to the practical scheme \eqref{eq:2dnum}.

\begin{thm}[Well-posedness for \eqref{eq:2dnum}]\label{Well-posed:2d:num}
Assume that the discrete curve $\Gamma^m_h$ satisfies the following mild conditions
\begin{enumerate}
 \item[$(\mathcal{A}1^h)$] The curve elements are nondegenerate over the quadrature points, i.e., for each $\widehat{\sigma}_j$, $j=1,\ldots,N$, it holds
	\begin{equation*}
	|\p_{\widehat{\bm{\tau}}} \bY^m_h|(x_{j_i})> 0,\quad \ i=1,\ldots,Q,
\end{equation*}
where $\{x_{j_i}\}_{i=1}^{Q}\subset \widehat{\sigma}_j$ is the image of $\{\qx_i\}_{i=1}^{Q}$ under the affine linear map $A_j$. Equivalently,
\begin{equation*}
   |\nabla \bF^m_{h,j}(\qx_i)|>0,\quad j=1,\ldots,N,\ i=1,\ldots,Q.
\end{equation*}
	\item[$(\mathcal{A}2^h)$] Nonparallel assumption for numerical integration, that is
\begin{equation*}
		\mathrm{dim}\ \mathrm{span} \l\{\l(\bn^m_h\varphi_h,1 \r)^h_{\Gamma^m_h}| \, \varphi_h\in \V_\ell^m \r\}= 2.
	\end{equation*}
\end{enumerate}
We further assume that the quadrature rule \eqref{quadrature rule 1d} satisfies 
\begin{itemize}
     \item[$(\mathcal{A}3^h)$] It holds that $Q\ge \ell+1$ and $\qw_i>0$ for $i=1,\ldots,Q$.
\end{itemize}
Then, the system \eqref{eq:2dnum} has a unique solution.
\end{thm}

\begin{proof}
	We first notice that for any $f_h\in \V^m_{\ell}$, we have  $(f_h,f_h)^h_{\Gamma^m_h}\ge 0$ with equality if and only if $f_h=0.$ Indeed, by definition \eqref{Numerical integration}, assumption $(\mathcal{A}1^h)$ and $(\mathcal{A}3^h)$, we have
\begin{align*}
    (f_h,f_h)^h_{\Gamma^m_h}
    &=\sum_{j=1}^N\sum_{i=1}^Q |(f_h\circ\bF^m_{h,j})(\qx_i)|^2\  |\nabla \bF^m_{h,j}(\qx_i)|\  \qw_i\ge 0,
\end{align*}
with equality if and only if $(f_h\circ\bF^m_{h,j})(\qx_i)=0$ for all quadrature points $\qx_i$. Since $f_h\circ\bF^m_{h,j}$ is a polynomial of degree $\ell$ and $Q\ge \ell+1$, we have $f_h=0$. Similarly, for any $f_h\in \V^m_{\ell}$, we have  $(\nabla_{s} f_h,\nabla_{s} f_h)^h_{\Gamma^m_h}\ge 0$ with equality if and only if $f_h$ is constant. The rest of the proof proceeds as the proof of \Cref{Well-posed:2d:exact}.
\end{proof}

\begin{remark}\label{Well-posed:remark}
The conditions above are natural generalizations of the mild conditions considered in \cite{BGN07A}. Obviously, the assumption  $(\mathcal{A}1)$ implies  $(\mathcal{A}1^h)$. On the other hand, the assumptions $(\mathcal{A}2)$ and  $(\mathcal{A}2^h)$ are equivalent if the quadrature rule \eqref{quadrature rule 1d} is exact to degree $2\ell-1$.  Indeed, we first note that the normal vector can be represented as
\[
\bn^{m}_h \circ \bF_{h,j}^m =\frac{(\nabla  \bF_{h,j}^m)^\perp }{|\nabla  \bF_{h,j}^m|}
\]
 over the reference simplex. Then, for a quadrature rule \eqref{quadrature rule 1d} that is exact to degree $2\ell-1$ we have
\begin{align*}
	\l(1,\varphi_h\bn^{m}_h \r)_{\Gamma^m_h}
    &= \int_{\widehat{\Gamma}_h}(\varphi_h\circ\bY_h^m)\,  \frac{(\p_{\widehat{\bm{\tau}}} \bY_h^m)^\perp}{|\p_{\widehat{\bm{\tau}}}\bY_h^m|}\, |\p_{\widehat{\bm{\tau}}}\bY^m_h| \ \d s 
	=\sum_{j=1}^N \int_{\widehat{\sigma}_j}(\varphi_h\circ\bY_h^m)\, (\p_{\widehat{\bm{\tau}}} \bY_h^m)^\perp \ \d s  \\
    &=\sum_{j=1}^N \int_{\widehat{\K}}(\varphi_h\circ\bF_{h,j}^m)\, (\nabla  \bF_{h,j}^m)^\perp \ \d\widehat{x}
  =\sum_{j=1}^N \int_{\widehat{\K}}^h(\varphi_h\circ\bF_{h,j}^m)\, (\nabla  \bF_{h,j}^m)^\perp \ \d\widehat{x} \\&
	=\sum_{j=1}^N \int_{\widehat{\K}}^h(\varphi_h\circ\bF_{h,j}^m)\, \frac{(\nabla  \bF_{h,j}^m)^\perp}{|\nabla  \bF_{h,j}^m|} |\nabla  \bF_{h,j}^m| \ \d\widehat{x} =\l(1,\varphi_h\bn^{m}_h \r)_{\Gamma^m_h}^h,
\end{align*}
where we have observed that each component of the vector $(\varphi_h\circ\bF_{h,j}^m)\, (\nabla  \bF_{h,j}^m)^\perp$ is a polynomial of degree $\le 2\ell-1$ over $\widehat{\K}$.

\end{remark}

The following lemma will be crucial to prove a generalization of \cite[Lemma~57]{BGN20} to curved elements. In later applications we will choose $\widehat{\sigma}$ to be flat and $\bY$ and $\widetilde\bY$ to be polynomial.

\begin{lemma}\label{lemma:curve}
	Let $\widehat{\sigma},\sigma,\widetilde{\sigma}\subset \R^2$ be three $C^1$-curves, with or without boundary. Let $\sigma,\widetilde{\sigma}\subset \R^2$ be parametrized by $\bY:\widehat{\sigma}\rightarrow \R^2$ and $\widetilde{\bY}=\bX\circ \bY:\widehat{\sigma}\rightarrow \R^2$, respectively, where $\bX\in C^1(\sigma,\R^2)$. Then the following inequality holds
	\begin{equation} \label{eq:lemma}
		\big[\nabla_s\bX\cdot \nabla_s(\bX-\id)\big]\circ \bY \, |\p_{\widehat{\bm{\tau}}}\bY|\ge |\p_{\widehat{\bm{\tau}}}\widetilde{\bY}|-|\p_{\widehat{\bm{\tau}}}\bY|,\quad \text{on } \widehat{\sigma},
	\end{equation}
where the equality holds if $\bX=\id\in C^1(\sigma,\R^2)$. Here, the curve gradient $\nabla_s$ is defined over $\sigma$.
	
\end{lemma}
	
\begin{proof}
	Let $\bm{\tau}$ be the unit tangent of $\sigma$. Then, by definition, we have
\[
\nabla_{s} \bX= \p_{\bm{\tau}}\bX\otimes \bm{\tau},\quad \text{and}\quad \nabla_{s} \id=\bm{\tau} \otimes \bm{\tau}.
\]
Then we can estimate by using the chain rule
\begin{align*}
	[\nabla_{s}\bX\cdot \nabla_{s}(\bX-\id)]\circ \bY |\p_{\widehat{\bm{\tau}}}\bY|
	&\ge \frac{1}{2}(|\nabla_{s}\bX|^2-|\nabla_{s}\id|^2)\circ \bY |\p_{\widehat{\bm{\tau}}}\bY|\\
	&=  \frac{1}{2}(|\p_{\bm{\tau}}\bX|^2-1)\circ \bY |\p_{\widehat{\bm{\tau}}}\bY|\\
	&\ge (|\p_{\bm{\tau}}\bX|-1)\circ \bY |\p_{\widehat{\bm{\tau}}}\bY|\\
	&=|\p_{\widehat{\bm{\tau}}}\widetilde{\bY}|-|\p_{\widehat{\bm{\tau}}}\bY|,
\end{align*}
where in the first line and the third line we use the basic inequalities for $\bu,\bv\in \R^2$ and $f\in \R$
\[
\bu\cdot (\bu-\bv)\ge \frac{1}{2}(|\bu|^2-|\bv|^2),\quad\text{and}\quad  \frac{f^2-1}{2}\ge f-1,
\]
respectively. This proves \eqref{eq:lemma}. Clearly, \eqref{eq:lemma} holds with equality if $\bX=\id\in C^1(\sigma,\R^2)$. 
\end{proof}	

\begin{thm}[Unconditional energy stability with exact integration]\label{Thm: energy stability for curve}
	Under assumption $(\mathcal{A}1)$, suppose that $(\bX^{m+1}_h,\kappa^{m+1}_h)$ $\in [\V^m_\ell]^2\times \V^m_\ell$ is a solution of \eqref{eq:2d}. Then for any $M\ge 0$, we have
	\begin{equation}\label{Unconditional energy stability:curve}
		\cH^{1}(\Gamma^{M+1}_h) \le \cH^{1}(\Gamma^{0}_h)-\begin{cases}
			\tau\sum_{m=0}^{M}\l(\kappa^{m+1}_h,\kappa^{m+1}_h\r)_{\Gamma^{m}_h},\quad &\text{MCF},\\
			\tau\sum_{m=0}^{M}\l(\nabla_s\kappa^{m+1}_h,\nabla_s\kappa^{m+1}_h\r)_{\Gamma^{m}_h},\quad &\text{SD},
		\end{cases}
	\end{equation}
	where $\cH^{1}(\Gamma^{m}_h) \colonequals \int_{\Gamma^{m}_h}\d s=\int_{\widehat{\Gamma}_h}| \p_{\widehat{\bm{\tau}}}\bY^{m}_h|\ \d s$ is the exact energy (perimeter) of the discrete curve $\Gamma^m_h$.
	
\end{thm}

\begin{proof}
	Take $\varphi_h=\kappa_h^{m+1}\in \V^m_\ell$ in \eqref{eq:2d1} and $\bm{\eta}_h=\bX_h^{m+1}-\id_{|_{\Gamma^m_h}}\in [\V^m_\ell]^2$ in \eqref{eq:2d2}, we obtain
\begin{align*}
	\l(\nabla_s\bX_h^{m+1},\nabla_s(\bX_h^{m+1}-\id)\r)_{\Gamma^m_h}=-\begin{cases}
		\tau\l(\kappa^{m+1}_h,\kappa^{m+1}_h\r)_{\Gamma^m_h}\le 0,\quad &\text{MCF},\\
		\tau\l(\nabla_s\kappa^{m+1}_h,\nabla_s\kappa^{m+1}_h\r)_{\Gamma^m_h}\le 0,\quad &\text{SD}.
	\end{cases}
\end{align*}
Denote $\Gamma^{m+1}_h=\bigcup_{j=1}^N\sigma_j^{m+1}$ and $\Gamma^{m}_h=\bigcup_{j=1}^N\sigma_j^{m}$ be the decomposition of curved elements. For each $j=1,\ldots,N$, we apply \Cref{lemma:curve} for $\widehat{\sigma}=\widehat{\sigma}_j$, $\sigma=\sigma^{m}_j$, $\bY=\bY^{m}_h|_{\widehat{\sigma}_j}$,  $\widetilde{\sigma}=\sigma^{m+1}_j$, $\widetilde{\bY}=\bY^{m+1}_h|_{\widehat{\sigma}_j}=\bX^{m+1}_h\circ \bY^{m}_h|_{\widehat{\sigma}_j}$ and  obtain
\[
\big[\nabla_s\bX^{m+1}_h\cdot \nabla_s(\bX^{m+1}_h-\id)\big]\circ \bY^{m}_h \, |\p_{\widehat{\bm{\tau}}}\bY^{m}_h|\ge |\p_{\widehat{\bm{\tau}}}\bY^{m+1}_h|-|\p_{\widehat{\bm{\tau}}}\bY^{m}_h|,\quad \text{on}\quad \widehat{\sigma}_j.
\]
Integrating over $\widehat{\sigma}_j$ and taking the summation over $j=1,\ldots,N$, using the definition of integration, we arrive at
\begin{multline*}
	\begin{aligned}
 	\l(\nabla_s\bX^{m+1}_h,\nabla_s\l(\bX^{m+1}_h-\id\r)\r)_{\Gamma^m_h}
	&=\int_{\widehat{\Gamma}_h}\big[\nabla_s\bX^{m+1}_h\cdot \nabla_s(\bX^{m+1}_h-\id)\big]\circ \bY^{m}_h \, |\p_{\widehat{\bm{\tau}}}\bY^{m}_h| \ \d s\\
	&\ge \int_{\widehat{\Gamma}_h}(|\p_{\widehat{\bm{\tau}}}\bY^{m+1}_h|-|\p_{\widehat{\bm{\tau}}}\bY^{m}_h|)\ \d s\\
	&=\int_{\widehat{\Gamma}_h}|\p_{\widehat{\bm{\tau}}}\bY^{m+1}_h|\ \d s-\int_{\widehat{\Gamma}_h}|\p_{\widehat{\bm{\tau}}}\bY^{m}_h|\ \d s=\cH^{1}(\Gamma^{m+1}_h)-\cH^{1}(\Gamma^{m}_h).
	\end{aligned}
\end{multline*}
We then complete the proof by induction for $m$. 
\end{proof}

\begin{thm}[Unconditional energy stability with numerical integration]\label{Thm: energy stability for curve num}
Assume that $(\mathcal{A}1^h)$ holds and that the quadrature rule \eqref{quadrature rule 1d} has positive weights $\qw_i>0$ for $i=1,\ldots,Q$. Let $(\bX^{m+1}_h,\kappa^{m+1}_h)$ $\in [\V^m_\ell]^2\times \V^m_\ell$ be a solution of \eqref{eq:2dnum}. Then for any $M\ge 0$, we have
	\begin{equation}\label{Unconditional energy stability:num}
		L^{M+1}_h \le L_h^{0}-\begin{cases}
			\tau\sum_{m=0}^{M}\l(\kappa^{m+1}_h,\kappa^{m+1}_h\r)^h_{\Gamma^{m}_h},\quad &\text{MCF},\\
			\tau\sum_{m=0}^{M}\l(\nabla_s\kappa^{m+1}_h,\nabla_s\kappa^{m+1}_h\r)^h_{\Gamma^{m}_h},\quad &\text{SD},
		\end{cases}
	\end{equation}
where $L_h^{m} \colonequals \l(1,1\r)^h_{\Gamma^m_h}$ is the numerical energy (perimeter) of the discrete curve $\Gamma^m_h$. 
	
\end{thm}

\begin{proof}
	Similarly to the proof of \Cref{Thm: energy stability for curve}, we can obtain the following inequality at each quadrature point $x_{j_i}=\bF^m_{h,j}(\qx_i)$,
\begin{equation*} 
	\big[\nabla_s\bX^{m+1}_h\cdot \nabla_s(\bX^{m+1}_h-\id)\big]\circ \bY^{m}_h \, |\p_{\widehat{\bm{\tau}}}\bY^{m}_h|\ge |\p_{\widehat{\bm{\tau}}}\bY^{m+1}_h|-|\p_{\widehat{\bm{\tau}}}\bY^{m}_h|,\quad \text{at} \ x_{j_i}=\bF^m_{h,j}(\qx_i).
\end{equation*}
Therefore, multiplying both sides by $|\nabla A_j|\,\qw_i$, taking the sum first over $i=1,\ldots,Q$, and then over $j=1,\ldots,N$, wo obtain on recalling \eqref{quadrature rule 1d} that
\begin{multline*}
	\begin{aligned}
 	\l(\nabla_s\bX^{m+1}_h,\nabla_s\l(\bX^{m+1}_h-\id\r)\r)^h_{\Gamma^m_h}
	&=\int_{\widehat{\Gamma}_h}^h\big[\nabla_s\bX^{m+1}_h\cdot \nabla_s(\bX^{m+1}_h-\id)\big]\circ \bY^{m}_h \, |\p_{\widehat{\bm{\tau}}}\bY^{m}_h| \ \d s\\
	&\geq\int_{\widehat{\Gamma}_h}^h|\p_{\widehat{\bm{\tau}}}\bY^{m+1}_h|\ \d s-\int_{\widehat{\Gamma}_h}^h|\p_{\widehat{\bm{\tau}}}\bY^{m}_h|\ \d s=L^{m+1}_h-L^m_h.
	\end{aligned}
\end{multline*}
The remainder of the proof proceeds as the proof of \Cref{Thm: energy stability for curve}.\end{proof}

\subsection{Structure-preserving isoparametric finite element scheme for surface diffusion of curves}\label{sec:SD2d:SP}

In this subsection, we propose an implicit isoparametric scheme for surface diffusion of curves that conserves the enclosed area exactly, and which can be seen as an extension of the work of Bao and Zhao \cite{Bao-Zhao} to higher-order. Here, it is possible to consider schemes with both exact and numerical integration. In fact, for the area enclosed by $\Gamma_h^m$ it holds that
\begin{align}
	A_h^m& \colonequals \frac{1}{2}\int_{\Gamma_h^m}\id \cdot \bn_h^m \ \d s 
	=\frac{1}{2}\int_{\widehat{\Gamma}_h}\bY_h^m\cdot (\p_{\widehat{\bm{\tau}}} \bY_h^m)^\perp \ \d s 
	=\frac{1}{2}\sum_{j=1}^N \int_{\widehat{\K}} \bF^m_{h,j} \cdot (\nabla \bF^m_{h,j})^\perp \ \d\widehat{x},\label{Definition of area}
\end{align}
where we notice that $\bF^m_{h,j} \cdot (\nabla \bF^m_{h,j})^\perp$ is a polynomial of degree $\le 2\ell-1$. Hence if the quadrature rule \eqref{quadrature rule 1d} is exact to degree $2\ell-1$, then \eqref{Definition of area} implies that
$A_h^m = \frac12 \l( \id , \bn_h^m \r)^h_{\Gamma^m_h}$. For simplicity, we only consider the scheme with numerical integration for the structure-preserving approximation in this subsection.

\smallskip

\noindent \textbf{Structure-preserving isoparametric finite element method with numerical integration for surface diffusionin 2d:} For $m\ge 0$, we find  $\bX_h^{m+1}\in [\V^m_\ell]^2$ and $\kappa^{m+1}_h\in \V^m_\ell$ such that the following equations hold for any $\varphi_h\in \V^m_\ell$ and $\bm{\eta}_h\in [\V^m_\ell]^2$
\begin{subequations}\label{SP:curve}
	\begin{align}
		\l(\frac{\bX^{m+1}_h-\id}{\tau},\varphi_h\bn^{m+\frac{1}{2}}_h \r)^h_{\Gamma^m_h}-\l(\nabla_s\kappa^{m+1}_h,\nabla_s\varphi_h\r)^h_{\Gamma^m_h}&=0,\label{SDF1:SP} \\
			\l(\kappa^{m+1}_h\bn^{m+\frac{1}{2}}_h,\bm{\eta}_h \r)^h_{\Gamma^m_h} +\l(\nabla_s\bX^{m+1}_h,\nabla_s \bm{\eta}_h \r)_{\Gamma^m_h}^{h}&=0.\label{SDF2:SP}
	\end{align}
\end{subequations}
Here, the intermediate normal vector is piecewise defined by
\begin{equation*} 
	\bn^{m+\frac{1}{2}}_h \circ \bY_h^m\colonequals \frac{(\p_{\widehat{\bm{\tau}}} \bY_h^{m+1}+\p_{\widehat{\bm{\tau}}} \bY_h^{m})^\perp}{2|\p_{\widehat{\bm{\tau}}} \bY_h^m|},\qquad \text{on}\quad \widehat{\sigma}_j,
\end{equation*}
where $\bY^{m+1}_h \colonequals \bX^{m+1}_h\circ \bY^{m}_h$ as usual. Hence, in contrast to our previous schemes, the approximation \eqref{SP:curve} leads to a system of nonlinear equations at each time level.

\begin{thm}[Structure-preserving property with numerical integration]\label{Structure-preserving for curve}
Assume that $(\mathcal{A}1^h)$ holds, that the quadrature rule \eqref{quadrature rule 1d} has positive weights $\qw_i>0$ for $i=1,\ldots,Q$, and that it is exact to degree $2\ell-1$. Let $(\bX^{m+1}_h,\kappa^{m+1}_h)\in [\V^m_\ell]^2\times \V^m_\ell$ be a solution of \eqref{SP:curve}. Then for any $M\ge 0$, we have
	\begin{equation}
		L^{M+1}_h+\tau\sum_{m=0}^{M}\l(\nabla_s\kappa^{m+1}_h,\nabla_s\kappa^{m+1}_h\r)_{\Gamma^{m}_h}^h\le L^{0}_h,\qquad A^{M+1}_h=A^{M}_h=\cdots A^{1}_h= A^0_h,
	\end{equation}
	where $A^{m}_h$ is the area enclosed by $\Gamma^m_h$.	
\end{thm}

\begin{proof}
	The energy stability can be shown as in the proof of \Cref{Thm: energy stability for curve num}. 
For the area-preserving property, we first define the linear interpolation $\Gamma_h(\alpha)$ between the discrete curves $\Gamma^{m+1}_h=\bY^{m+1}_h(\widehat{\Gamma}_h)$ and $\Gamma^{m}_h=\bY^{m}_h(\widehat{\Gamma}_h)$ via $\Gamma_h(\alpha) =\bY_h(\widehat\Gamma_h,\alpha)$ for
\begin{equation*}
	\bY_h(\cdot,\alpha) \colonequals \alpha\bY^{m+1}_h+(1-\alpha)\bY^{m}_h,\quad 0\le \alpha\le 1.
\end{equation*}
Note that $\bY_h(\cdot,\alpha)$ is also a piecewise polynomial of degree $\le \ell$ and is uniquely determined by the values $(1-\alpha)\bY^{m+1}_h(x_{j_i})+\alpha\bY^{m}_h(x_{j_i})$ at the Lagrange nodes $\{x_{j_i}\}_{j=1,\ldots,N,i=1,\ldots,\ell+1}\subset \widehat{\Gamma}_h$. Then, for the family of curves $\Gamma_h(\alpha)$, we define the normal
\[
\bn_h(\alpha) \circ \bY_h(\alpha) \colonequals \frac{(\p_{\widehat{\bm{\tau}}}\bY_h(\alpha))^\perp}{|\p_{\widehat{\bm{\tau}}}\bY_h(\alpha)|},\qquad \text{in}\quad \widehat{\sigma}_j,\quad j=1,\ldots,N.
\]
Thanks to our assumptions on the quadrature rule,
the area $A_h(\alpha)$ enclosed by $\Gamma_h(\alpha)$ can be computed exactly by the natural analogue of the formula \eqref{Definition of area}.

Applying the Reynolds transport theorem \cite{Bansch} to $\Gamma_h(\alpha)$, integrating with respect to $\alpha$ from $0$ to $1$, and setting $\varphi_h=\tau$ in \eqref{SDF1:SP}, we have
\begin{align*}
	A^{m+1}_h-A^{m}_h
	&=\int^1_0\frac{\d }{\d \alpha}A_h(\alpha)\ \d\alpha\\
	&=\int^1_0\int_{\Gamma_h(\alpha)}\p_\alpha \bY_h(\alpha)\cdot   \bn_h(\alpha)\ \d s\ \d\alpha\\
	&=\int^1_0\int_{\widehat{\Gamma}_h}[\bY^{m+1}_h-\bY^m_h] \cdot (\p_{\widehat{\bm{\tau}}} \bY_h(\alpha))^{\perp} \ \d s\ \d\alpha\\
	&= \int_{\widehat{\Gamma}_h}[\bY^{m+1}_h-\bY^m_h] \cdot \int^1_0\l[\alpha\p_{\widehat{\bm{\tau}}}\bY^{m+1}_h+(1-\alpha)\p_{\widehat{\bm{\tau}}}\bY^{m}_h\r]^{\perp}\ \d \alpha\ \d s\\
	&=\int_{\widehat{\Gamma}_h}[\bY^{m+1}_h-\bY^m_h]\cdot \l[\frac{\p_{\widehat{\bm{\tau}}} \bY_h^{m+1}+\p_{\widehat{\bm{\tau}}} \bY_h^{m}}{2} \r]^\perp \ \d s \\
    &=\int_{\widehat{\Gamma}_h}^h[\bY^{m+1}_h-\bY^m_h]\cdot \l[\frac{\p_{\widehat{\bm{\tau}}} \bY_h^{m+1}+\p_{\widehat{\bm{\tau}}} \bY_h^{m}}{2} \r]^\perp \ \d s \\
	&=
    \big(\bX^{m+1}_h-\id,\bn^{m+\frac{1}{2}}_h\big)_{\Gamma^{m}_h}^h=\big(\nabla_s\kappa^{m+1}_h,\nabla_s\tau\big)_{\Gamma^{m}_h}^h=0,
\end{align*}
where we have observed that the integrand $[\bY^{m+1}_h-\bY^{m}_h] \cdot \l[\frac{\p_{\widehat{\bm{\tau}}} \bY_h^{m+1}+\p_{\widehat{\bm{\tau}}} \bY_h^{m}}{2} \r]^\perp$
is a polynomial of degree $2\ell-1$ on $\widehat{\sigma}_j$. We then complete the proof of the area-preservation property by induction.\end{proof}

\section{Surface evolution}

In this section, we consider mean curvature flow and surface diffusion for a family of surfaces $\Gamma(t)$, determined by
\begin{align*}
	\mathcal{V}=\begin{cases}
		H,\qquad &\text{mean curvature flow of surfaces},\\
		-\Delta_{\Gamma}H ,\qquad\ &\text{surface diffusion of surfaces},
	\end{cases}
\end{align*}
where $\mathcal{V}$ is the normal velocity along the surface $\Gamma$, $H$ is the mean curvature and $\Delta_{\Gamma}$ represents the Laplacian--Beltrami operator on $\Gamma$.

\subsection{Isoparametric finite element approximation}

In this section, we use the surface finite element method setup proposed in \cite{Demlow,Elliott,Heine04}. We first define the reference element space $(\widehat{\K},\widehat{\V},\widehat{\Sigma})$ in $\R^2$, where $\widehat{\K}$ is the classical $2$-simplex in $\R^2$, $\widehat{\V}=\P^{\ell}(\widehat{\K})$ and $\widehat{\Sigma}$ is the classical set of Lagrange nodes \cite{Brenner-Scott}. Similar to the case of the curves, we define a reference grid $\widehat{\Gamma}_h$ consisting of flat triangles $\widehat{\sigma}_j$ as the image of the reference simplex by affine linear maps $A_j$, and set $h=\max_{j=1}^J\mathrm{diam}(\widehat{\sigma}_j)$, where $\mathrm{diam}(\widehat{\sigma}_j)$ is the diameter of the triangle $\widehat{\sigma}_j$. Subsequently, we find the parametrization $\bY^m_h$ over the reference grid by the isoparametric finite element method; equivalently, we find the piecewise approximation $\bF^m_{h,j}$ over the reference simplex.

For the full discretization, suppose for $\Gamma(t_m)$ we already have an approximating surface
\[\Gamma^m_h=\bigcup_{j=1}^J\sigma^m_j =\bY^m_h(\widehat{\Gamma}_h)\]
for some function $\l.\bY^m_h\r|_{\widehat{\sigma}_j}\in \P^{\ell}(\widehat{\sigma}_j)$ and $m\ge 0$, we can similarly define the evolving surface finite element space $\W^m_\ell$ over $\Gamma^m_h$ as
\begin{align*}
	\W_\ell^m
	& \colonequals \l\{v_h\in H^1(\Gamma_h^m)\,:\,v_h\circ  \l.\bY_h^m\r|_{\widehat{\sigma}_j}\in \P^{\ell}(\widehat{\sigma}_j), \; j=1,\ldots,J.\r\}.
\end{align*}
We consider the following two formulations for geometric evolution equations of surfaces \cite{BGN07A,BGN07B}:  mean curvature flow (MCF) and  surface diffusion (SD)
	\begin{align*}
		\p_t\bY\circ \bY^{-1}\cdot \bn&=\begin{cases}
			H,\quad &\text{MCF in 3d},\\
			-\Delta_{\Gamma}H ,\quad &\text{SD in 3d},
		\end{cases}\\
		H\bn&=\Delta_{\Gamma}\id ,
	\end{align*}
where $\bY(t)\colon\Gamma(0) \to \Gamma(t)$ is the flow map and can be seen as the parametrization of $\Gamma(t)$ over $\Gamma(0)$. Here the coupled equations are defined on the evolving surface $\Gamma(t)$ and $\bn$ denotes the unit outer normal to $\Gamma(t)$. We then consider the following fully discrete higher-order finite element method:
\smallskip

\noindent \textbf{Isoparametric finite element method with exact integration:} Find the solution (called the local flow map) $\bX_h^{m+1}\in [\W_\ell^m]^3$ and $H^{m+1}_h\in \W_\ell^m$ that satisfy the following equations
\begin{subequations}\label{eq:3d}
	\begin{align}
		\l(\frac{\bX^{m+1}_h-\id}{\tau},\varphi_h\bn^{m}_h \r)_{\Gamma^m_h}&=\begin{cases}
			\l(H^{m+1}_h,\varphi_h\r)_{\Gamma^m_h},\quad &\text{MCF in 3d},\\
			\l(\nabla_{\Gamma^m_h}H^{m+1}_h,\nabla_{\Gamma^m_h}\varphi_h\r)_{\Gamma^m_h},\quad &\text{SD in 3d},
		\end{cases} \label{eq:3d1} \\
			\l(H^{m+1}_h\bn^m_h,\bm{\eta}_h \r)_{\Gamma^m_h} &+\l(\nabla_{\Gamma^m_h} \bX^{m+1}_h,\nabla_{\Gamma^m_h} \bm{\eta}_h \r)_{\Gamma^m_h}=0\label{eq:3d2}
	\end{align}
\end{subequations}
for any $\varphi_h\in \W_\ell^m$ and $\bm{\eta}_h\in [\W_\ell^m]^3$. 
Here $\bn^m_h$ and $\nabla_{\Gamma^m_h}$ are the piecewise defined unit normal and surface gradient on $\Gamma^m_h$.
We then define $\bY^{m+1}_h \colonequals \bX^{m+1}_h\circ \bY^{m}_h$, and $\Gamma^{m+1}_h \colonequals \bY^{m+1}_h(\widehat{\Gamma}_h)$ is the desired isoparametric higher-order approximation of $\Gamma(t_{m+1})$.

As before, the inner product $(\cdot,\cdot)_{\Gamma^m_h}$ is the exact integration defined as
\begin{equation} \label{exact integration:3d}
	(u,v)_{\Gamma^m_h}
	  \colonequals \int_{\Gamma^{m}_h} u\cdot v\ \d s=\int_{\widehat{\Gamma}_h}(u\circ \bY^m_h)\cdot (v\circ \bY^m_h)\, |\mathcal{J}(\bY^{m}_h)| \ \d s,
\end{equation}
where $\cJ(\bY^m_h)=\p_{\widehat{\bm{\tau}}_1}\bY^m_h\times \p_{\widehat{\bm{\tau}}_2}\bY^m_h$ on $\widehat{\sigma}_j$,
 with $\{\widehat{\bm{\tau}}_1,\widehat{\bm{\tau}}_2\}$ being an orthonormal basis of the tangent space of $\widehat{\sigma}_j$.
Similarly, the unit normal vector is piecewise represented as
\begin{equation} \label{eq:nmh}
    \bn^m_h \circ \bY^m_h =\frac{ \cJ(\bY^m_h)}{|\cJ(\bY^m_h)|} \quad \text{on}\quad \widehat{\sigma}_j.
\end{equation}
Here we assume that the
orientations of $\widehat{\sigma}_j$ and $\bY^m_h$ are chosen such that $n^m_h$
almost everywhere on $\Gamma^m_h$ is the outer unit normal.

Similar to the curve case, in practice, we choose a suitable quadrature rule on the reference simplex $\widehat{\K}$ to compute the integrals. Let $\{\qbx_{i}\}_{i=1}^Q\subset\widehat{\mathbb{K}}$ be the distinct quadrature points and let $\{\qw_{i}\}_{i=1}^Q$ be the corresponding quadrature weights. We introduce
\begin{equation}\label{quadrature rule 2d}
	\int_{\widehat{\K}}^hf(\bxHat)\ \d \bxHat \colonequals  \sum_{i=1}^{Q}f(\qbx_{i})\ \qw_{i}.
\end{equation}
The numerical integration over each flat triangle $\widehat{\sigma}_j$ of the reference grid $\widehat{\Gamma}_h$ is then defined by
\begin{equation*}
	\int^h_{\widehat{\sigma}_j}f(\bm{x})\ \d s \colonequals \sum_{i=1}^Q (f\circ A_j)(\qbx_i)\, |\cJ(A_j)|(\qbx_i) \ \qw_{i}.
\end{equation*}
where $\cJ(A_j)=\p_{\bm e_1} A_j \times \p_{\bm e_2} A_j$, with $\bm e_1$ and $\bm e_2$ denoting the unit vectors in $\mathbb R^2$.
We say that the quadrature rule \eqref{quadrature rule 2d} is exact to degree $p$ if, for any $f\in \P^{p}(\widehat{\K})$, the numerical integral is equal to the integral of $f$. We define the quadrature-based inner product $(u,v)_{\Gamma^m_h}^h$ over the surface $\Gamma^m_h$ as
\begin{align}
    (u,v)^h_{\Gamma^m_h}
    & \colonequals \int_{\widehat{\Gamma}_h}^h(u\circ \bY^m_h)\cdot (v\circ \bY^m_h)\, |\mathcal{J}(\bY^{m}_h)| \ \d s \notag\\
    &\colonequals \sum_{j=1}^J\int_{\widehat{\sigma}_j}^h(u\circ \bY^m_h)\cdot (v\circ \bY^m_h)\, |\mathcal{J}(\bY^{m}_h)| \ \d s\notag \\
    &= \sum_{j=1}^J\int_{\widehat{\K}}^h(u\circ \bF^m_{h,j})\cdot (v\circ \bF^m_{h,j})\, |\mathcal{J}(\bF^m_{h,j})| \ \d\widehat{\bm{x}}\notag\\
    &=\sum_{j=1}^J\sum_{i=1}^Q(u\circ \bF^m_{h,j})(\qbx_i)\cdot (v\circ \bF^m_{h,j})(\qbx_i)|\cJ(\bF^m_{h,j})|(\qbx_i)\ \qw_i
    \label{Numerical integration 3d}
\end{align}
for functions $u,v$ defined on the discrete surface $\Gamma^m_h$, where $\cJ(\bF^m_{h,j})=\p_{\bm e_1}\bF_{h,j}^m\times \p_{\bm e_2}\bF_{h,j}^m$.
In the above we have used that $|\cJ(\bF^{m}_{h,j})|=|\cJ(\bY^{m}_{h})||\cJ(A_j)|$, which follows from \Cref{lemma:appendix} in \Cref{sec:appendix:transform}.
The following  higher-order method is what we consider in practice:

\noindent \textbf{Isoparametric finite element method with numerical integration:} Find the solution (called the local flow map) $\bX_h^{m+1}\in [\W_\ell^m]^3$ and $H^{m+1}_h\in \W_\ell^m$ that satisfy the following equations
\begin{subequations}\label{eq:3dnum}
	\begin{align}
		\l(\frac{\bX^{m+1}_h-\id}{\tau},\varphi_h\bn^{m}_h \r)^h_{\Gamma^m_h}&=\begin{cases}
			\l(H^{m+1}_h,\varphi_h\r)^h_{\Gamma^m_h},\quad &\text{MCF in 3d},\\
			\l(\nabla_{\Gamma^m_h}H^{m+1}_h,\nabla_{\Gamma^m_h}\varphi_h\r)^h_{\Gamma^m_h},\quad &\text{SD in 3d},
		\end{cases} \label{eq:3d1num} \\
			\l(H^{m+1}_h\bn^m_h,\bm{\eta}_h \r)^h_{\Gamma^m_h} &+\l(\nabla_{\Gamma^m_h} \bX^{m+1}_h,\nabla_{\Gamma^m_h} \bm{\eta}_h \r)_{\Gamma^m_h}^h=0,\label{eq:3d2num}
	\end{align}
\end{subequations}
for any $\varphi_h\in \W_\ell^m$ and $\bm{\eta}_h\in [\W_\ell^m]^3$. 

\subsection{Well-posedness and unconditional energy stability}

Similar to the curve evolution case, we prove the well-posedness under some mild conditions.

\begin{thm}[Well-posedness for \eqref{eq:3d} and \eqref{eq:3dnum}]
Assume that the discrete surface $\Gamma^m_h$ satisfies
\begin{enumerate}
	\item[$(\mathcal{B}1)$] The surface elements are nondegenerate,   i.e., for each $\widehat{\sigma}_j$, $j=1,\ldots,J$, it holds
	\begin{equation*}
	|\mathcal{J}(\bY^m_h)|> 0,\quad \text{on}\,\,\, \widehat{\sigma}_j.
\end{equation*}
\item[$(\mathcal{B}2)$] Nonparallel assumption, that is
\begin{equation*} 
		\mathrm{dim}\ \mathrm{span} \l\{\l(\bn^m_h\varphi_h,1 \r)_{\Gamma^m_h}| \, \varphi_h\in \W_\ell^m \r\}= 3.
\end{equation*}
\end{enumerate}
Then the system \eqref{eq:3d} has a unique solution. Moreover, assume that the discrete surface $\Gamma^m_h$ satisfies
\begin{enumerate}
   \item[$(\mathcal{B}1^h)$] The surface elements are nondegenerate over the quadrature points,  i.e., for each $\widehat{\sigma}_j$, $j=1,\ldots,J$, it holds
	\begin{equation*}
	|\mathcal{J}(\bY^m_h)|(\bx_{j_i})> 0,\quad  i=1,\ldots,Q,
\end{equation*}
where $\{\bx_{j_i}\}_{i=1}^{Q}\subset \widehat{\sigma}_j$ is the image of $\{\qbx_i\}_{i=1}^{Q}$ under the affine linear map $A_j$. Equivalently,
\begin{equation*}
   |\nabla_{\widehat{\bx}}\bF^m_{h,j}(\qbx_i)|>0,\quad j=1,\ldots,N,\ i=1,\ldots,Q.
\end{equation*}
\item[$(\mathcal{B}2^h)$] Nonparallel assumption, that is
\begin{equation*} 
		\mathrm{dim}\ \mathrm{span} \l\{\l(\bn^m_h\varphi_h,1 \r)^h_{\Gamma^m_h}| \, \varphi_h\in \W_\ell^m \r\}= 3.
\end{equation*}
\end{enumerate}
And we further assume that the quadrature rule \eqref{quadrature rule 2d} satisfies 
\begin{itemize}
     \item[$(\mathcal{B}3^h)$] It holds that $\qw_i>0$ for $i=1,\ldots,Q$. Moreover, if $f\in \P^{\ell}(\widehat{\K})$ with $f(\qbx_i)=0$ for $i=1,\ldots,Q$,
     then $f=0$.
\end{itemize}
Then the system \eqref{eq:3dnum} has a unique solution.\end{thm}

\begin{proof}
	The argument is similar to the proofs of \Cref{Well-posed:2d:exact,Well-posed:2d:num}, and so we omit it here.~\end{proof}

\begin{remark} \label{rem:assB}
	Similar to \Cref{Well-posed:remark},  the assumption $(\mathcal{B}1)$ implies  $(\mathcal{B}1^h)$. Moreover, by writing the normal vector over the reference simplex as
	\[
   \bn^m_h \circ \bF_{h,j}^m =\frac{\cJ(\bF_{h,j}^m)}{|\cJ(\bF_{h,j}^m)|},
	\]
we see that the assumption $(\mathcal{B}2^h)$ is equivalent to $(\mathcal{B}2)$ if the quadrature rule \eqref{quadrature rule 2d} is exact to degree $3\ell-2$, recall \eqref{Numerical integration 3d}.
\end{remark}

The following lemma is the natural extension of \Cref{lemma:curve} to surfaces. Recall also \cite[Lemma~57]{BGN20}.

\begin{lemma}\label{lemma:surface}
	Let $\widehat{\sigma},\sigma,\widetilde{\sigma}\subset \R^3$ be three $C^1$-surfaces, possibly with boundary. Let $\sigma$ and $\widetilde{\sigma}$ be parametrized by $\bY\colon\widehat{\sigma} \to \R^3$ and $\widetilde{\bY}=\bX\circ \bY\colon\widehat{\sigma} \to \R^3$, respectively, where $\bX\in C^1(\sigma,\R^3)$. The following inequality holds everywhere in $\widehat{\sigma}$
	\begin{equation} \label{eq:lemmasurface}
		\big[\nabla_{\sigma} \bX: \nabla_{\sigma}(\bX-\id)\big]\circ \bY \, |\mathcal{J}(\bY)| \ge |\mathcal{J}(\widetilde{\bY})| -|\mathcal{J}(\bY)|,
	\end{equation}
with equality if $\bX=\id\in C^1(\sigma,\R^3)$. Here, $\nabla_{\sigma}$ is the surface gradient of $\sigma$, and  $|\mathcal{J}(\widetilde{\bY})|$, $|\mathcal{J}(\bY)|$ are the surface elements of $\widetilde{\sigma}$ and  $\sigma$, respectively.

\end{lemma}

\begin{proof}
	Let  $\{\bm{\tau}_1,\bm{\tau}_2\}$ be a local orthonormal basis of the tangent space of the surface $\sigma$. With this basis, we can write
\[
\nabla_{\sigma} \bX= \sum_{k=1}^2\p_{\bm{\tau}_k}\bX\otimes \bm{\tau}_k,\
|\nabla_{\sigma} \bX|^2= \sum_{k=1}^2 |\p_{\bm{\tau}_k}\bX|^2,\quad 
\text{and}\quad \nabla_{\sigma} \id=\sum_{k=1}^2\bm{\tau}_k \otimes \bm{\tau}_k,\
|\nabla_{\sigma} \id|^2 = 2,
\]
where $\p_{\bm{\tau}_k}$ is the directional derivative with respect to $\bm{\tau}_k$. The change of variable gives the equality
\begin{equation}\label{change of variable}
	|\mathcal{J}(\widetilde{\bY})|= \sqrt{\det\l(\p_{\bm{\tau}_i}\bX\cdot \p_{\bm{\tau}_j}\bX \r)_{1\le i,j\le 2} }\circ \bY\, |\mathcal{J}(\bY)|,
\end{equation}
which follows from \Cref{lemma:appendix} in \Cref{sec:appendix:transform}.
Combining \eqref{change of variable} and the inequality
\[
\sqrt{\det\l(\p_{\bm{\tau}_i}\bX\cdot \p_{\bm{\tau}_j}\bX \r)_{1\le i,j\le 2} }=\sqrt{|\p_{\bm{\tau}_1}\bX|^2|\p_{\bm{\tau}_2}\bX|^2-(\p_{\bm{\tau}_1}\bX\cdot \p_{\bm{\tau}_2}\bX)^2}\le |\p_{\bm{\tau}_1}\bX||\p_{\bm{\tau}_2}\bX|,
\]
we can then estimate
\begin{align*}
	\big[\nabla_{\sigma} \bX: \nabla_{\sigma}(\bX-\id)\big]\circ \bY \, |\mathcal{J}(\bY)|&\ge \frac{1}{2}\l(|\nabla_{\sigma} \bX|^2-|\nabla_{\sigma}\id|^2 \r)\circ \bY \, |\mathcal{J}(\bY)|\\
	&= \frac{1}{2}\l(|\p_{\bm{\tau}_1}\bX|^2+|\p_{\bm{\tau}_2}\bX|^2\r)\circ \bY \, |\mathcal{J}(\bY)|-|\mathcal{J}(\bY)|\\
	&\ge |\p_{\bm{\tau}_1}\bX||\p_{\bm{\tau}_2}\bX|\circ\bY\ |\mathcal{J}(\bY)|-|\mathcal{J}(\bY)|\\
	&= \sqrt{\det\l(\p_{\bm{\tau}_i}\bX\cdot \p_{\bm{\tau}_j}\bX \r)_{1\le i,j\le 2} }\circ\bY\ |\mathcal{J}(\bY)|-|\mathcal{J}(\bY)|\\
	&=|\mathcal{J}(\widetilde{\bY})|-|\mathcal{J}(\bY)|,
\end{align*}
where we have used the basic inequalities
\[
\bA:(\bA-\bB)\ge \frac{1}{2}\l(|\bA|^2-|\bB|^2\r),\quad\text{and}\quad  a^2+b^2\ge 2ab.
\]
This completes the proof, on noting that \eqref{eq:lemmasurface} holds with equality if $\bX=\id\in C^1(\sigma,\R^3)$.\end{proof}

\begin{remark}
  We observe that the analogue of the estimate we proved in Lemma~\ref{lemma:curve} for curves in $\R^2$, and in Lemma~\ref{lemma:surface} for surfaces in $\R^3$, no longer holds true for hypersurfaces in $\R^d$, $d\ge 4$. Indeed, by a scaling argument we see that there cannot be a constant $c_0>0$ such that 
    \begin{equation}
        c_0\sum_{i=1}^{d-1} |\p_{\bm{\tau}_i}\bX|^2 \ge \sqrt{\det\l(\p_{\bm{\tau}_i}\bX\cdot \p_{\bm{\tau}_j}\bX \r)_{1\le i,j\le d-1} },\qquad d\ge 4,
    \end{equation}
    since the left hand side scales quadratically in $\bX$, while the right hand side scales to the power $d-1$.
    See also \cite[Remark~56]{BGN20}.
\end{remark}

\begin{thm}[Unconditional energy stability with exact integration]\label{Thm: energy stability for surface}
Under assumption $(\mathcal{B}1)$, suppose  $(\bX^{m+1}_h,\kappa^{m+1}_h)$ $\in [\W^{m}_\ell]^3\times \W^{m}_\ell$ is a solution of \eqref{eq:3d}. Then for any $M\ge 0$, we have
	\begin{equation*}
		\cH^2(\Gamma^{M+1}_h) \le \cH^2(\Gamma^{0}_h)-\begin{cases}
			\tau\sum_{m=0}^{M}\l(H^{m+1}_h,H^{m+1}_h\r)_{\Gamma^{m}_h},\quad &\text{MCF},\\
			\tau\sum_{m=0}^{M}\l(\nabla_{\Gamma^m_h}H^{m+1}_h,\nabla_{\Gamma^m_h}H^{m+1}_h\r)_{\Gamma^{m}_h},\quad &\text{SD},
		\end{cases}
	\end{equation*}
	where $\cH^2(\Gamma^{m}_h) \colonequals \int_{\Gamma^{m}_h}\d s=\int_{\widehat{\Gamma}_h}|\mathcal{J}(\bY^m_h)|\ \d s$ is the exact energy (surface area) of the discrete surface $\Gamma^m_h$.
\end{thm}

\begin{proof}
The proof proceeds as the proof of \Cref{Thm: energy stability for curve}, by taking $\varphi_h=H^{m+1}_h\in \W^{m}_\ell$ in \eqref{eq:3d1} and $\bm{\eta}_h=\bX^{m+1}_h-\id_{|_{\Gamma^m_h}}\in [\W^{m}_\ell]^3$ in \eqref{eq:3d2}, and using  \Cref{lemma:surface}.
\end{proof}

\begin{thm}[Unconditional energy stability with numerical integration]
Assume that $(\mathcal{B}1^h)$ holds and that the quadrature rule \eqref{quadrature rule 2d} has positive weights $\qw_i>0$ for $i=1,\ldots,Q$. Let  $(\bX^{m+1}_h,\kappa^{m+1}_h)$ $\in [\W^{m}_\ell]^3\times \W^{m}_\ell$ be a the solution of \eqref{eq:3dnum}. 
Then for any $M\ge 0$, we have
	\begin{equation*} 
		S^{M+1}_h \le S_h^{0}-\begin{cases}
			\tau\sum_{m=0}^{M}\l(H^{m+1}_h,H^{m+1}_h\r)^h_{\Gamma^{m}_h},\quad &\text{MCF},\\
			\tau\sum_{m=0}^{M}\l(\nabla_{\Gamma^m_h}H^{m+1}_h,\nabla_{\Gamma^m_h}H^{m+1}_h\r)^h_{\Gamma^{m}_h},\quad &\text{SD},
		\end{cases}
	\end{equation*}
where $S_h^{m} \colonequals \l(1,1\r)^h_{\Gamma^m_h}$ is the numerical energy (surface area) of the discrete surface $\Gamma^m_h$.

\end{thm}

\begin{proof}
	The proof proceeds as the proof of \Cref{Thm: energy stability for surface}, where for the application of the quadrature rule with positive weights
we argue as in the proof of \Cref{Thm: energy stability for curve num}.
\end{proof}

\begin{remark}
	It is worth mentioning another isoparametric finite element method for the simulation of mean curvature flow. The method was first introduced for linear finite elements by Dziuk \cite{Dziuk1991} and further developed for finite elements of higher spatial order by Li \cite{Li1,Li2}. By utilizing the classical formula $H\bn=\Delta_{\Gamma}\id$, Dziuk's formulation for surfaces can be summarized to find the local flow map  $\bX^{m+1}_h$ such that
	\begin{equation}\label{Dziuk-type}
		\l(\frac{\bX^{m+1}_h-\id}{\tau},\bm{\eta}_h \r)_{\Gamma^{m}_h}+\l(\nabla_{\Gamma^m_h} \bX^{m+1}_h,\nabla_{\Gamma^m_h} \bm{\eta}_h \r)_{\Gamma^m_h}=0, \quad \bm{\eta}_h\in [\W_\ell^m]^3.
	\end{equation}
	In practice, a suitable numerical integration 
    is necessary when implementing \eqref{Dziuk-type} for $\P^{\ell}$-elements with $\ell\ge 2$. Li proved convergence results for \eqref{Dziuk-type} with order $\ell\ge 3$ for the curve case and with order $\ell\ge 6$ for the surface case. Very recently, Bai and Li \cite{Bai2023} proved a convergence result under the projected distance with order $\ell\ge 3$.

	We point out that the proof of unconditional energy stability in this paper holds for the above Dziuk-type scheme \eqref{Dziuk-type}, as well as its numerical integration form
    \begin{equation}\label{Dziuk-typeh}
		\l(\frac{\bX^{m+1}_h-\id}{\tau},\bm{\eta}_h \r)_{\Gamma^{m}_h}^h+\l(\nabla_{\Gamma^m_h} \bX^{m+1}_h,\nabla_{\Gamma^m_h} \bm{\eta}_h \r)_{\Gamma^m_h}^h=0 , \quad \bm{\eta}_h\in [\W_\ell^m]^3.
	\end{equation}
    Indeed, by taking $\bm{\eta}_h=\bX^{m+1}_h-\id_{|_{\Gamma^m_h}}$ in \eqref{Dziuk-type} and applying \Cref{lemma:surface}, we can derive
	\begin{align*}
		\cH^2(\Gamma^{m+1}_h)-\cH^2(\Gamma^{m}_h)
		&\le \l(\nabla_{\Gamma^m_h}\bX^{m+1}_h,\nabla_{\Gamma^m_h}(\bX^{m+1}_h-\id)\r)_{\Gamma^m_h}\\
		&=- \frac{1}{\tau}\l(\bX^{m+1}_h-\id,\bX^{m+1}_h-\id\r)_{\Gamma^{m}_h}\le 0,
	\end{align*}
    and similarly for \eqref{Dziuk-typeh}. There is no direct analogue of the Dziuk method for surface diffusion. Moreover, the schemes \eqref{Dziuk-type} and \eqref{Dziuk-typeh} may suffer from mesh distortions, especially near singularities. We refer to \Cref{sec:mcf} for a comparison of the Dziuk-type schemes and the BGN-type schemes proposed in this paper.
\end{remark}

\subsection{A structure-preserving isoparametric finite element  method for surface diffusion of surfaces}

As an analogue of the curve case, we consider the following volume-conserving scheme for surface diffusion in three dimensions:

\smallskip

\noindent\textbf{Structure-preserving isoparametric finite element method with numerical integration for surface diffusion in 3d:} For $m\ge 0$, find   $\bX_h^{m+1}\in [\W^{m}_\ell]^3$ and $H^{m+1}_h\in \W^{m}_\ell$ that satisfy the following equations for any $\varphi_h\in \W^{m}_\ell$ and $\bm{\eta}_h\in [\W^{m}_\ell]^3$
\begin{subequations}\label{SP:surface}
	\begin{align}
		\l(\frac{\bX^{m+1}_h-\id}{\tau},\varphi_h\bn^{m+\frac{1}{2}}_h \r)^h_{\Gamma^m_h}-\l(\nabla_{\Gamma^m_h}H^{m+1}_h,\nabla_{\Gamma^m_h}\varphi_h\r)^h_{\Gamma^m_h}&=0,\label{SDF3d1:SP} \\
			\l(H^{m+1}_h\bn^{m+\frac{1}{2}}_h,\bm{\eta}_h \r)^h_{\Gamma^m_h} +\l(\nabla_{\Gamma^m_h}\bX^{m+1}_h,\nabla_{\Gamma^m_h} \bm{\eta}_h \r)_{\Gamma^m_h}^{h}&=0.\label{SDF3d2:SP}
	\end{align}
\end{subequations}
Here, the intermediate normal vector is piecewise defined by
\begin{equation*}
	\bn^{m+\frac{1}{2}}_h \circ\bY^m_h\colonequals \frac{\cJ(\bY^m_h)+4\cJ(\bY^{m+1/2}_h)+\cJ(\bY^{m+1}_h)}{6|\cJ(\bY^m_h)|},\qquad \text{on } \widehat{\sigma}_j,
\end{equation*}
where $\bY^{m+1/2}_h$ is piecewise parametrized over $\widehat{\Gamma}_h$ by the average value of $\bY^{m}_h$ and $\bY^{m+1}_h$ at the Lagrange nodes. Note that the definition the volume enclosed by $\Gamma^m_h$ can also be exactly computed by taking the quadrature rule \eqref{quadrature rule 2d} exact to degree $3\ell-2$ as follows
\begin{align}
	V_h^m& \colonequals \frac{1}{3}\int_{\Gamma_h^m}\id\cdot \bn_h^m \ \d s=\frac{1}{3}\int_{\widehat{\Gamma}_h}\bY^m_h\cdot \mathcal{J}(\bY^m_h) \ \d s=\frac{1}{3}\sum_{j=1}^J\int_{\widehat{\K}}\bF^m_{h,j}\cdot  \cJ(\bF^m_{h,j}) \ \d\bxHat=\frac{1}{3}(\id,\bn_h^m)^h_{\Gamma^m_h} ,\label{Definition of volume}
\end{align}
where we have recalled \eqref{exact integration:3d}, \eqref{eq:nmh}, \eqref{Numerical integration 3d} and \Cref{rem:assB}.

We now prove the following structure-preserving property.
\begin{thm}[Structure-preserving property with numerical integration]\label{Structure-preserving for surface}
	Assume that $(\mathcal{B}1^h)$ holds, that the quadrature rule \eqref{quadrature rule 2d} has positive weights $\qw_i>0$ for $i=1,\ldots,Q$, and that it is exact to degree $3\ell-2$. Let $(\bX^{m+1}_h,H^{m+1}_h)\in [\W^{m}_\ell]^3\times \W^{m}_\ell$ be a solution of \eqref{SP:surface}. Then for any $M\ge 0$, we have
	\begin{equation}
		S^{M+1}_h+\tau\sum_{m=0}^{M}\l(\nabla_{\Gamma^{m}_h} H^{m+1}_h,\nabla_{\Gamma^{m}_h}H^{m+1}_h\r)_{\Gamma^{m}_h}^h\le S^{0}_h,\qquad V^{M+1}_h=V^{M}_h=\cdots V^{1}_h= V^0_h,
	\end{equation}
	where $V^{m}_h$ is the volume enclosed by $\Gamma^m_h$.	
\end{thm}

\begin{proof}
	The energy stability can be proved similarly to before, and we omit it here for brevity. For the volume-preserving property, we first define the linear interpolation $\Gamma_h(\alpha)$ between the dsicrete surfaces $\Gamma^{m+1}_h=\bY^{m+1}_h(\widehat{\Gamma}_h)$ and $\Gamma^{m}_h=\bY^{m}_h(\widehat{\Gamma}_h)$ via $\Gamma_h(\alpha) = \bY_h(\widehat{\Gamma}_h,\alpha)$, where
\begin{equation}
	\bY_h(\cdot,\alpha) \colonequals \alpha\bY^{m+1}_h+(1-\alpha)\bY^{m}_h,\quad 0\le \alpha\le 1.
\end{equation}
Note that $\bY_h(\cdot,\alpha)$ is uniquely determined by the values $\alpha\bY^{m+1}_h(\bx_{j_i})+(1-\alpha)\bY^{m}_h(\bx_{j_i})$ at the Lagrange nodes $\{\bx_{j_i}\}\subset \widehat{\Gamma}_h$. For $\Gamma_h(\alpha)$ we define the normal vector
\[
\bn_h(\alpha) \circ \bY_h(\alpha)\colonequals \frac{\cJ(\bY_h(\alpha))}{|\cJ(\bY_h(\alpha))|},\qquad \text{on } \widehat{\sigma}_j,
\]
and the volume $V_h(\alpha)$ enclosed by $\Gamma_h(\alpha)$ can be computed as in \eqref{Definition of volume}.

Applying the Reynolds transport theorem \cite{Bansch} to $\Gamma_h(\alpha)$, integrating with respect to $\alpha$ from $0$ to $1$ and taking $\varphi_h=\tau$ in \eqref{SDF3d1:SP}, we have
\begin{multline*}
	\begin{aligned}
	V^{m+1}_h-V^{m}_h
	&=\int^1_0\frac{\d }{\d \alpha}V_h(\alpha)\ \d\alpha\\
	&=\int^1_0\int_{\Gamma_h(\alpha)}\p_\alpha \bY_h(\alpha)\cdot   \bn_h(\alpha)\ \d s\ \d\alpha\\
	&=\int^1_0\int_{\widehat{\Gamma}_h}[\bY^{m+1}_h-\bY^m_h] \cdot \cJ(\bY_h(\alpha)) \ \d s\ \d\alpha\\
	&=\int_{\widehat{\Gamma}_h}[\bY^{m+1}_h-\bY^m_h]\cdot \int^1_0\cJ(\bY_h(\alpha))\ \d \alpha\ \d s\\
    &=\int_{\widehat{\Gamma}_h}[\bY^{m+1}_h-\bY^m_h]\cdot \int^1_0\l[\p_{\widehat{\bm{\tau}}_1}\l(\alpha\bY^{m+1}_h+(1-\alpha)\bY^{m}_h\r)\times \p_{\widehat{\bm{\tau}}_2}\l(\alpha\bY^{m+1}_h+(1-\alpha)\bY^{m}_h\r) \r]\ \d \alpha\ \d s\\
    &=\int_{\widehat{\Gamma}_h}[\bY^{m+1}_h-\bY^m_h]\cdot \l(\int^1_0\alpha^2\d \alpha \ \cJ(\bY^{m+1}_h)+ \int^1_0(1-\alpha)^2\d \alpha \ \cJ(\bY^{m}_h)\r.\\
    &\qquad\qquad\qquad\qquad \l.+\int^1_0\alpha(1-\alpha)\d \alpha\l[4\cJ(\bY^{m+1/2}_h)-\cJ(\bY^{m+1}_h)-\cJ(\bY^{m}_h)\r]\r)\ \d s\\
	&=\int_{\widehat{\Gamma}_h}[\bY^{m+1}_h-\bY^m_h]\cdot \frac{1}{6}\l[\cJ(\bY^m_h)+4\cJ(\bY^{m+1/2}_h)+\cJ(\bY^{m+1}_h)\r]\ \d s \\
        &=\int_{\widehat{\Gamma}_h}^h[\bY^{m+1}_h-\bY^m_h]\cdot \frac{1}{6}\l[\cJ(\bY^m_h)+4\cJ(\bY^{m+1/2}_h)+\cJ(\bY^{m+1}_h)\r]\ \d s \\	
    &=\big(\bX^{m+1}_h-\id,\bn^{m+\frac{1}{2}}_h\big)_{\Gamma^{m}_h}^h=\big(\nabla_{\Gamma^m_h}H^{m+1}_h,\nabla_{\Gamma^m_h}\tau\big)^h_{\Gamma^m_h}=0,
	\end{aligned}
\end{multline*}
where we have noted that the integrand
\[[\bY^{m+1}_h-\bY^m_h]\cdot\frac{1}{6}\l[\cJ(\bY^m_h)+4\cJ(\bY^{m+1/2}_h)+\cJ(\bY^{m+1}_h)\r]\]
is a piecewise polynomial of degree $3\ell-2$ on $\widehat{\Gamma}_h$. We can then complete the proof of the volume-preservation property by induction.\end{proof}

\section{Numerical results}

We have implemented the proposed isoparametric finite element schemes in the framework \textsc{Dune} \cite{Bastian2021,Sander2020,dune} using parametrized grids based on the \textsc{dune-curvedgrid} and \textsc{dune-curvedgeometry} \cite{Praetorius2022,dune-curvedgeometry,dune-curvedgrid} modules. As the underlying grid manager for the reference grid, we have used \textsc{dune-foamgrid} \cite{Sander2017,dune-foamgrid}. Our implementations and all the parameter files necessary to perform the presented numerical computations are available at the GitHub repository \cite{higher-order-bgn}. Unless specified explicitly, the reference grids are chosen as poly\-gonal/poly\-hedral approximations of the initial curves/surfaces with good mesh quality. For computations with $\P^{\ell}$-elements, the Gauss quadrature rules are uniformly chosen as exact to degree $2\ell-1$ for curves and exact to degree $3\ell-2$ for surfaces, respectively. Hence these quadrature rules always satisfy the assumptions in our theorems.

\subsection{Experimental order of convergence}

We test the experimental order of convergence for both mean curvature flow and surface diffusion of curves and surfaces. First, for mean curvature flow, the initial shape is chosen as the unit circle $x^2+y^2=1$ for the curve case and as the unit sphere  $x^2+y^2+z^2=1$ for the surface case. 
In order to test the expected convergence order $\mathcal{O}(\tau+h^{\ell+1})$ for $\P^\ell$-finite elements, we choose for the time step size the parameter $\tau=ch^{\ell+1}$ and define the approximate $L^\infty$-error between the numerical solution and the exact solution at time $T$ as
\begin{equation*}
    \cE_{L^\infty}(\ell,h,\tau)\colonequals d_{L^{\infty}}\l(\Gamma^{T/\tau}_h,\Gamma(T)\r).
\end{equation*}
For the curve case, $\Gamma(T)$ is a circle with radius $\sqrt{1-2T}$, and for the surface case, $\Gamma(T)$ is a sphere with radius $\sqrt{1-4T}$. The distance function $d_{L^{\infty}}$ is the maximum of the distance between the true circle/sphere and $\Gamma^{T/\tau}_h$ at all the quadrature points. The experimental order of convergence is then given by
\[
\mathrm{Order} \colonequals \frac{\log(\cE_{L^\infty}(\ell,h,\tau)/\cE_{L^\infty}(\ell,h/2,\tau/2^{\ell+1})) }{\log 2}.
\]
The numerical results are depicted in \Cref{Tab:mcf:curves,Tab:mcf:surfaces}
for the computations in 2d and 3d, respectively. As expected, we can observe that the spatial convergence order for $\P^\ell$-finite elements is $\ell+1$.

\begin{table}[h!]
\centering
\renewcommand{\arraystretch}{1.7}
\def\temptablewidth{0.8\textwidth}
\vspace{-2pt}
\caption{Convergence order test for curvature flow for a shrinking circle, using $\P^\ell$-elements with $\ell=1,2,3$. The parameters are chosen as $T=0.05$, $(h_0,\tau_0)=(0.2,0.05)$.}
{\rule{\temptablewidth}{1pt}}
\begin{tabular*}{\temptablewidth}{@{\extracolsep{\fill}}ccccccc}
 $(h,\tau)$
 & $\cE_{L^\infty}(1,h,\tau)$&  Order & $\cE_{L^\infty}(2,h,\tau)$ & Order & $\cE_{L^\infty}(3,h,\tau)$ & Order  \\  \hline
$(h_0,\tau_0)$  &  2.04e-3 &  ---  &  3.70e-3  & --- &  3.70e-3   & --- \\   \hline
$(\frac{h_0}{2^1},\frac{\tau_0}{2^{\ell+1}})$  & 1.24e-3  & 0.72  & 5.14e-4     & 2.85 & 2.59e-4   & 3.84 \\   \hline
$(\frac{h_0}{2^2},\frac{\tau_0}{2^{2(\ell+1)}})$ & 3.72e-5  & 1.74  & 6.52e-5    & 2.98  &  1.63e-5  &  3.99\\  \hline
$(\frac{h_0}{2^3},\frac{\tau_0}{2^{3(\ell+1)}})$  &  9.80e-6 &  1.92 &  8.18e-6   & 3.00  &  1.02e-6  &  4.00\\ \hline
$(\frac{h_0}{2^4},\frac{\tau_0}{2^{4(\ell+1)}})$  & 2.50e-6  & 1.97  & 1.02e-6    & 3.00   &   6.36e-8 & 4.00  \\
 \end{tabular*}
{\rule{\temptablewidth}{1pt}}
\label{Tab:mcf:curves}
\end{table}

\begin{table}[h!]
\centering
\renewcommand{\arraystretch}{1.7}
\def\temptablewidth{0.8\textwidth}
\vspace{-2pt}
\caption{Convergence order test for mean curvature flow for a shrinking sphere, using $\P^\ell$-elements with $\ell=1,2,3$. The parameters are chosen as $T=0.05$, $(h_0,\tau_0)=(0.5,0.05)$.}
{\rule{\temptablewidth}{1pt}}
\begin{tabular*}{\temptablewidth}{@{\extracolsep{\fill}}ccccccc}
 $(h,\tau)$
 & $\cE_{L^\infty}(1,h,\tau)$&  Order & $\cE_{L^\infty}(2,h,\tau)$ & Order & $\cE_{L^\infty}(3,h,\tau)$ & Order  \\  \hline
$(h_0,\tau_0)$  &  1.60e-1  &  ---  & 1.43e-2   & --- & 2.37e-2  & --- \\   \hline
$(\frac{h_0}{2^1},\frac{\tau_0}{2^{\ell+1}})$  &  5.13e-2   & 1.64  & 3.48e-3  &  2.04  & 1.66e-3  & 3.84 \\   \hline
$(\frac{h_0}{2^2},\frac{\tau_0}{2^{2(\ell+1)}})$  & 1.44e-2    & 1.83  & 4.90e-4  &     2.83 & 9.33e-5 &  4.15\\ \hline
$(\frac{h_0}{2^3},\frac{\tau_0}{2^{3(\ell+1)}})$  & 3.67e-3    &  1.97  & 5.39e-5  &   3.18 & 5.67e-6  &  4.04 \\
 \end{tabular*}
{\rule{\temptablewidth}{1pt}}
\label{Tab:mcf:surfaces}
\end{table}

A similar convergence order experiment can be conducted for surface diffusion of curves and surfaces. Here, we choose as the initial shape the ellipse $x^2/4+y^2=1$ and the ellipsoid $x^2/4+y^2+z^2=1$, respectively. 
Since for surface diffusion we do not have exact solutions, we define the approximate $L^2$-error under the distance function at time $T$ as
\begin{equation}\label{L2error}
		\cE_{L^2}(\ell,h,\tau) \colonequals d_{L^2}\l(\Gamma^{T/\tau}_{h},\Gamma^{T/(2^{\ell+1}\tau)}_{h/2}\r),
\end{equation}
see \eqref{eq:dL2} in \Cref{sec:appendix:hausdorff-distance} for a precise definition. 
The experimental order of convergence is defined similarly. \Cref{Tab:sd:curves,Tab:sd:surfaces} show that the $\P^\ell$-finite element methods achieve the expected convergence orders for both curve and surface cases, respectively.

\begin{remark}
The $L^2$-distance function defined in \eqref{eq:dL2} is called the projected distance in \cite{Bai2023}. Similarly, the $L^1$-distance (called the manifold distance in \cite{Jiang23,Zhao2021}) and the $L^\infty$-distance (called the Hausdorff distance in \cite{Jiang23}) are also widely used to characterize the error for the evolution of curves and surfaces. Similar convergence results hold for these errors; we omit the numerical results for simplicity.
\end{remark}

\begin{table}[h!]
\centering
\renewcommand{\arraystretch}{1.7}
\def\temptablewidth{0.8\textwidth}
\vspace{-2pt}
\caption{Convergence order test for surface diffusion for an initial ellipse, using $\P^\ell$-elements with $\ell=1,2,3$. The parameters are chosen as $T=0.05$, $(h_0,\tau_0)=(0.3,0.05)$.}
{\rule{\temptablewidth}{1pt}}
\begin{tabular*}{\temptablewidth}{@{\extracolsep{\fill}}ccccccc}
 $(h,\tau)$
 & $\cE_{L^2}(1,h,\tau)$&  Order & $\cE_{L^2}(2,h,\tau)$ & Order & $\cE_{L^2}(3,h,\tau)$ & Order  \\  \hline
$(h_0,\tau_0)$  & 3.84e-3  &  ---  &  9.08e-5 & --- &  1.27e-5 & --- \\   \hline
$(\frac{h_0}{2^1},\frac{\tau_0}{2^{\ell+1}})$  & 9.68e-4  &  1.99  & 1.19e-5  & 2.93   &  7.99e-7 & 4.00 \\   \hline
$(\frac{h_0}{2^2},\frac{\tau_0}{2^{2(\ell+1)}})$  & 2.43e-4 & 2.00  &  1.48e-6  & 3.01   & 4.71e-8 & 4.09  \\  \hline
$(\frac{h_0}{2^3},\frac{\tau_0}{2^{3(\ell+1)}})$  & 6.07e-5  &  1.99 & 1.78e-7  & 3.05   & 3.02e-9  &  3.96 \\
 \end{tabular*}
{\rule{\temptablewidth}{1pt}}
\label{Tab:sd:curves}
\end{table}

\begin{table}[h!]
\centering
\renewcommand{\arraystretch}{1.7}
\def\temptablewidth{0.8\textwidth}
\vspace{-2pt}
\caption{Convergence order test for surface diffusion for an initial ellipsoid, using $\P^\ell$-elements with $\ell=1,2,3$. The parameters are chosen as $T=0.05$, $(h_0,\tau_0)=(0.5,0.05)$.}
{\rule{\temptablewidth}{1pt}}
\begin{tabular*}{\temptablewidth}{@{\extracolsep{\fill}}ccccccc}
 $(h,\tau)$
 & $\cE_{L^2}(1,h,\tau)$&  Order & $\cE_{L^2}(2,h,\tau)$ & Order & $\cE_{L^2}(3,h,\tau)$ & Order  \\  \hline
$(h_0,\tau_0)$  & 7.52e-2  &  ---  &  8.81e-3 & --- & 4.78e-3  & --- \\   \hline
$(\frac{h_0}{2^1},\frac{\tau_0}{2^{\ell+1}})$  & 2.33e-2  &  1.69   & 1.90e-3   & 2.37    & 3.54e-4  & 3.76 \\   \hline
$(\frac{h_0}{2^2},\frac{\tau_0}{2^{2(\ell+1)}})$  &  6.43e-3   &  1.86  &  2.19e-4  & 3.12   & 3.00e-5  & 3.56  \\
 \end{tabular*}
{\rule{\temptablewidth}{1pt}}
\label{Tab:sd:surfaces}
\end{table}

\subsection{Evolution by mean curvature flow}\label{sec:mcf}

In this subsection, we conduct some morphology and geometric quantity evolutions for mean curvature flow by our $\P^\ell$-finite element method.

\begin{figure}[h!]
    \begin{subfigure}{0.25\linewidth}
        \includegraphicswithlegend[1.0]{0.8\textwidth}{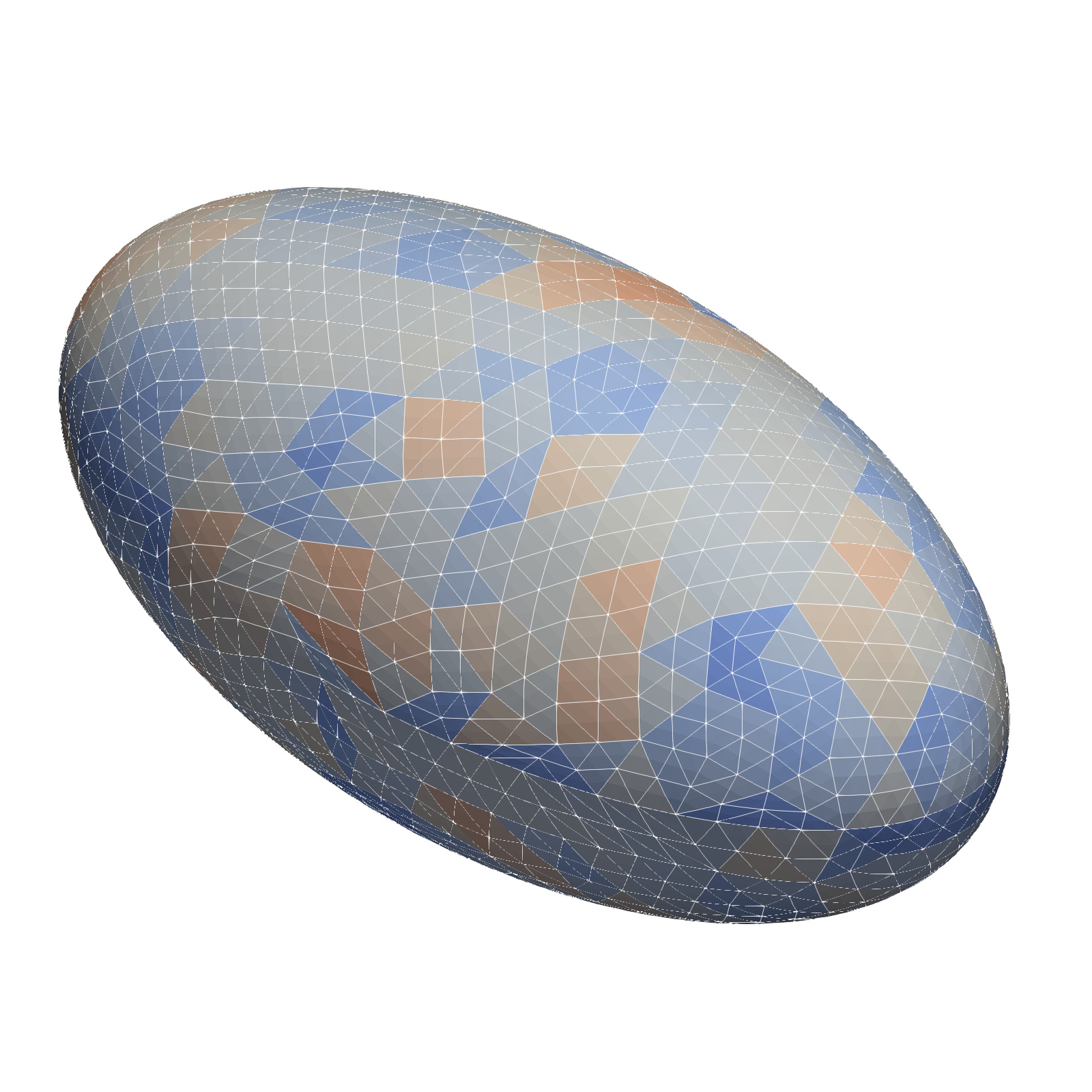}{1.4e-3}{1.7e-2}
    \end{subfigure}\hfill%
    \begin{subfigure}{0.25\linewidth}
        \includegraphicswithlegend[0.95]{0.8\textwidth}{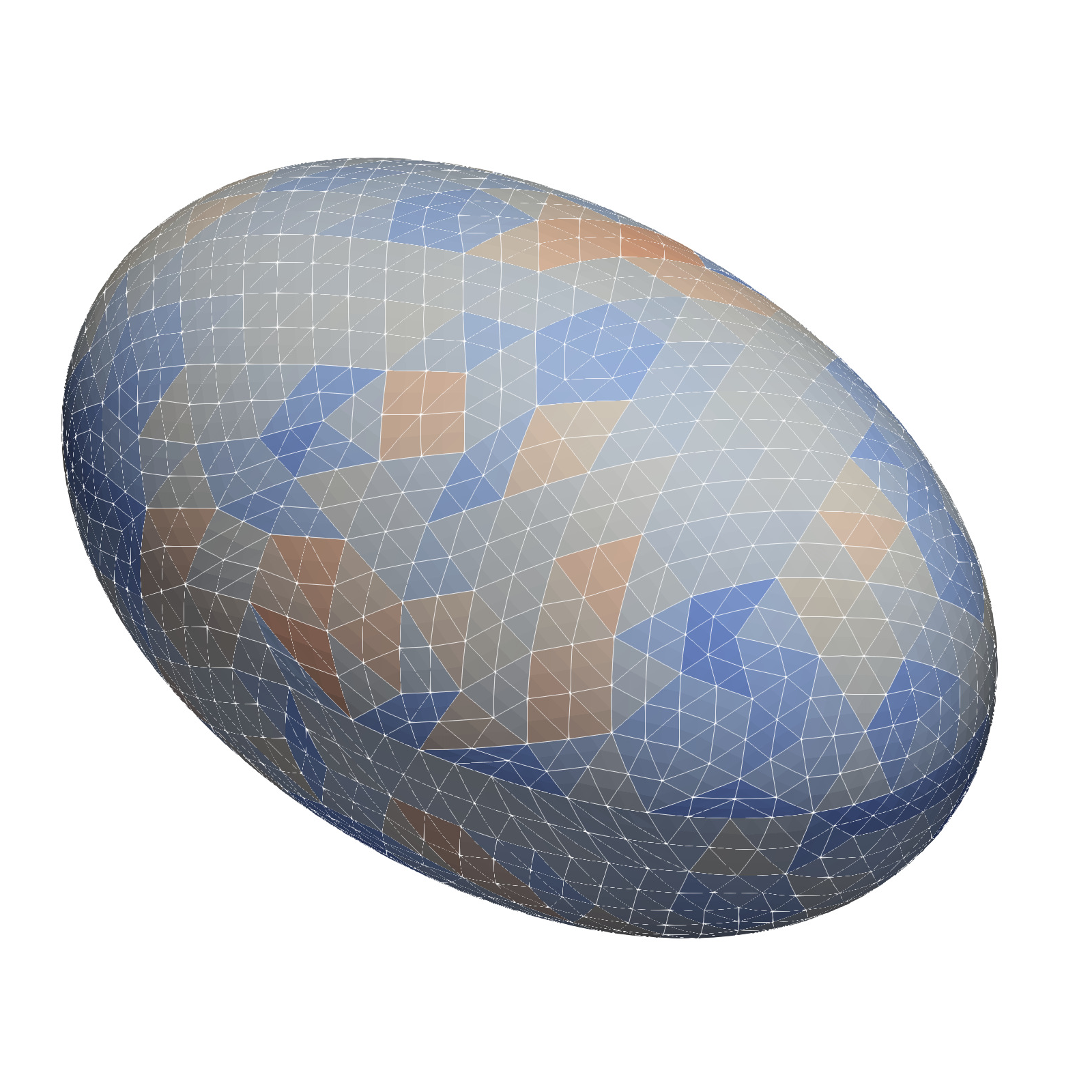}{4.6e-4}{1.7e-2}
    \end{subfigure}\hfill%
    \begin{subfigure}{0.25\linewidth}
        \includegraphicswithlegend[0.9]{0.8\textwidth}{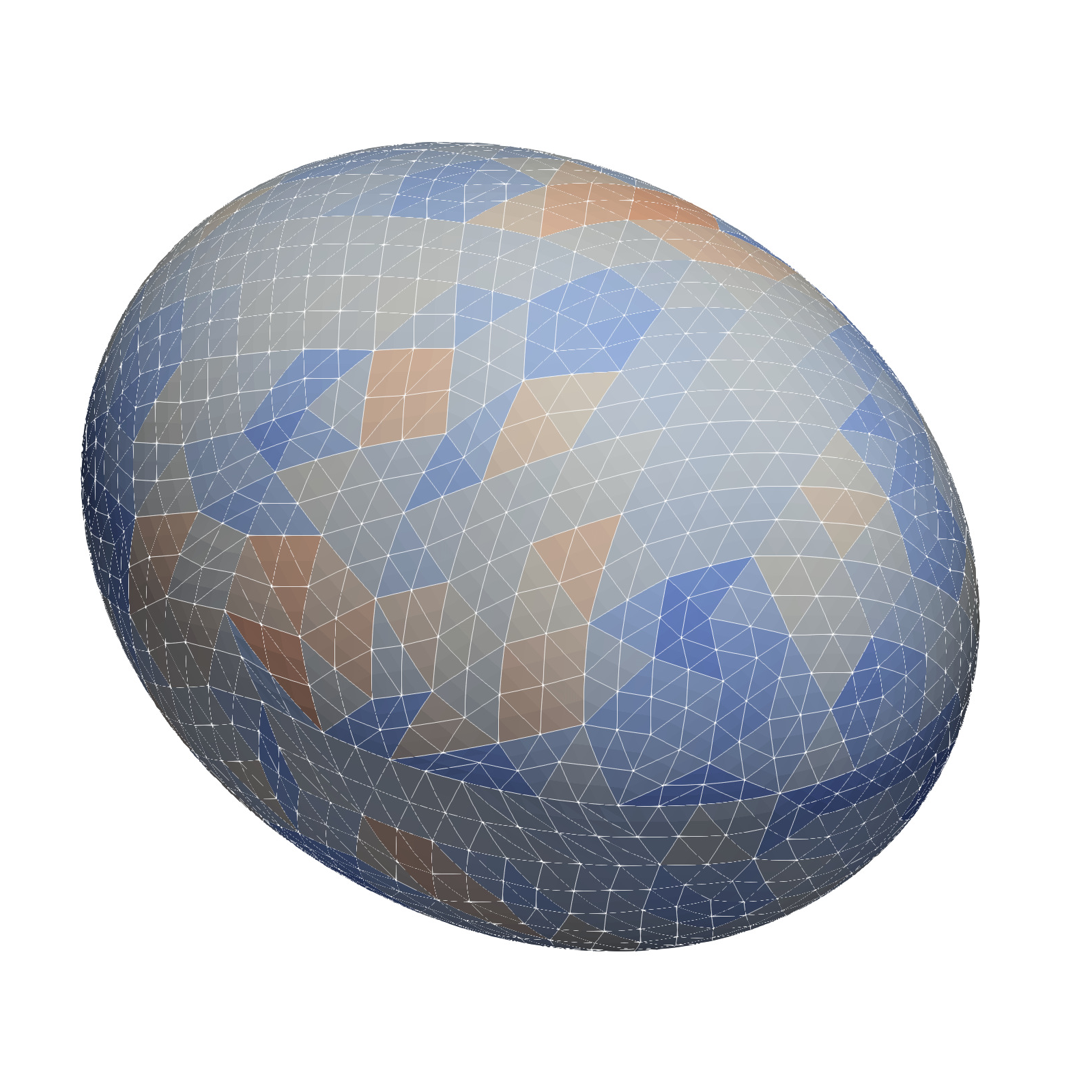}{9.4e-5}{1.1e-3}
    \end{subfigure}\hfill%
    \begin{subfigure}{0.25\linewidth}
        \includegraphicswithlegend[0.85]{0.8\textwidth}{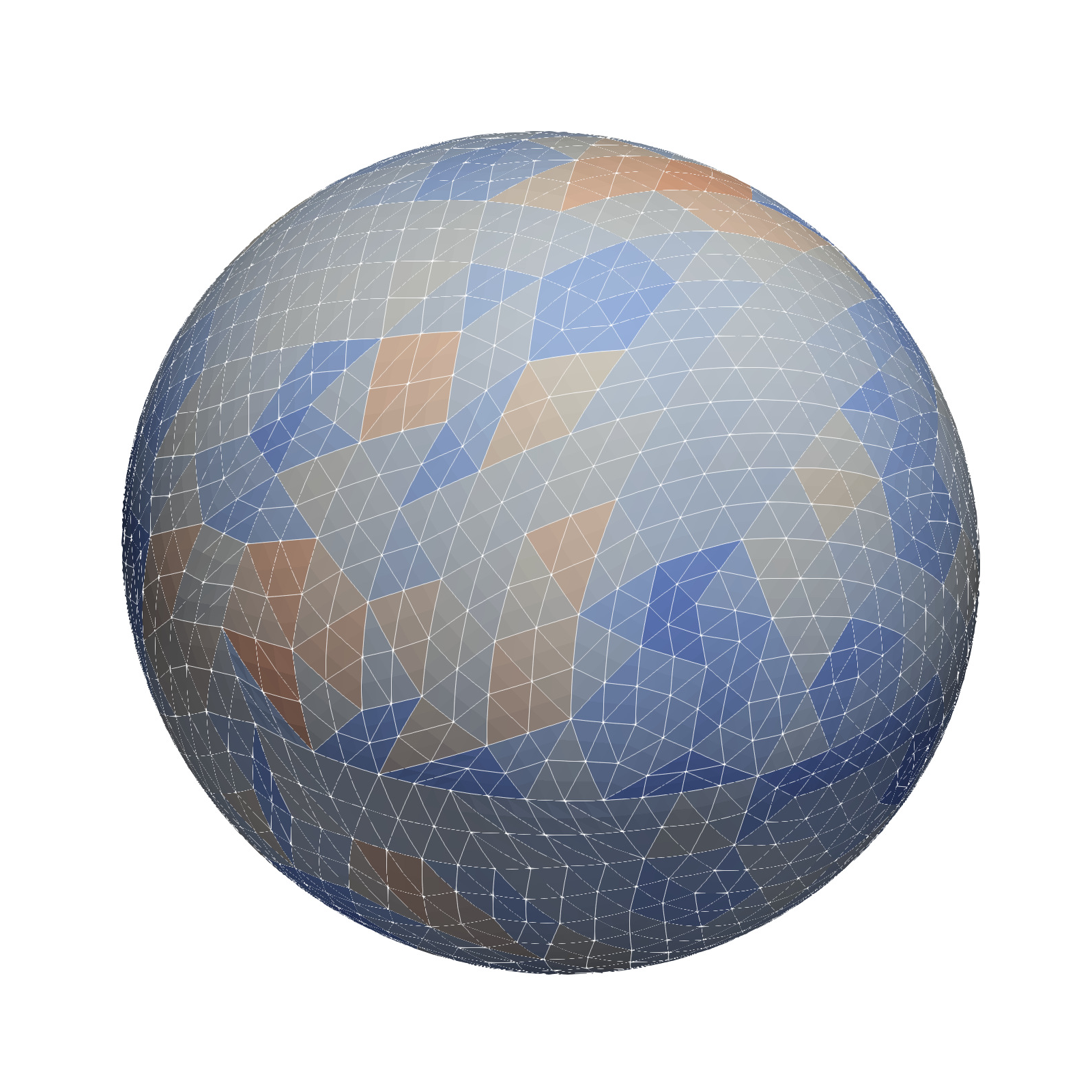}{1.4e-6}{1.7e-5}
    \end{subfigure}
    \caption{Evolution of an ellipsoid by mean curvature flow, using BGN-type scheme with $\ell=2$. The figures are rescaled, and the times are $t=0,0.3, 0.35, 0.377$, respectively. The colors indicate the area of the grid elements.}\label{mcf_ellipsoid_BGN}
\end{figure}

\begin{figure}[h!]
    \begin{subfigure}{0.25\linewidth}
        \includegraphicswithlegend[1.0]{0.8\textwidth}{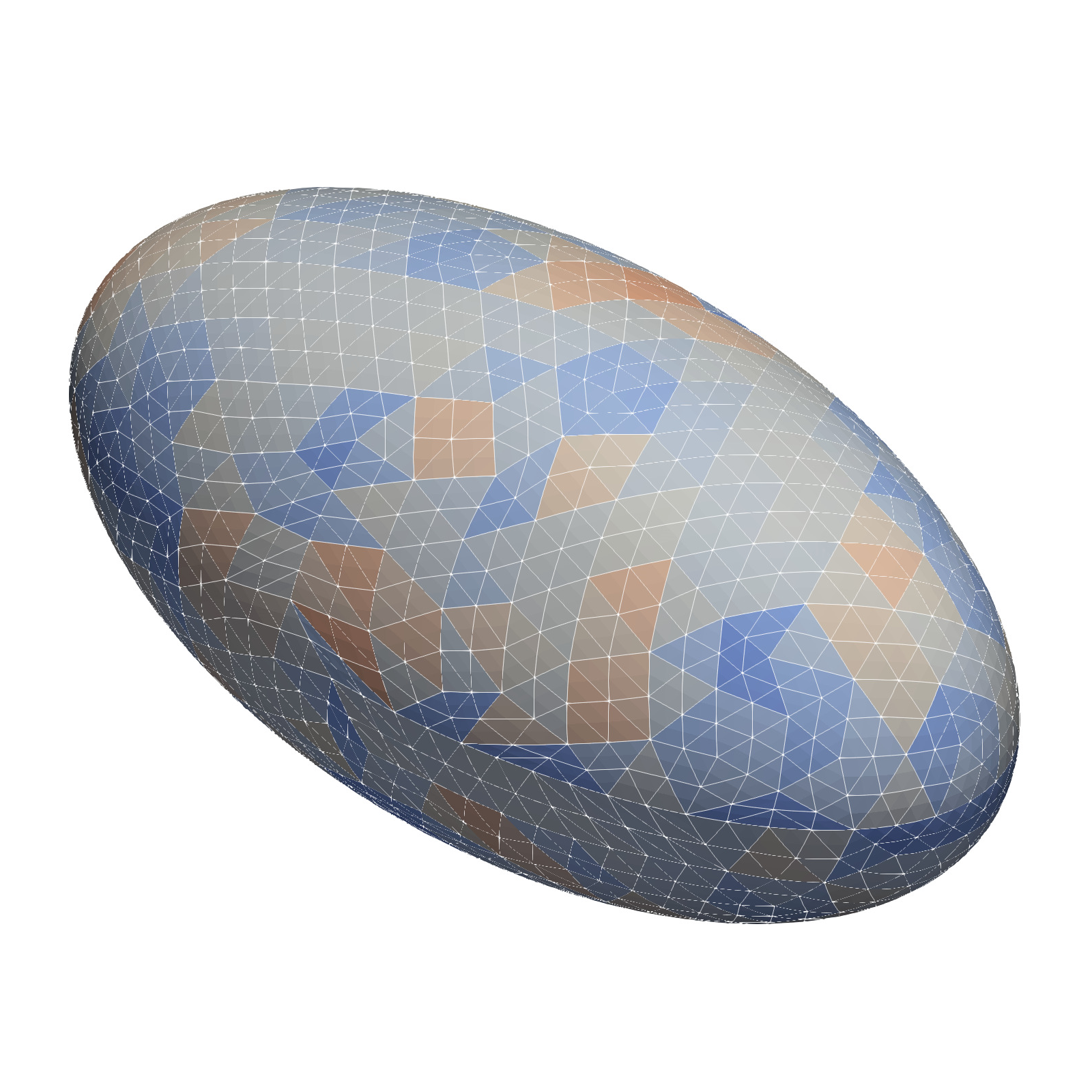}{1.4e-3}{1.7e-2}
    \end{subfigure}%
    \begin{subfigure}{0.25\linewidth}
        \includegraphicswithlegend[0.95]{0.8\textwidth}{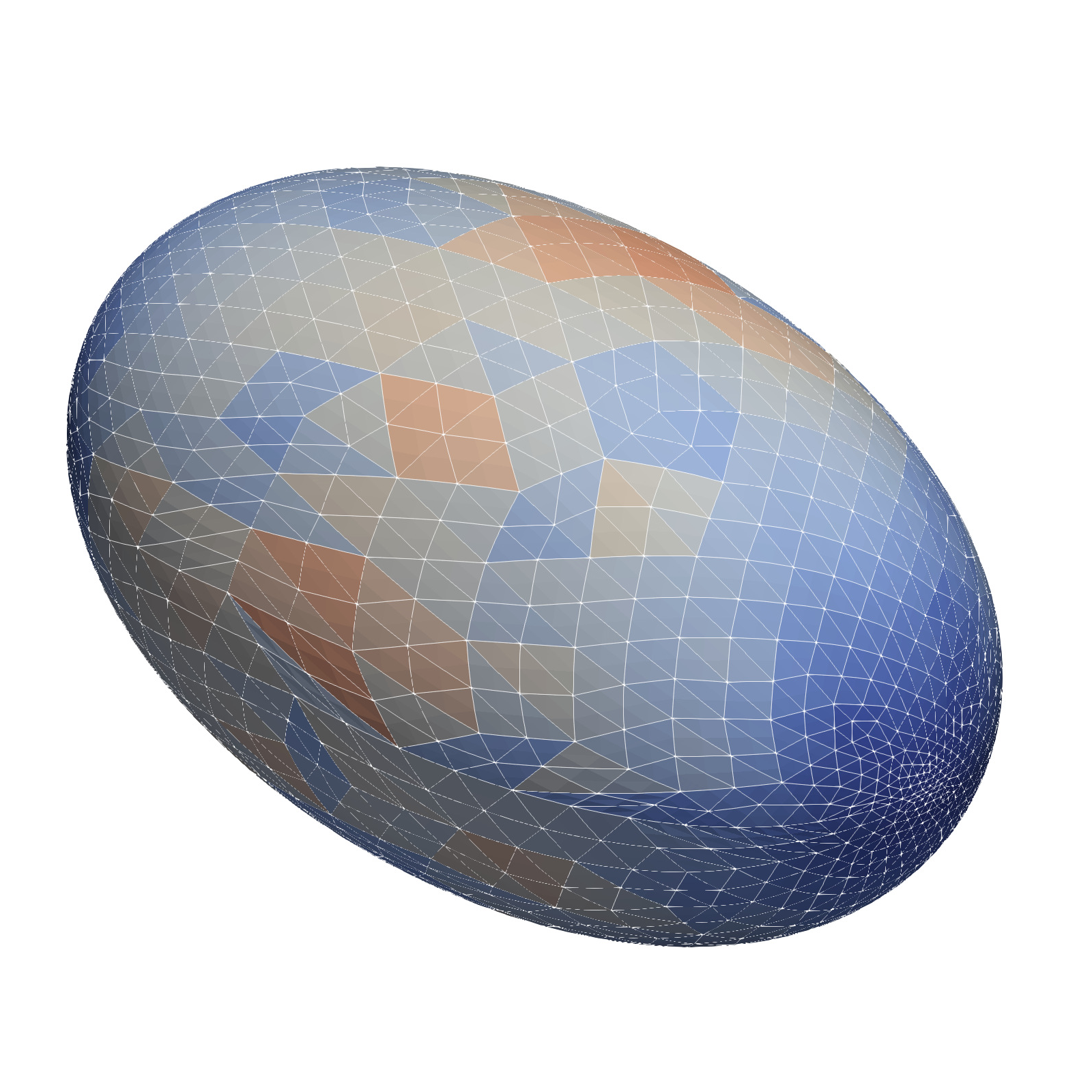}{7.3e-5}{9.0e-3}
    \end{subfigure}%
    \begin{subfigure}{0.25\linewidth}
        \includegraphicswithlegend[0.9]{0.8\textwidth}{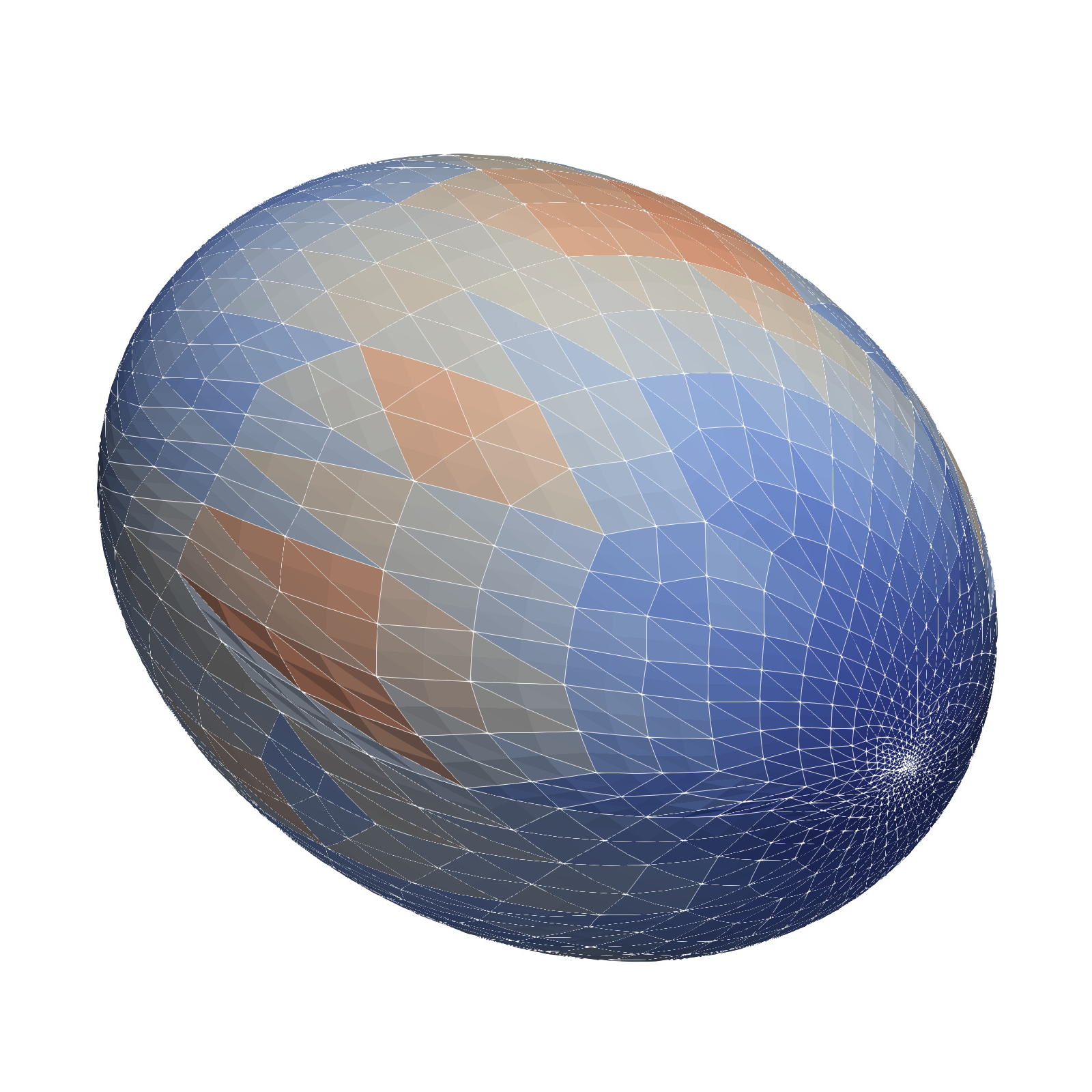}{1.7e-6}{3.2e-3}
    \end{subfigure}%
    \begin{subfigure}{0.25\linewidth}
        \includegraphicswithlegend[0.85]{0.8\textwidth}{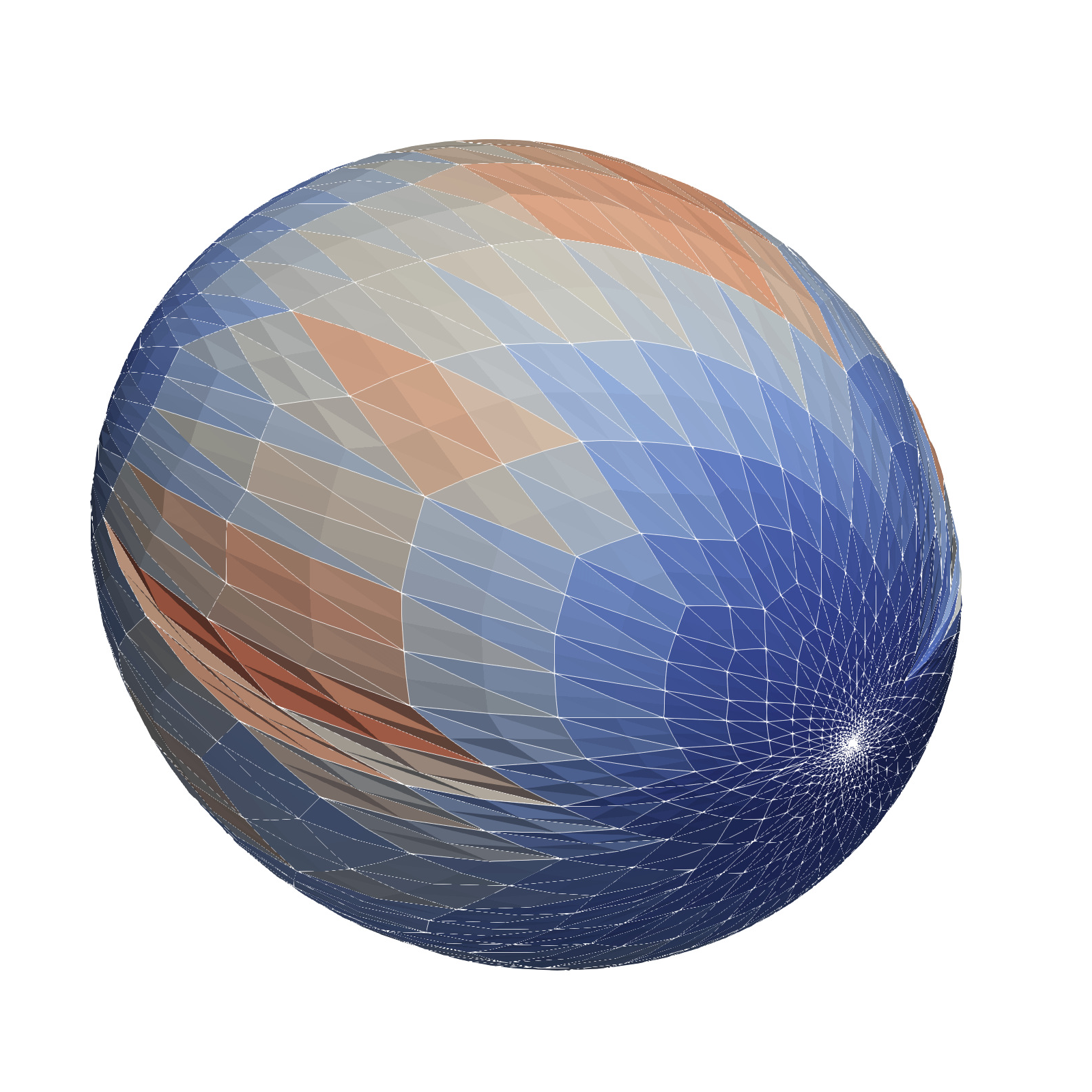}{6.0e-9}{1.1e-4}
    \end{subfigure}
    \caption{Evolution of an ellipsoid by mean curvature flow, using Dziuk-type schemes with $\ell=2$. The figures are rescaled, and the times are $t=0,0.3, 0.35, 0.377$, respectively. The colors indicate the area of the grid elements.}\label{mcf_ellipsoid_Dziuk}
\end{figure}

We first compare the mesh quality between our isoparametric method \eqref{eq:3dnum} and the Dziuk-type scheme \eqref{Dziuk-typeh}. As discussed in \cite{BGN08A} for $\P^1$-elements, Dziuk's method may lead to mesh distortion, while the BGN method achieves good mesh quality. This phenomenon also appears for higher spatial order methods proposed in this paper. As an example, we take an ellipsoid as the initial surface, with parameters $(J,K)=(2704,1354)$, where $J$ represents the number of triangles and $K$ represents the number of vertices. We use $\P^2$-elements for both the BGN-type scheme and the Dziuk-type scheme, and we monitor the mesh quality by comparing the area of each triangle $|\sigma_j^m|$ computed by a suitable quadrature rule. As illustrated in \Cref{mcf_ellipsoid_BGN,mcf_ellipsoid_Dziuk}, the BGN-type scheme successfully evolves the initial ellipsoid to a spherical point. At the same time, the mesh distortion phenomenon may appear for the Dziuk-type scheme (e.g., the last snapshot of \Cref{mcf_ellipsoid_Dziuk}). Moreover, the distribution of triangle sizes indicates that the Dziuk-type scheme leads to significantly different areas of triangles. As shown in the last snapshot of \Cref{mcf_ellipsoid_Dziuk}, the largest triangle is about $10^{5}$ times the area of the smallest triangle, while the last snapshot of \Cref{mcf_ellipsoid_BGN} shows the BGN-type method exhibits a similar area of triangles near the singularity. These numerical results indicate that the higher spatial order  BGN-type schemes can also adjust the mesh automatically, similar to the linear finite element methods \cite{BGN08A}.

\smallskip

\begin{figure}[h!]
\hspace{-12pt}
\begin{minipage}[t]{0.20\linewidth}
	\centering
\includegraphics[scale=0.055]{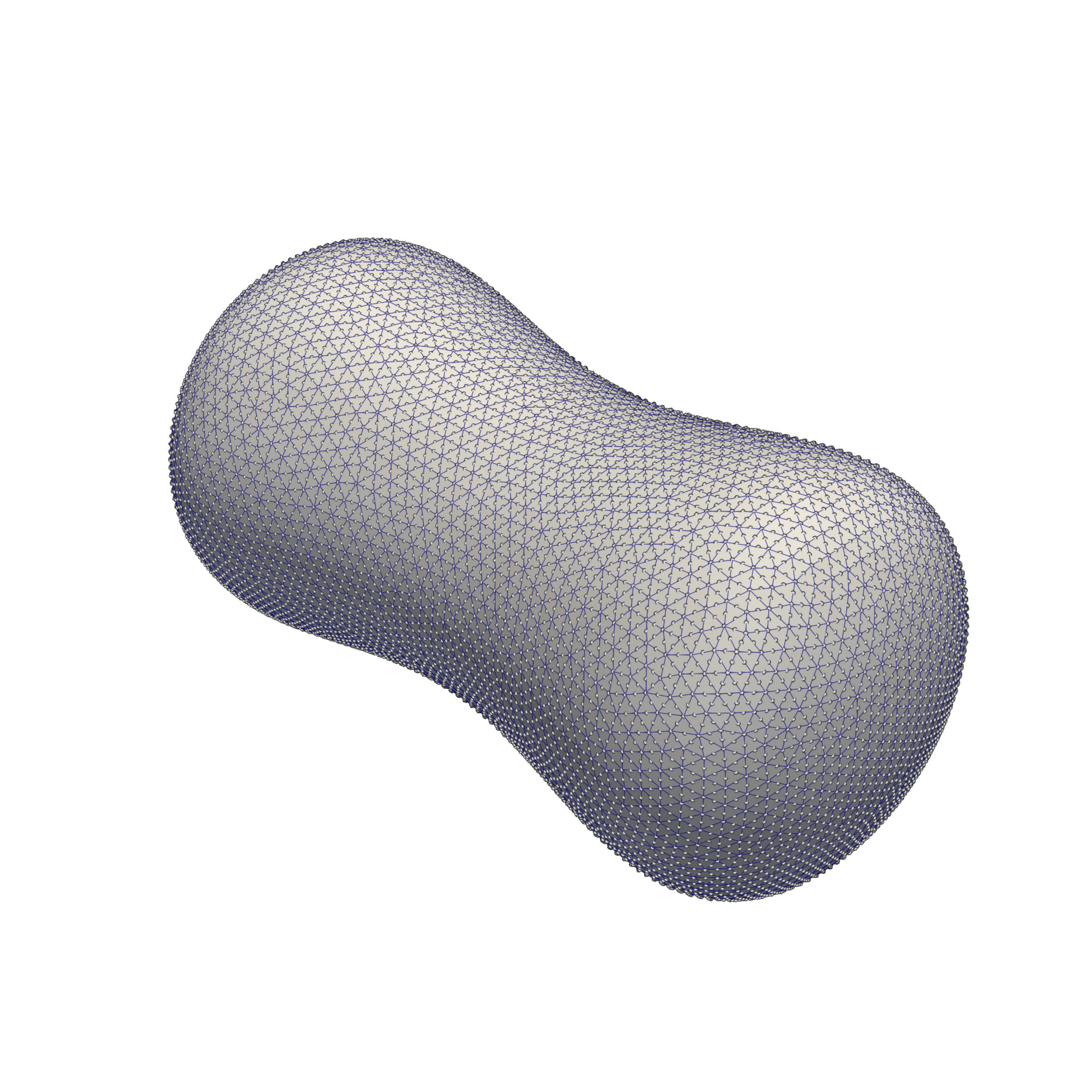}
\end{minipage}
\hspace{12pt}
\begin{minipage}[t]{0.20\linewidth}
	\centering
\includegraphics[scale=0.055]{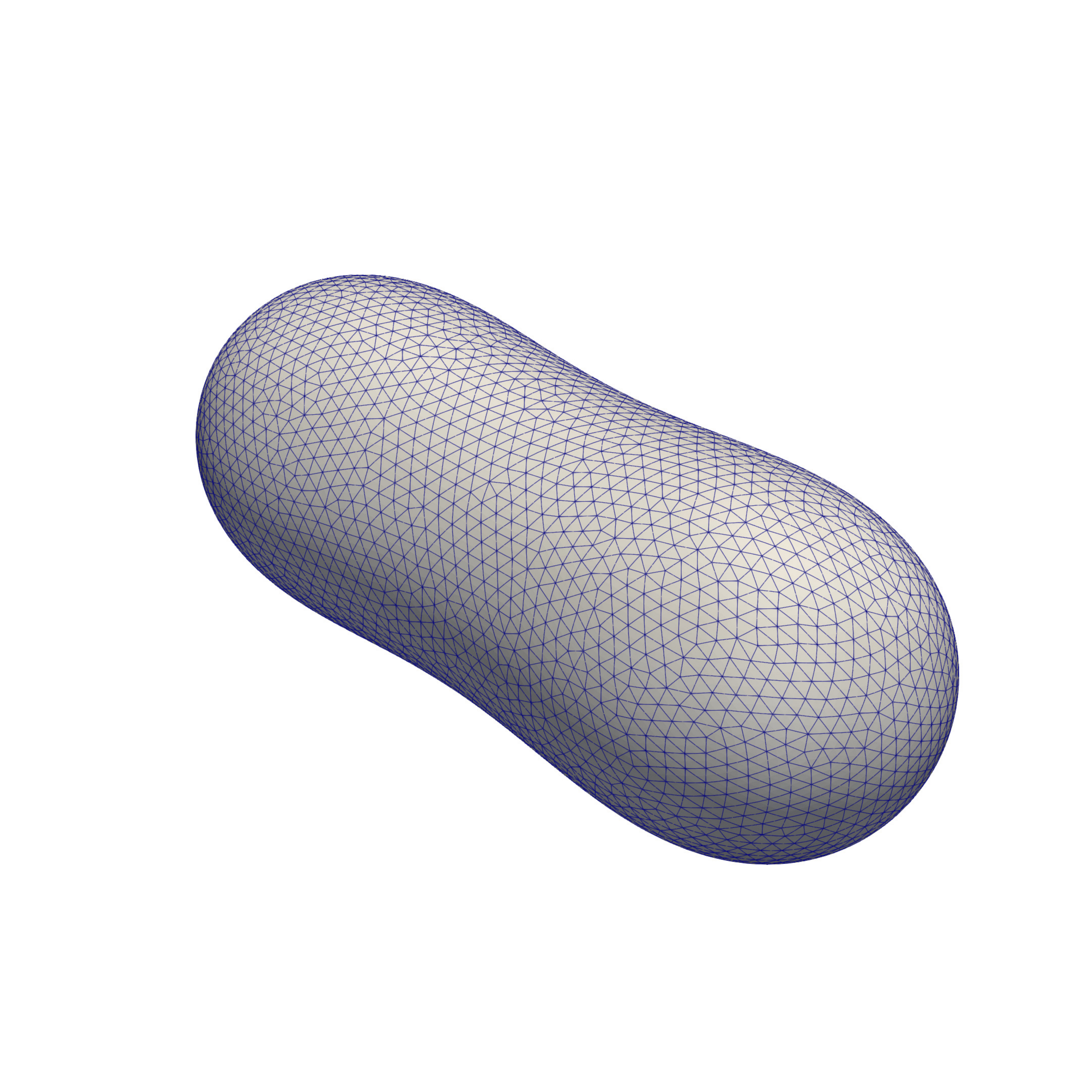}
\end{minipage}
\hspace{12pt}
\begin{minipage}[t]{0.20\linewidth}
	\centering
\includegraphics[scale=0.055]{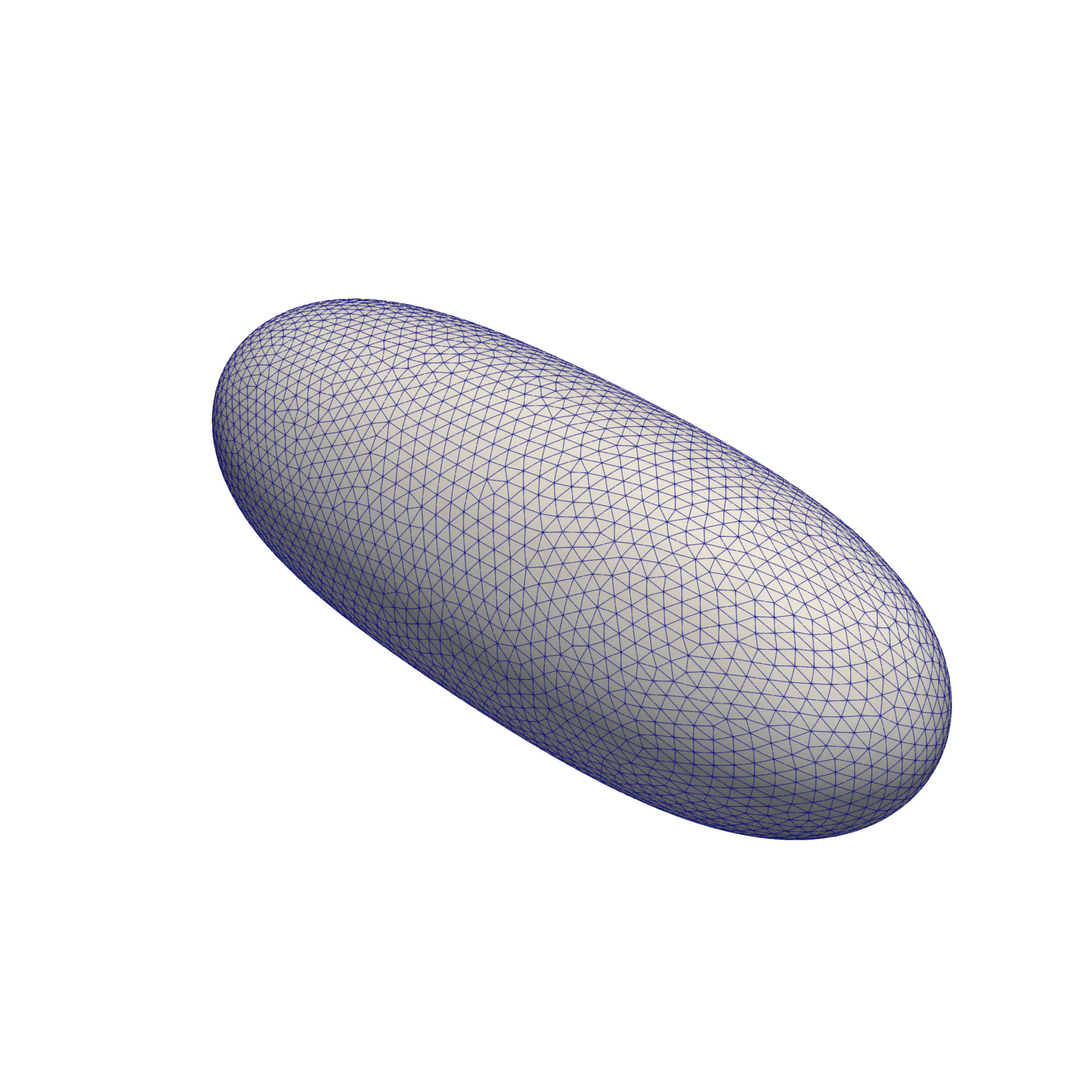}
\end{minipage}
\hspace{12pt}
\begin{minipage}[t]{0.20\linewidth}
	\centering
\includegraphics[scale=0.055]{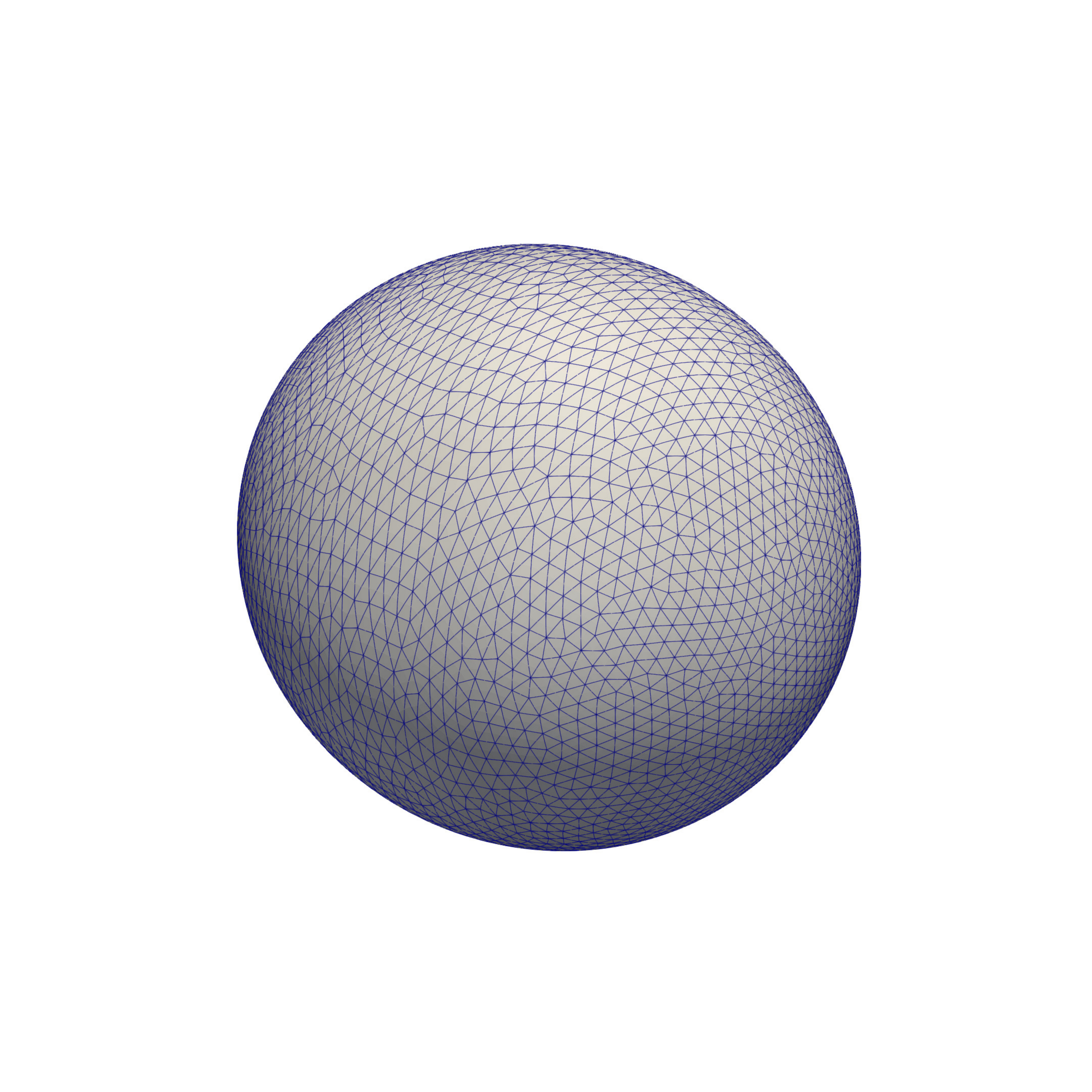}
\end{minipage}
\caption{Evolution of a mild dumbbell shape by mean curvature flow, using $\mathbb{P}^2$-elements, figures are rescaled. We used the discretization parameters $(J,K)=(7096,3550)$ and $\tau=0.0001$.  The times are  $t=0,0.040, 0.0800, 0.0916$, respectively.}
\label{mcf_dumbbell_fat}
\end{figure}

\begin{figure}[h!]
\hspace{-12pt}
\begin{minipage}[t]{0.20\linewidth}
	\centering
\includegraphics[scale=0.055]{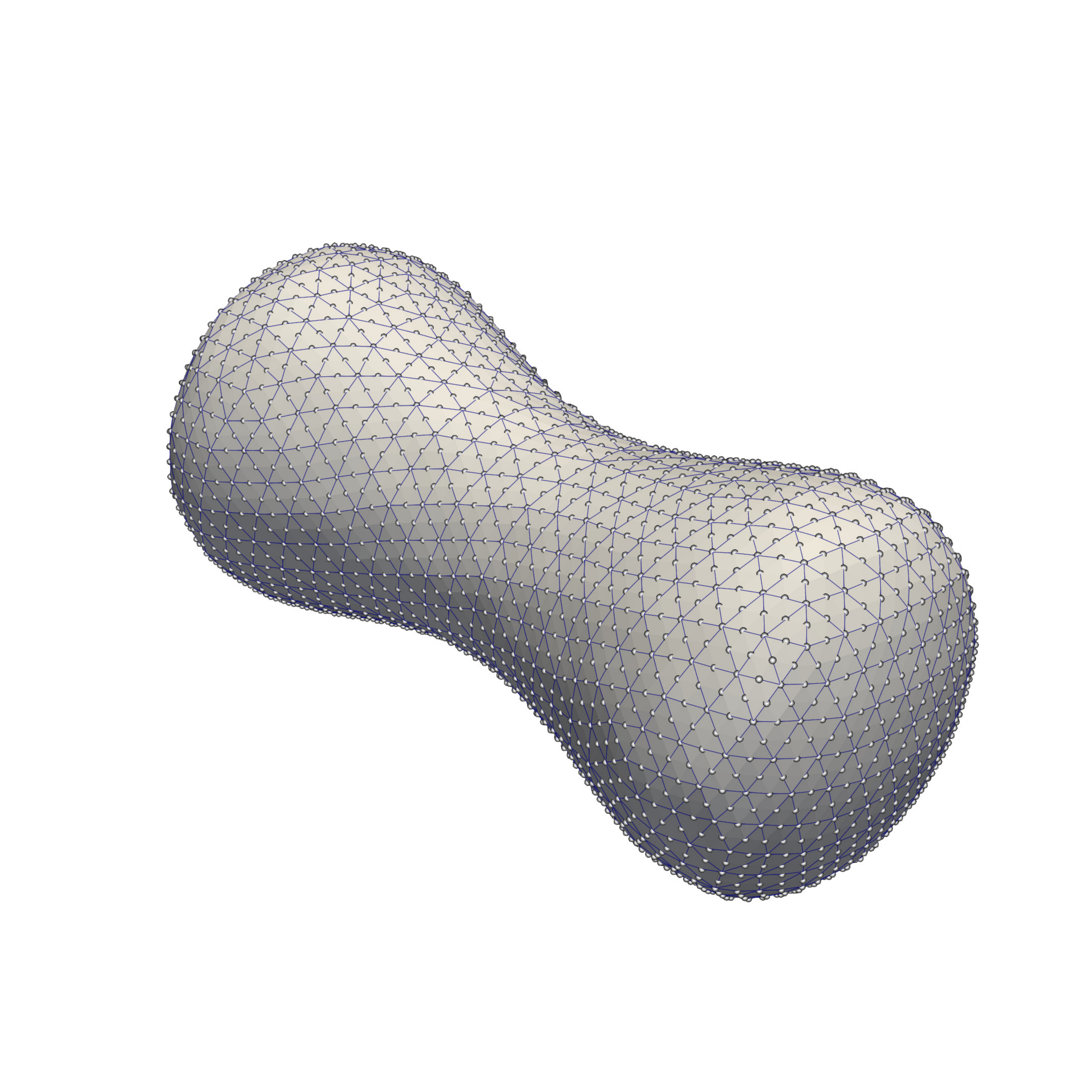}
\end{minipage}
\hspace{12pt}
\begin{minipage}[t]{0.20\linewidth}
	\centering
\includegraphics[scale=0.055]{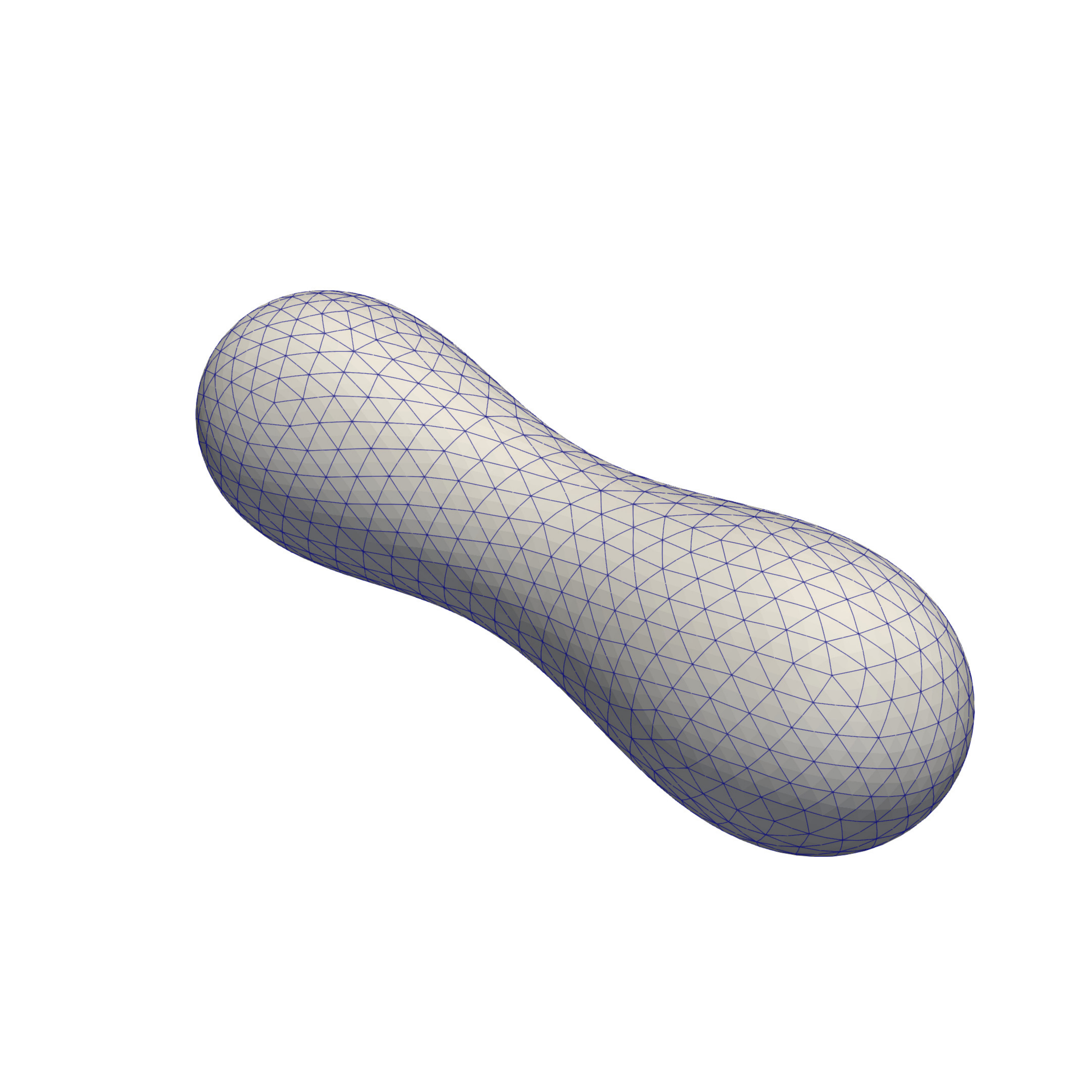}
\end{minipage}
\hspace{12pt}
\begin{minipage}[t]{0.20\linewidth}
	\centering
\includegraphics[scale=0.055]{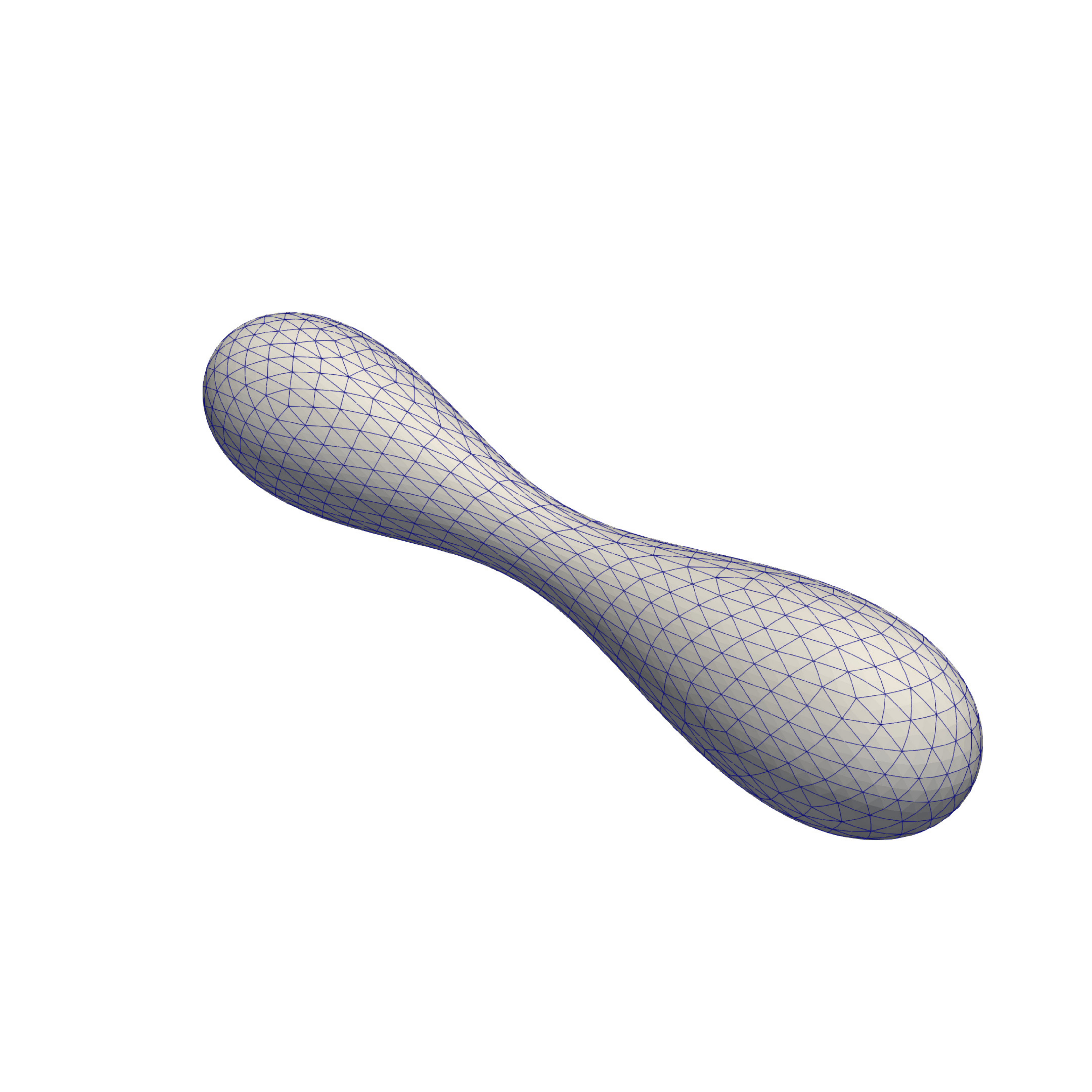}
\end{minipage}
\hspace{12pt}
\begin{minipage}[t]{0.20\linewidth}
	\centering
\includegraphics[scale=0.055]{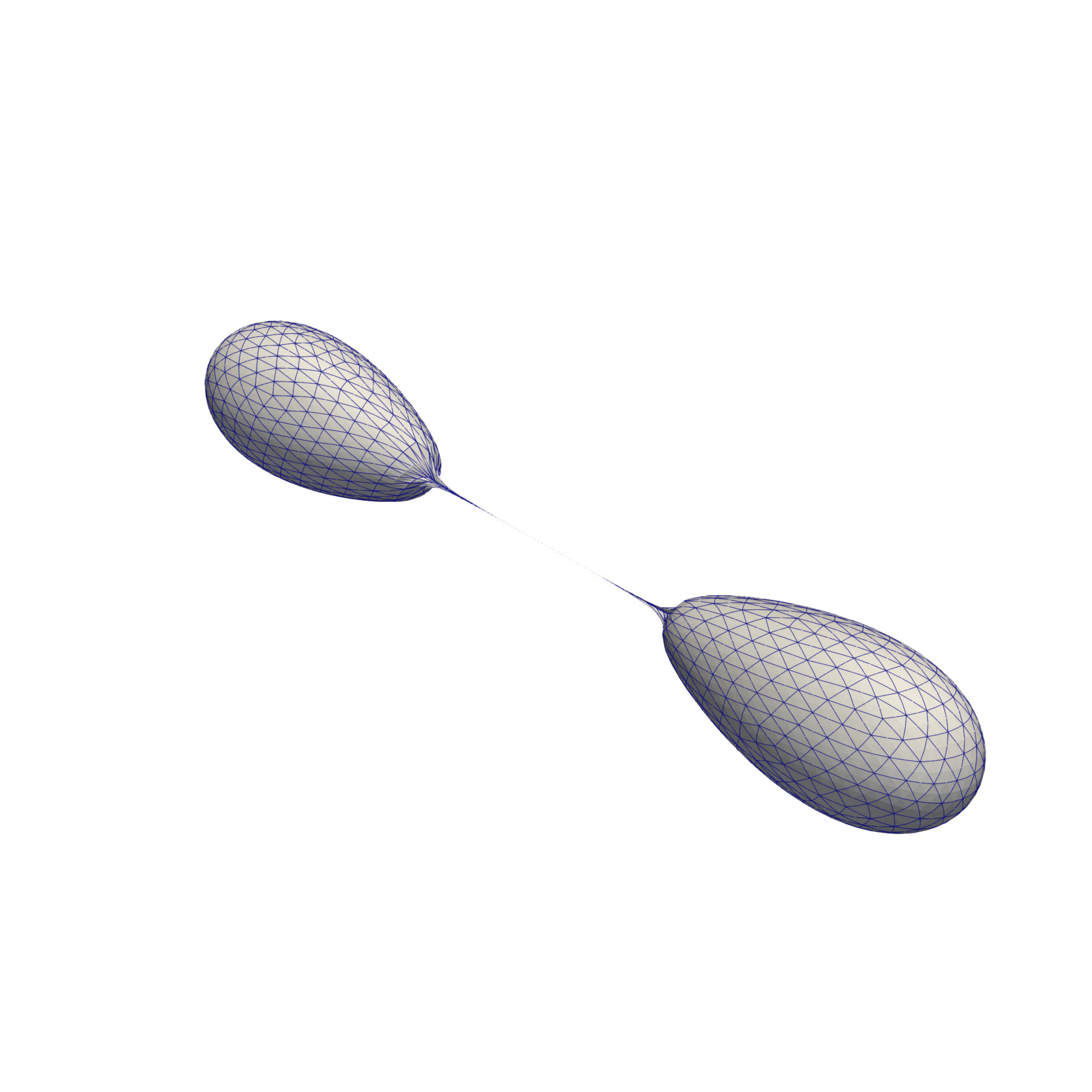}
\end{minipage}
\caption{Evolution of a dumbbell shape by mean curvature flow, using $\mathbb{P}^2$-elements. We used the discretization parameters $(J,K)=(3276,1430)$ and $\tau=0.0001$.  The times are  $t=0,0.040, 0.0555, 0.0585$, respectively.}
\label{mcf_dumbbell_thin}
\end{figure}

\begin{figure}[h!]
\hspace{0pt}
\begin{minipage}[t]{0.5\linewidth}
	\centering
\includegraphics[width=8.5cm,height=6cm]{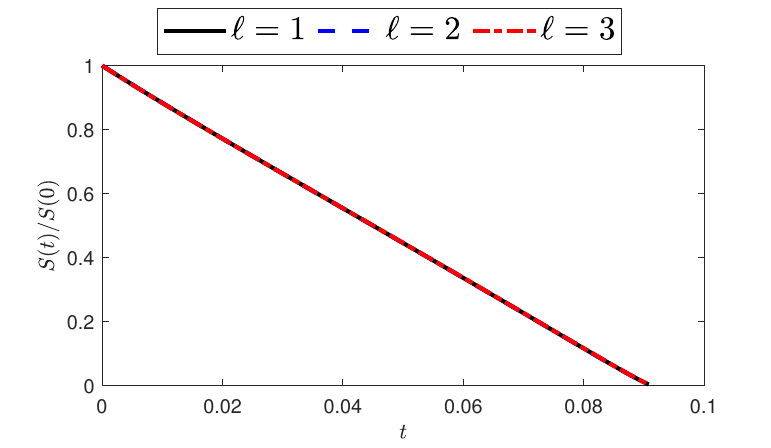}
\end{minipage}
\hspace{0pt}
\begin{minipage}[t]{0.5\linewidth}
	\centering
\includegraphics[width=8.5cm,height=6cm]{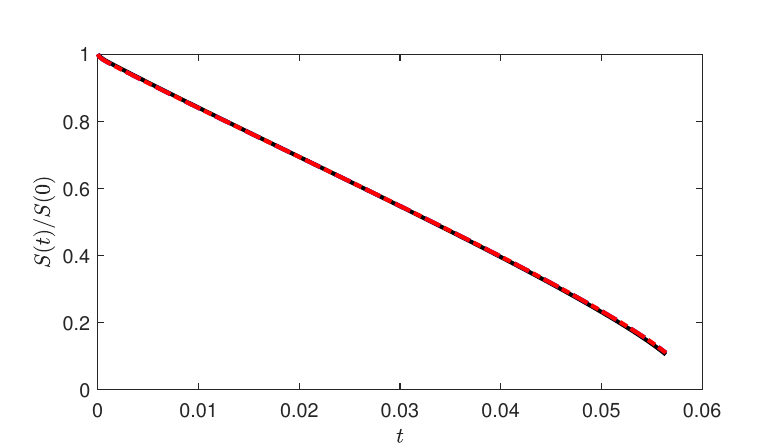}
\end{minipage}
\caption{Evolution of the normalized surface area for numerical simulations corresponding to \Cref{mcf_dumbbell_fat} (left figure) and  \Cref{mcf_dumbbell_thin} (right figure), respectively.}
\label{mcf_dumbbell_geo}
\end{figure}

We then apply our $\P^2$-finite element method to two dumbbell shape examples: a mild dumbbell shape defined by $(0.6x^2 + 0.4)^2(1 - x^2) - y^2 - z^2=0$, and a dumbbell shape given by $(0.7x^2 + 0.3)^2(1 - x^2) - y^2 - z^2=0$. The approximating surfaces employ the discretization parameters $(J,K)=(7096,3550)$ and $(J,K)=(3276,1430)$, respectively. As demonstrated in \Cref{mcf_dumbbell_fat,mcf_dumbbell_thin}, we observe that the mild dumbbell shape eventually shrinks to a spherical point, while the second dumbbell shape develops a singularity when the surface tries to pinch off into two. For these two examples, we also monitor the evolution of the surface area $S(t)/S(0)|_{t=t_m}=S_h^m/S_h^0$ in \Cref{mcf_dumbbell_geo} and notice that it is decreasing for both experiments until the singularity.

\subsection{Evolution by surface diffusion}

In this subsection, we present extensive numerical results of surface diffusion for both curves and surfaces. We consider  the evolution of the following geometric quantities: 
 \[
\left.\frac{L(t)}{L(0)}\right|_{t=t_m}=\frac{L_h^m}{L_h^0} ,\qquad \left.\frac{S(t)}{S(0)}\right|_{t=t_m}=\frac{S_h^m}{S_h^0},
\]
(2) the relative loss of area/volume defined as
\[
\l.\Delta A(t)\r|_{t=t_m}=\frac{A_h^m-A_h^0}{A_h^0},\qquad  \l.\Delta V(t)\r|_{t=t_m}=\frac{V_h^m-V_h^0}{V_h^0},
\]
and (3) the mesh quality function
\[
\left.\Psi(t)\right|_{t=t_m}=\frac{\max_{j}|\sigma_j^m|}{\min_{j}|\sigma_j^m|}.
\]

\begin{figure}[h!]
\hspace{16pt}
\begin{minipage}[t]{0.20\linewidth}
	\centering
\includegraphics[scale=0.05]{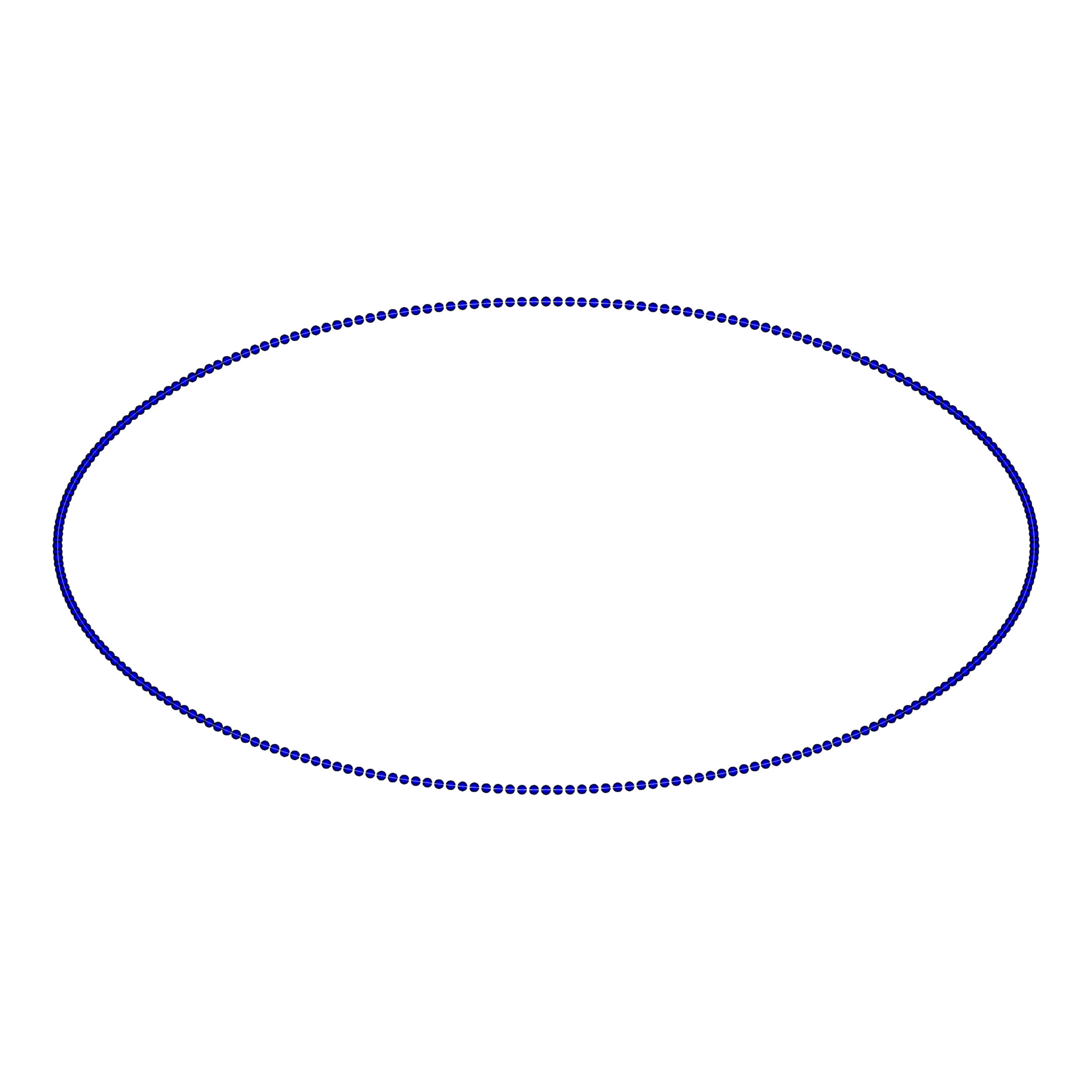}
\end{minipage}
\hspace{8pt}
\begin{minipage}[t]{0.20\linewidth}
	\centering
\includegraphics[scale=0.05]{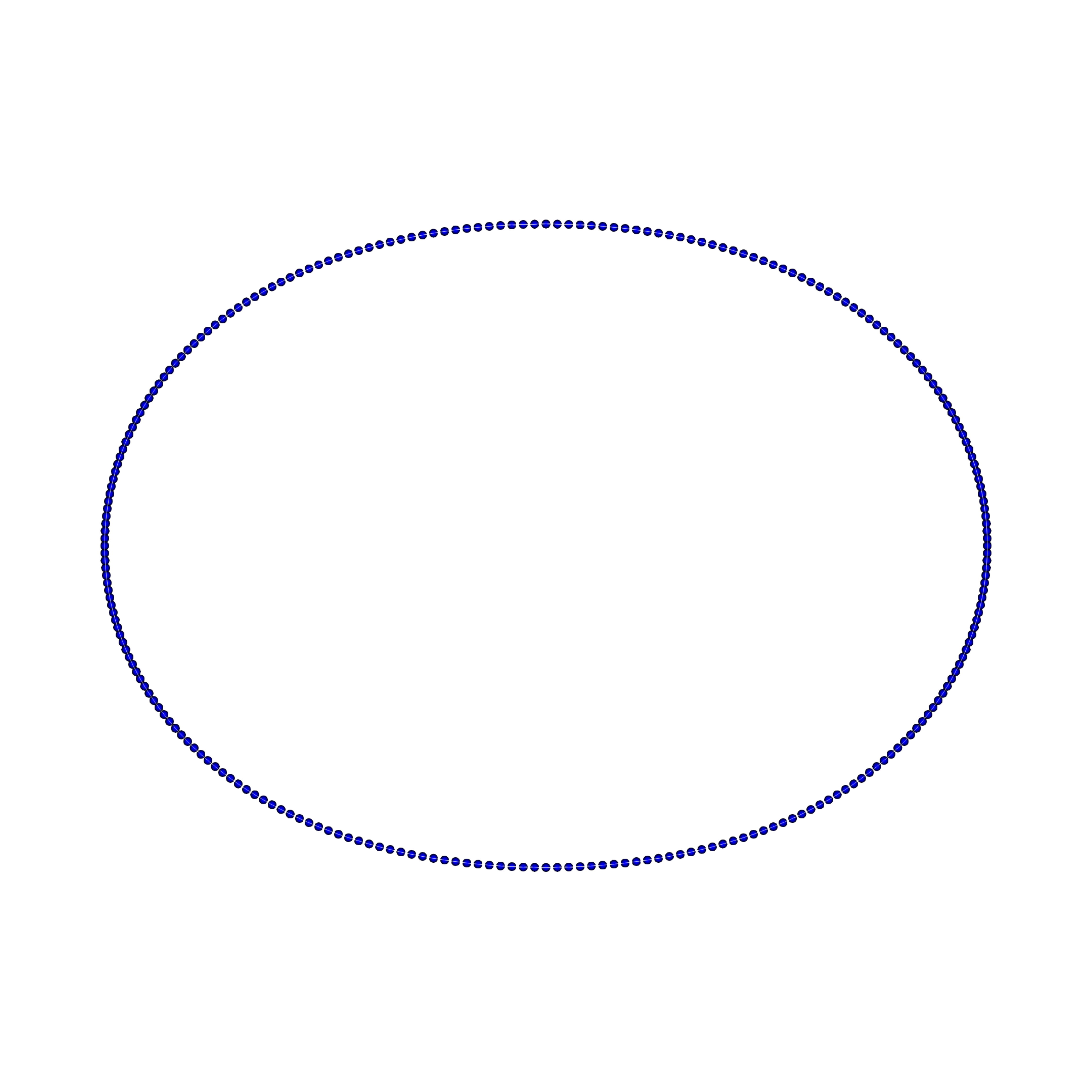}
\end{minipage}
\hspace{8pt}
\begin{minipage}[t]{0.20\linewidth}
	\centering
\includegraphics[scale=0.05]{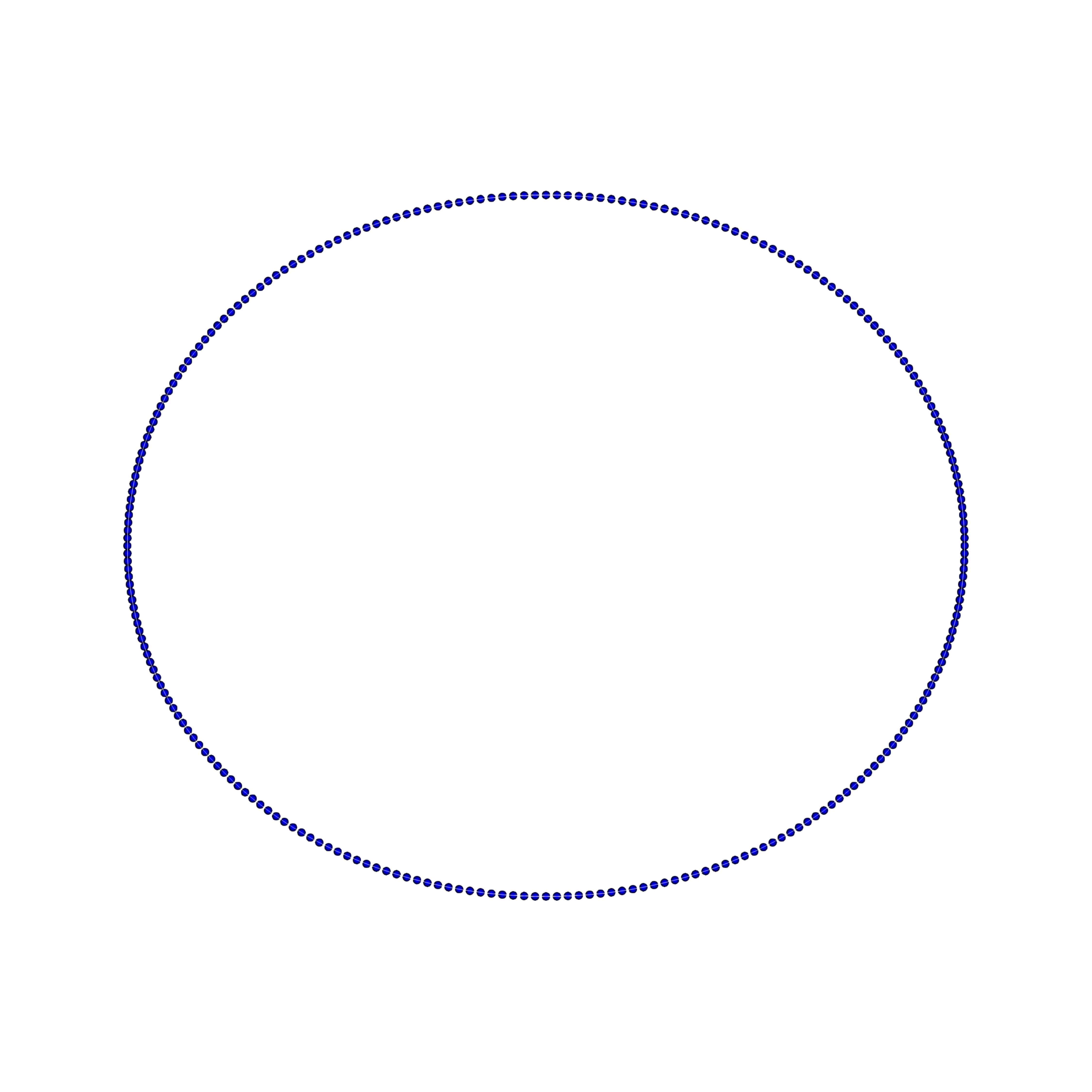}
\end{minipage}
\hspace{8pt}
\begin{minipage}[t]{0.20\linewidth}
	\centering
\includegraphics[scale=0.05]{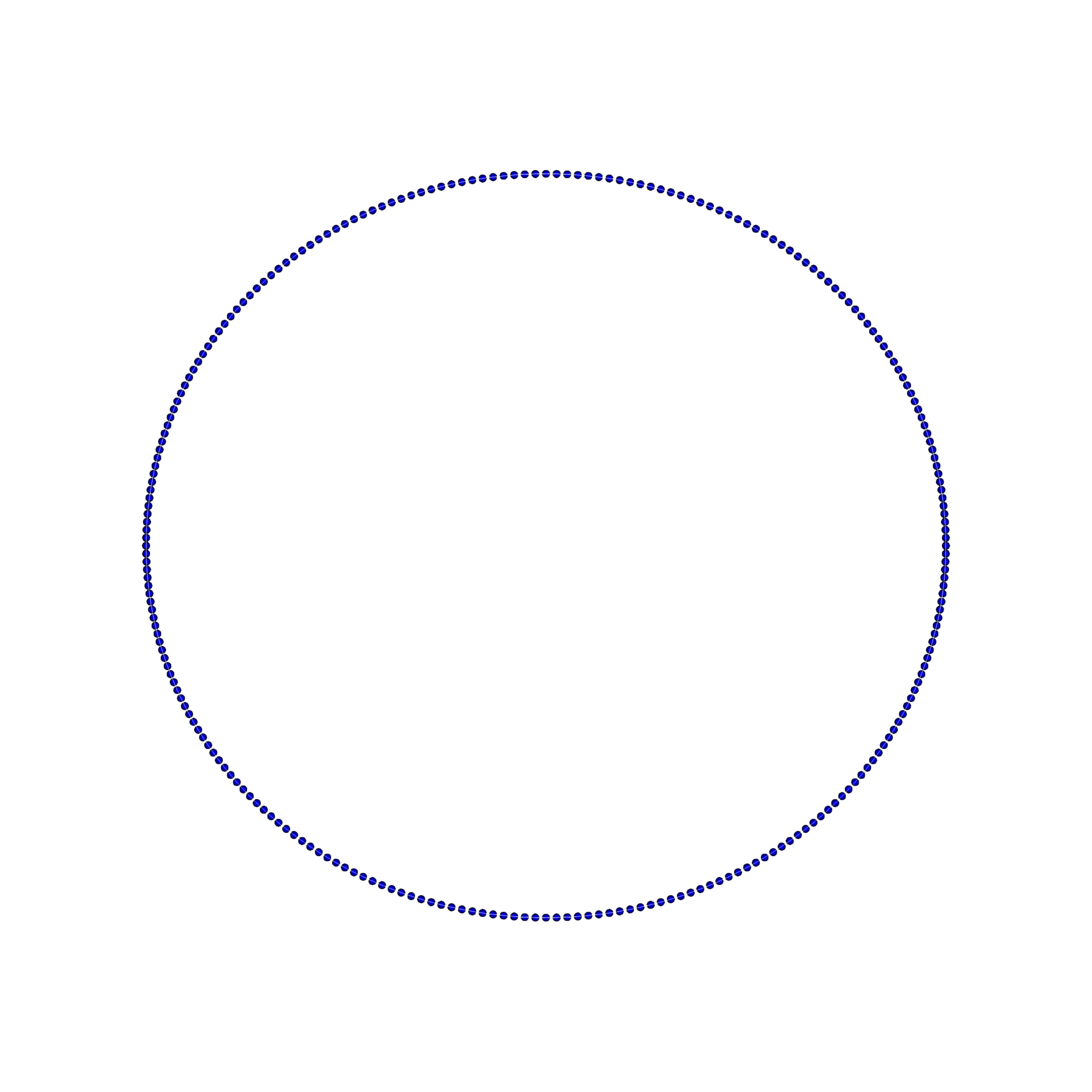}
\end{minipage}
\caption{Evolution of an ellipse by surface diffusion of curves, using $\P^2$-elements. 
We used the discretization parameters $N=128$ and $\tau=0.005$. The times are $t=0,0.3, 0.5, 0.8$, respectively.}
\label{cd_ellipse_l=2}
\end{figure}

\begin{figure}[h!]
	\centering
	\begin{subfigure}{.325\linewidth}
		\includegraphics[width=\textwidth]{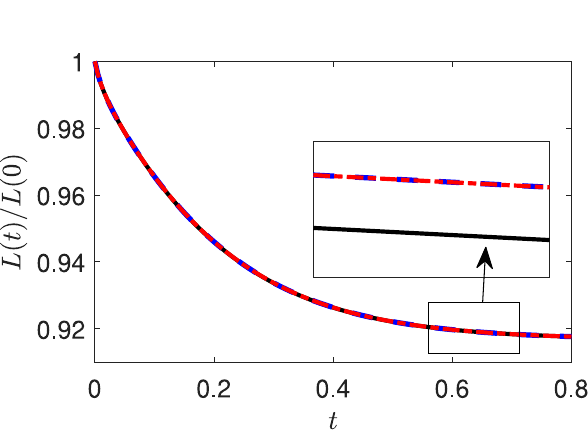}
		\caption{The normalized perimeter}
	\end{subfigure}\hfill%
	\begin{subfigure}{.325\linewidth}
		\includegraphics[width=\textwidth]{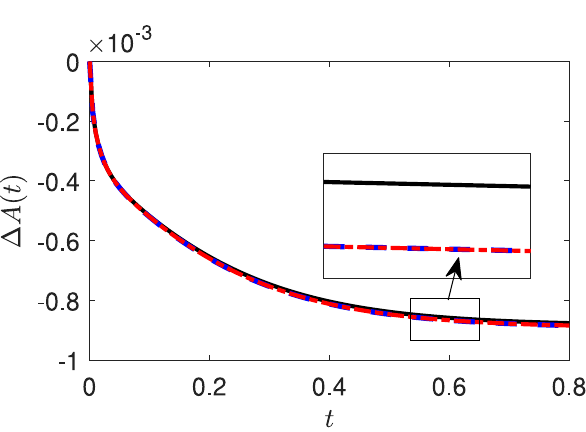}
		\caption{The relative area loss}
	\end{subfigure}\hfill%
	\begin{subfigure}{.325\linewidth}
		\includegraphics[width=\textwidth]{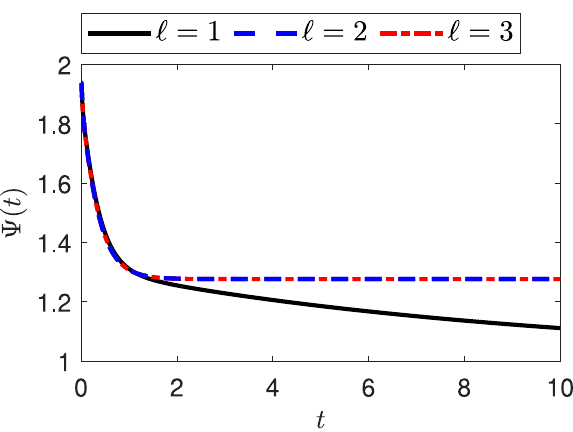}
		\caption{The mesh quality function}
	\end{subfigure}
	\caption{Evolution of geometric quantities for $\P^\ell$-elements, $\ell=1,2,3$.}
	\label{cd_ellipse_geo}
\end{figure}

\begin{figure}[h!]
\hspace{16pt}
\begin{minipage}[t]{0.20\linewidth}
	\centering
\includegraphics[scale=0.05]{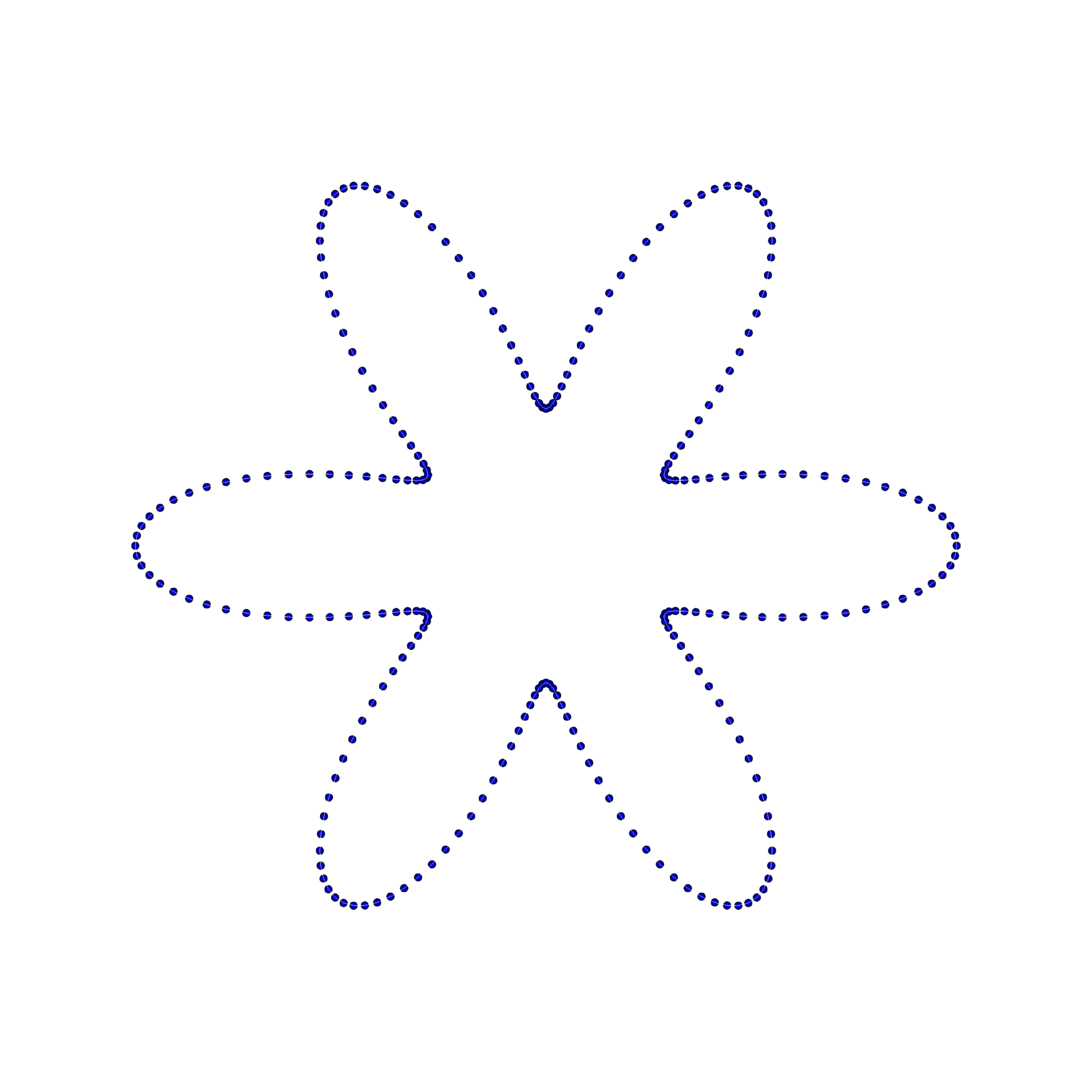}
\end{minipage}
\hspace{8pt}
\begin{minipage}[t]{0.20\linewidth}
	\centering
\includegraphics[scale=0.05]{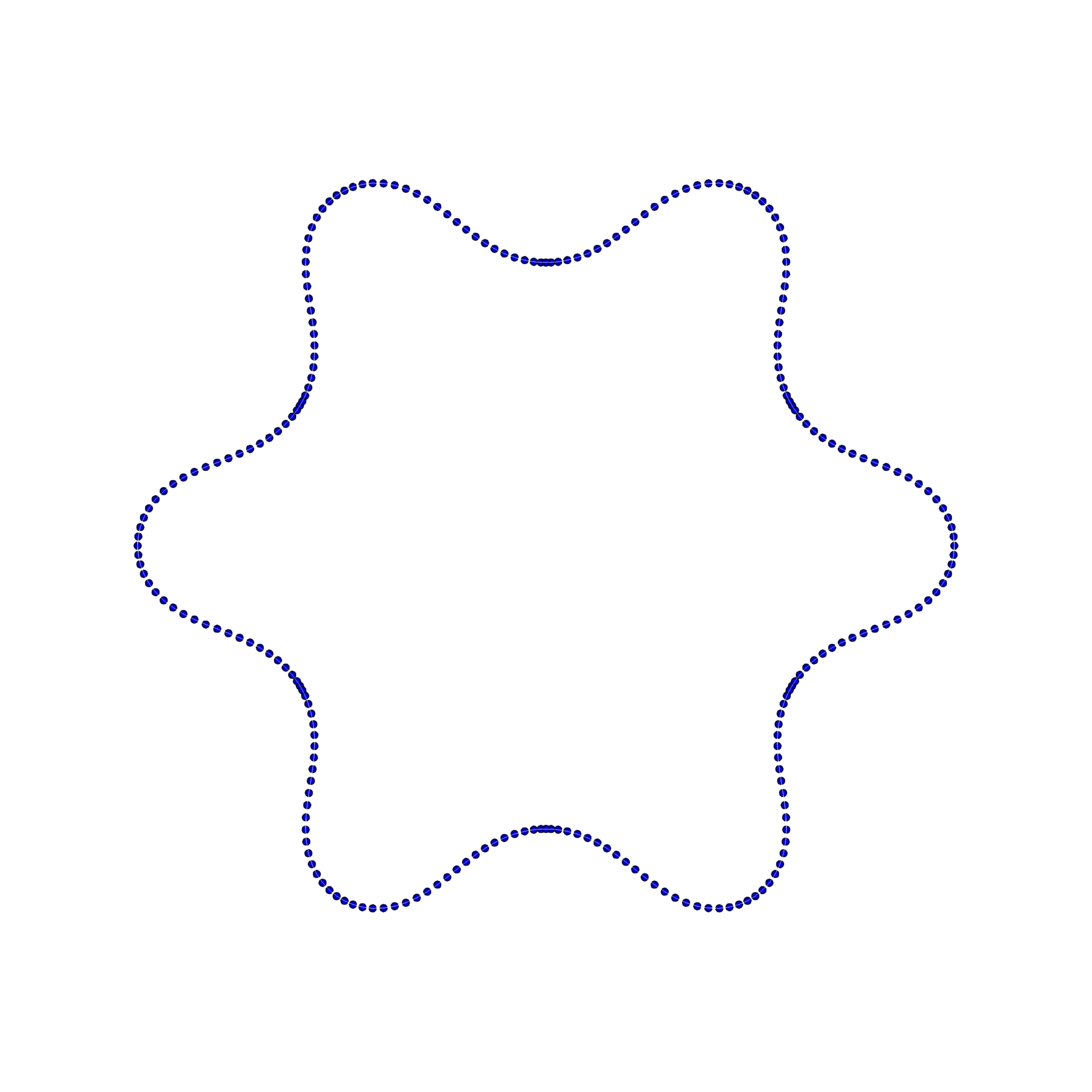}
\end{minipage}
\hspace{8pt}
\begin{minipage}[t]{0.20\linewidth}
	\centering
\includegraphics[scale=0.05]{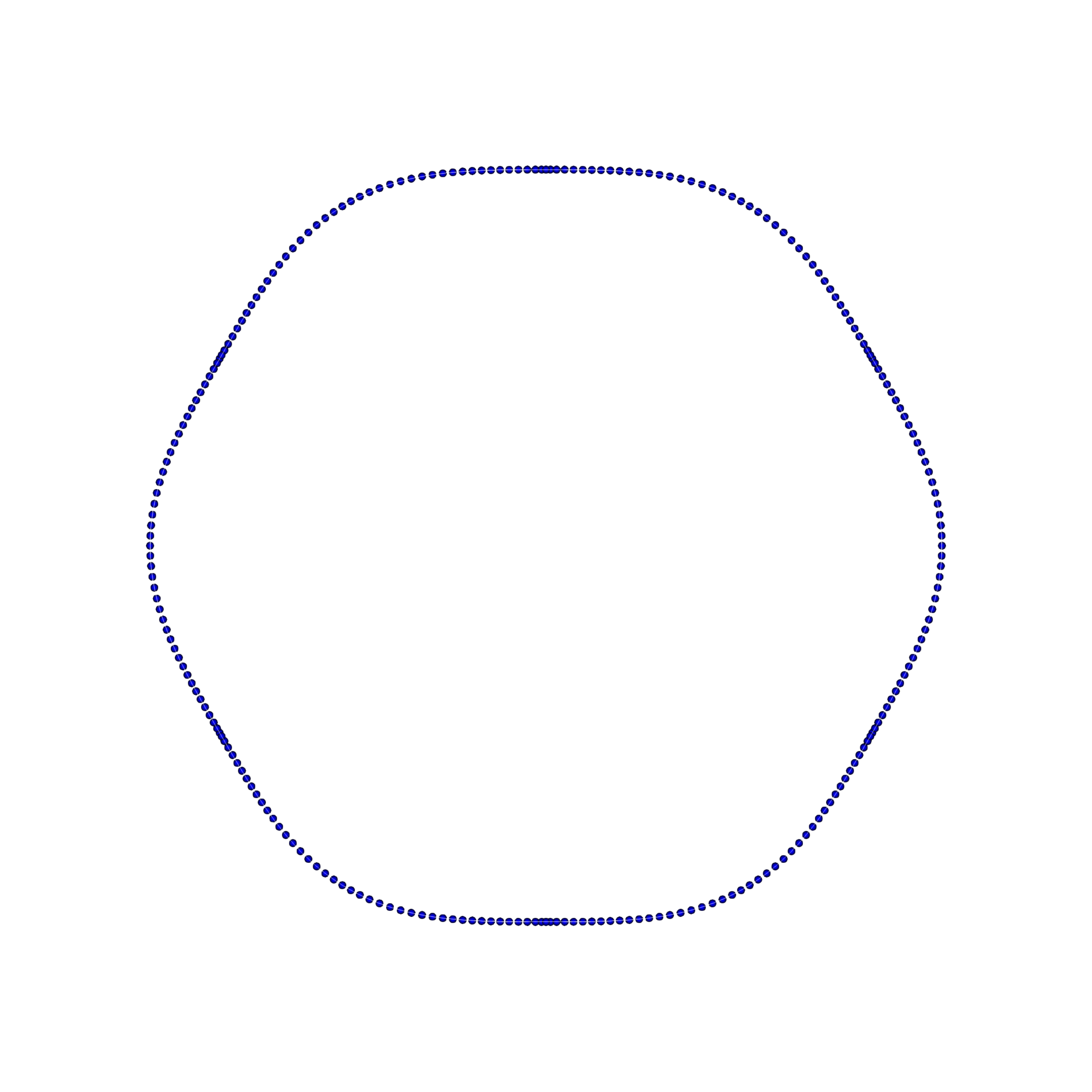}
\end{minipage}
\hspace{8pt}
\begin{minipage}[t]{0.20\linewidth}
	\centering
\includegraphics[scale=0.05]{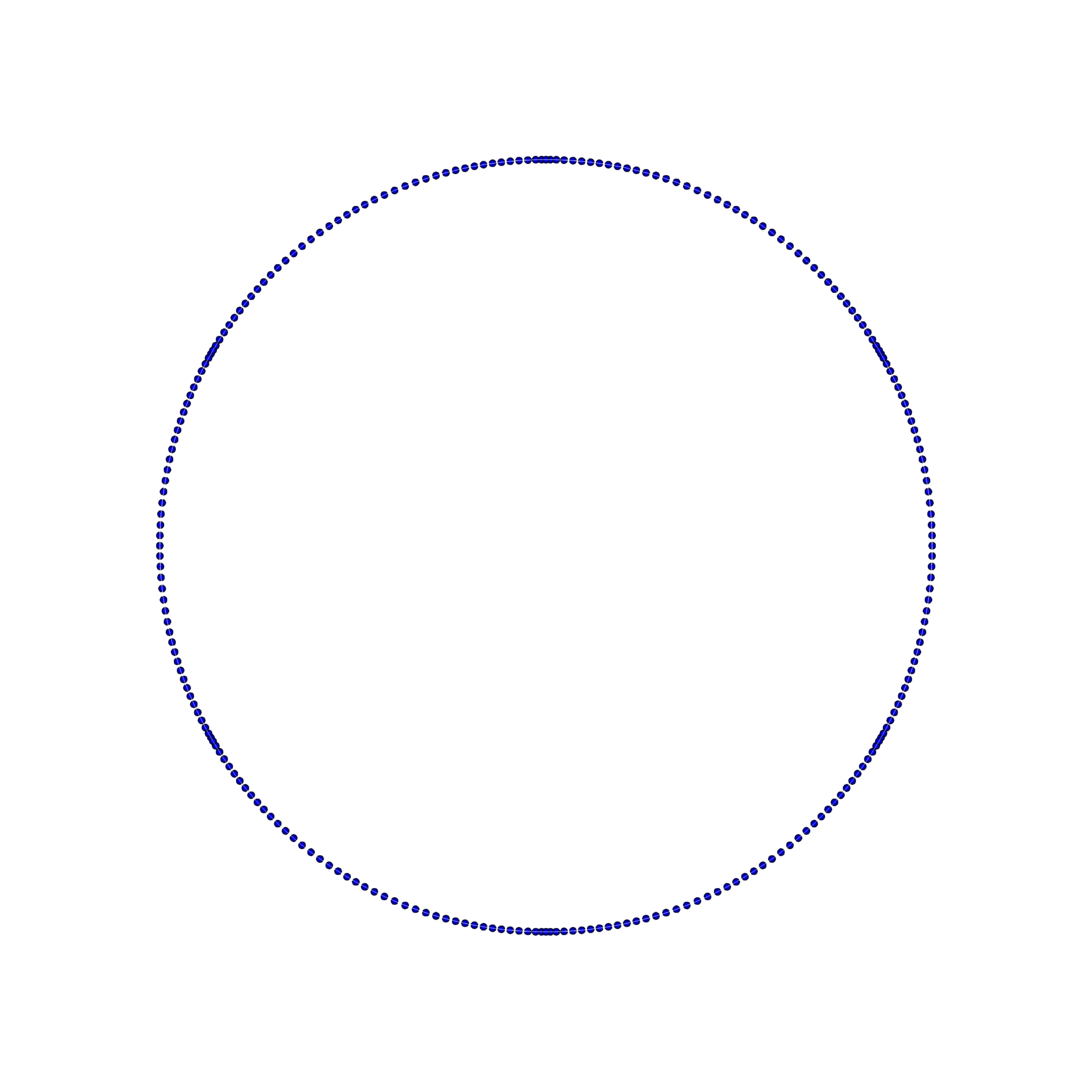}
\end{minipage}
\caption{Evolution of a flower curve by surface diffusion of curves, using $\P^2$-elements. The reference grid is an approximation polygon of an ellipse with parameter $N=128$, and the time step size is $\tau=0.001$. The times are $t=0,0.05, 0.1, 0.2$, respectively.}
\label{cd_flower_l=2}
\end{figure}

\begin{figure}[h!]
	\centering
	\begin{subfigure}{.325\linewidth}
		\includegraphics[width=\textwidth]{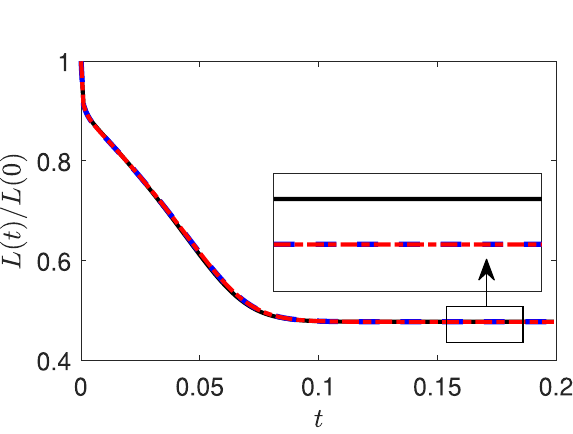}
		\caption{The normalized perimeter}
	\end{subfigure}\hfill%
	\begin{subfigure}{.325\linewidth}
		\includegraphics[width=\textwidth]{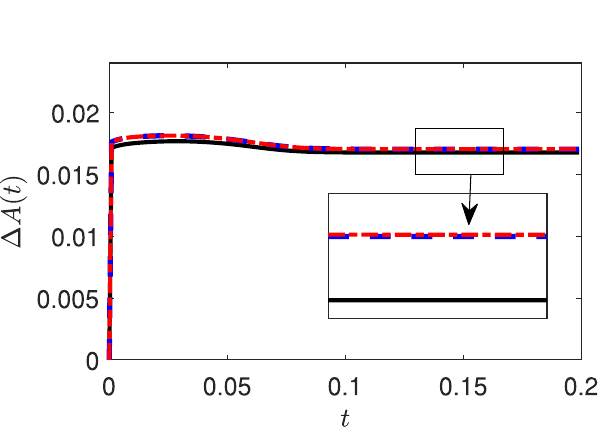}
		\caption{The relative area loss}
	\end{subfigure}\hfill%
	\begin{subfigure}{.325\linewidth}
		\includegraphics[width=\textwidth]{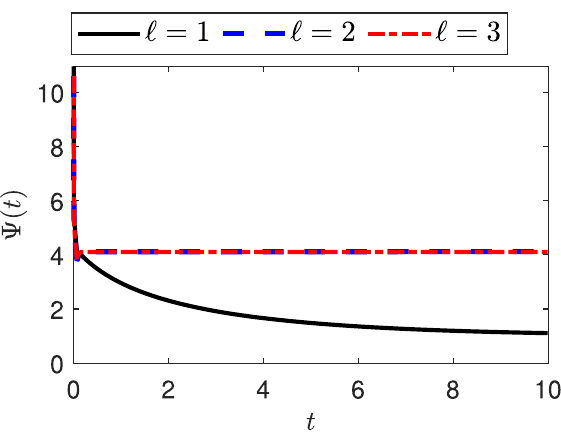}
		\caption{The mesh quality function}
	\end{subfigure}
	\caption{Evolution of geometric quantities for $\P^\ell$-elements, $\ell=1,2,3$.}
	\label{cd_flower_geo}
\end{figure}

We apply the isoparametric method for the curve case with $\ell=1,2,3$ for an ellipse and a flower curve, in each case approximating the curve
with discretizations using $N=128$ elements.
We plot the Lagrange points for the $\P^2$-elements in \Cref{cd_ellipse_l=2,cd_flower_l=2}; we observe that the isoparametric finite elements evolve an ellipse and a flower to their equilibrium, a circle. The evolutions of some geometric quantities are depicted in \Cref{cd_ellipse_geo,cd_flower_geo}. We can observe that for all $\ell=1,2$ and $3$, the isoparametric finite element methods decrease the perimeter, and the relative area loss is relatively small. In contrast to the linear finite element method $\ell=1$ case, it appears that the mesh quality of higher order BGN-type methods may not decrease to $1$ for long time evolutions.

\smallskip

We end this subsection by employing the isoparametric finite element methods for two surfaces with different topologies, i.e., an ellipsoid and a torus. The approximating surfaces employ the number of elements and vertices $(J,K)=(676,340)$ for the ellipsoid, and $(J,K)=(720,360)$ for the torus.

\begin{figure}[ht]
\hspace{-14pt}
\begin{minipage}[t]{0.20\linewidth}
	\centering
\includegraphics[scale=0.06]{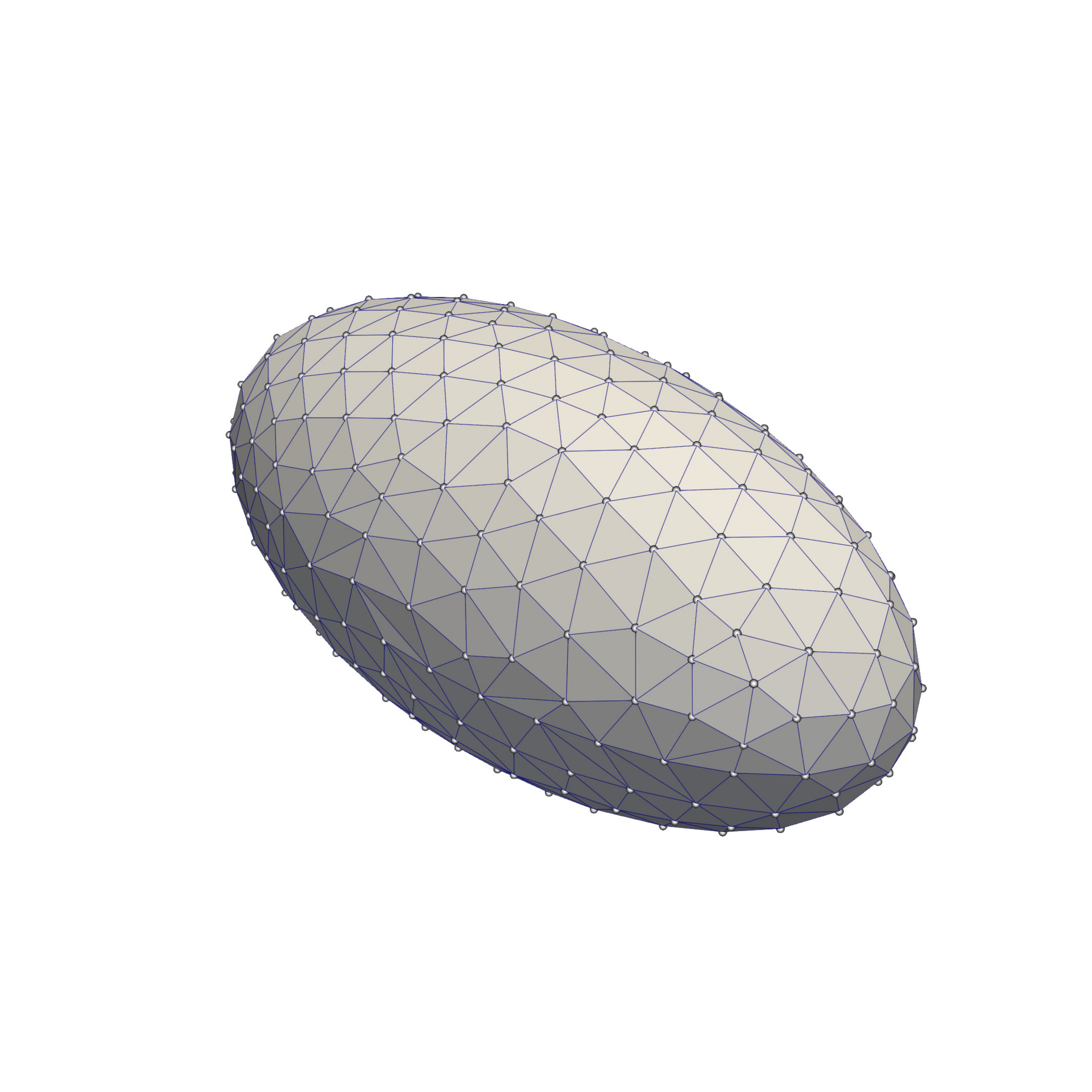}
\end{minipage}
\hspace{12pt}
\begin{minipage}[t]{0.20\linewidth}
	\centering
\includegraphics[scale=0.06]{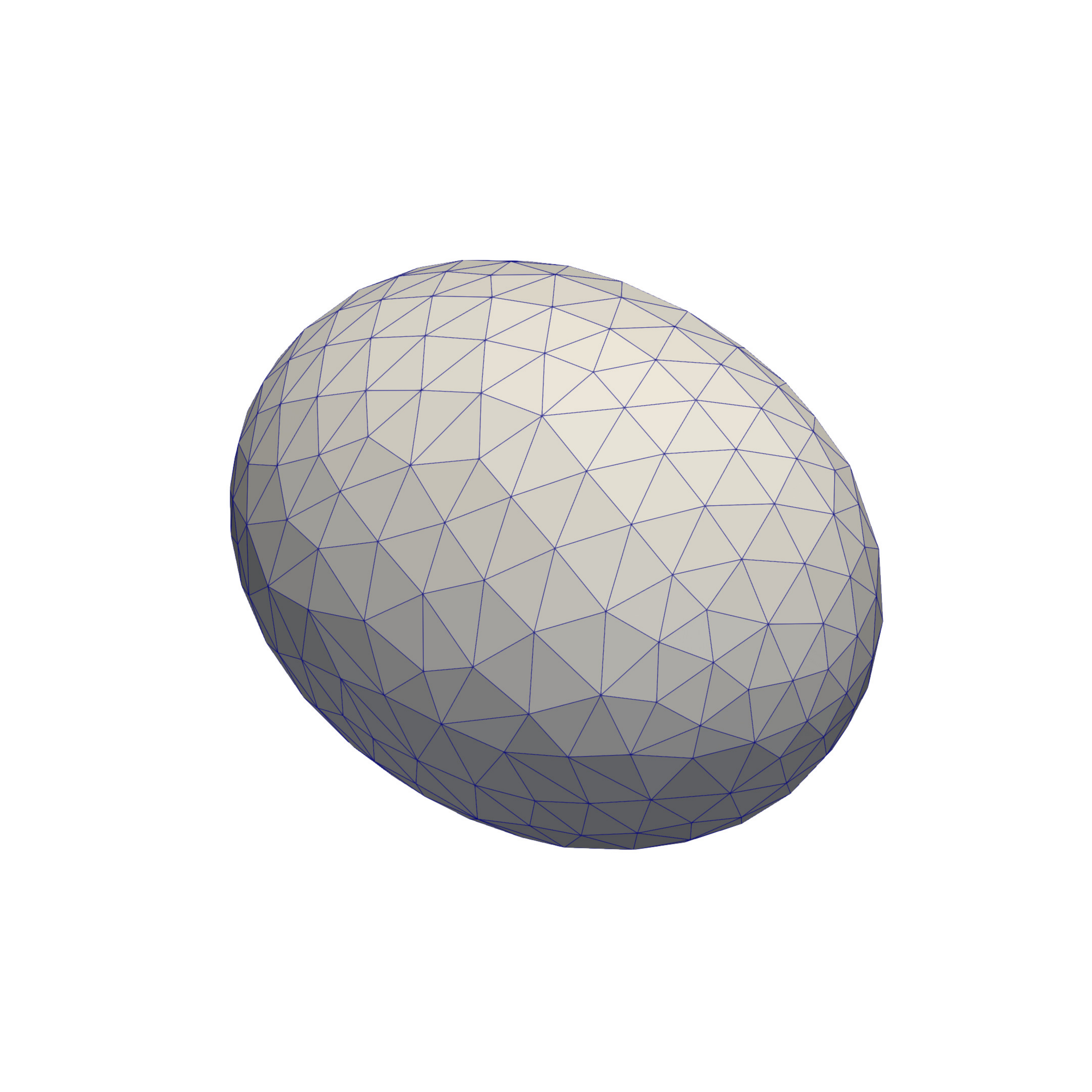}
\end{minipage}
\hspace{12pt}
\begin{minipage}[t]{0.20\linewidth}
	\centering
\includegraphics[scale=0.06]{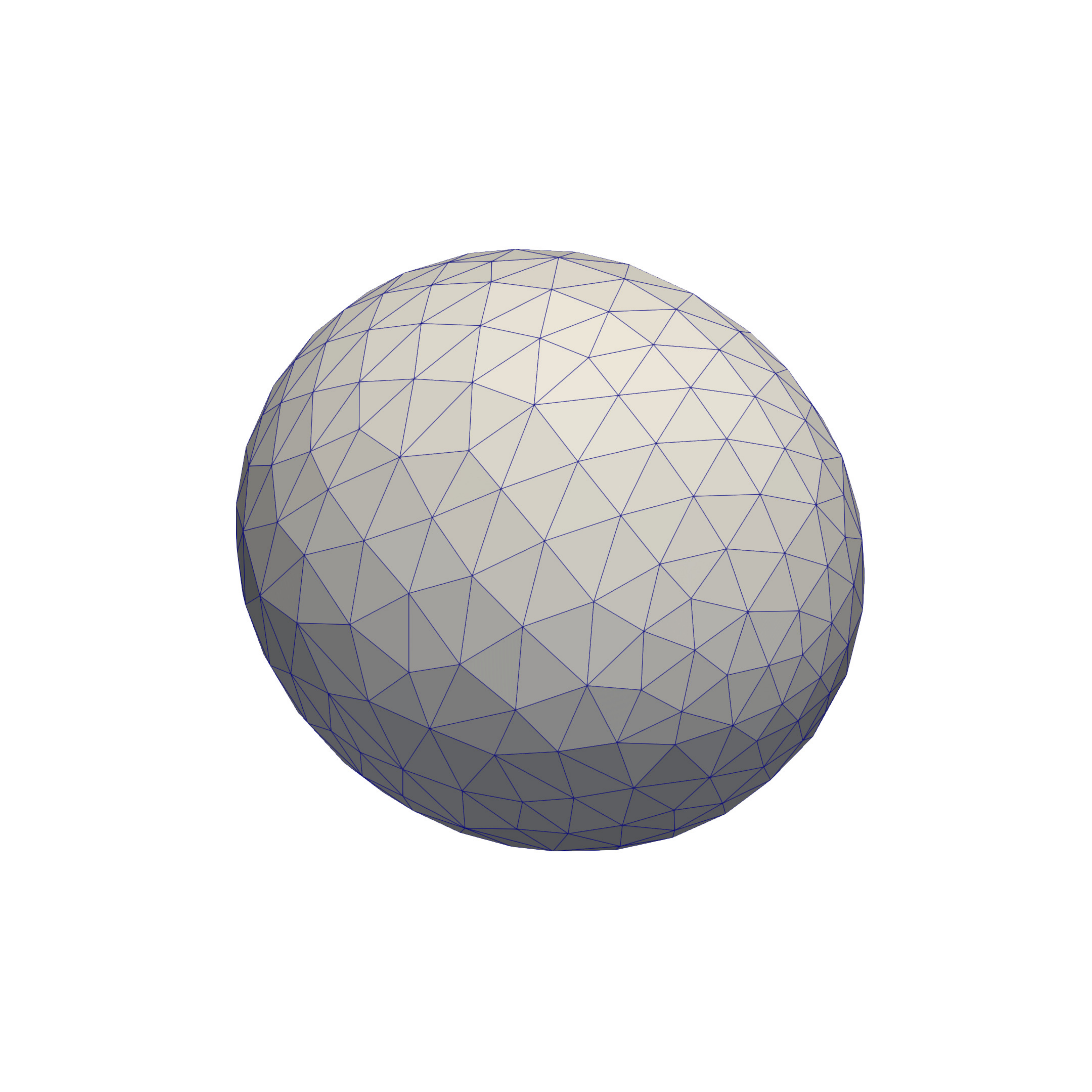}
\end{minipage}
\hspace{12pt}
\begin{minipage}[t]{0.20\linewidth}
	\centering
\includegraphics[scale=0.06]{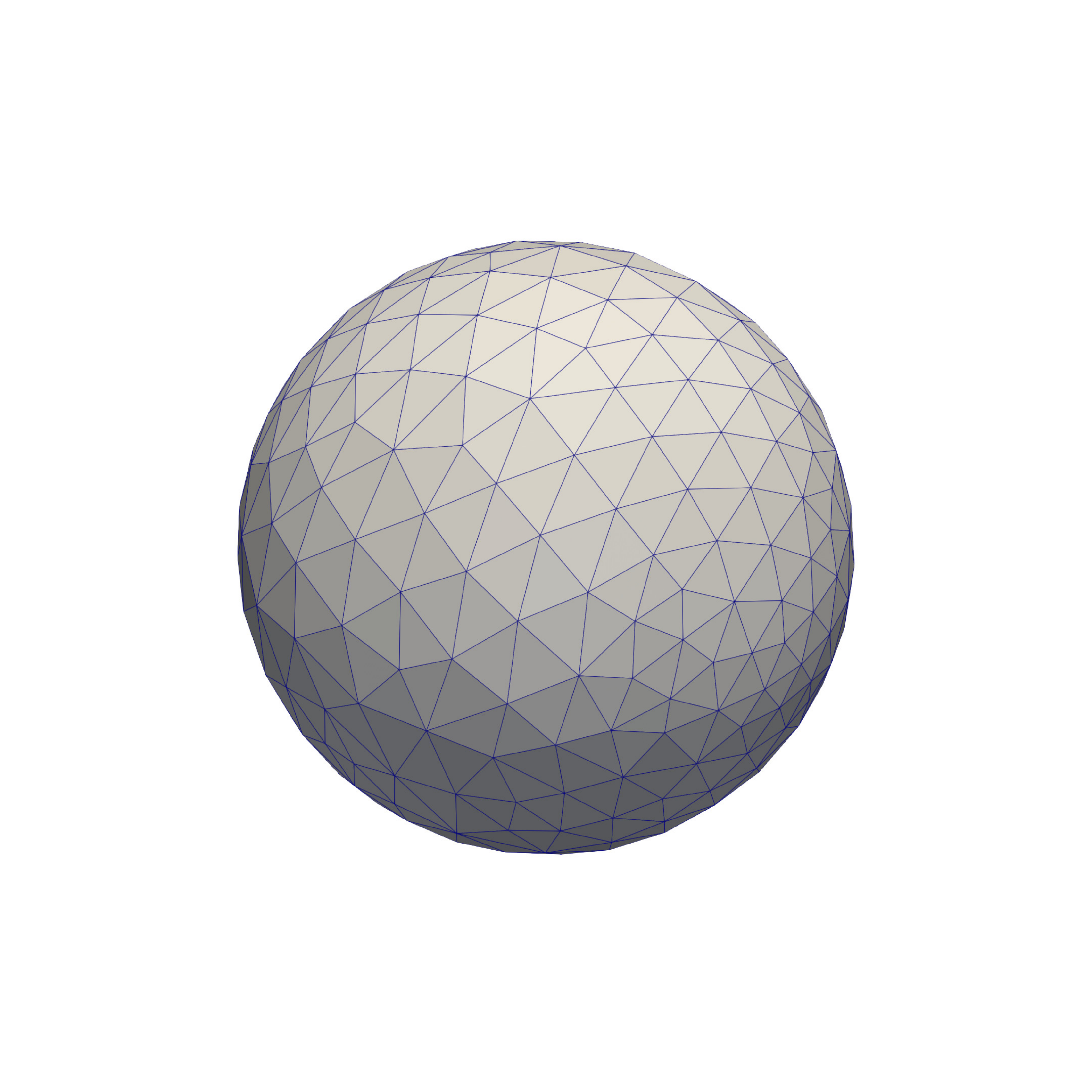}
\end{minipage}
\caption{Evolution of an ellipsoid by surface diffusion, using $\P^1$-elements. 
We used the discretization parameters $(J,K)=(676,340)$ and $\tau=0.005$. The number of total Lagrange nodes is $340$. The times are  $t=0,0.15, 0.25, 0.5$, respectively.}
\label{sd_ellipsoid_l=1}
\end{figure}

\begin{figure}[h!]
\hspace{-14pt}
\begin{minipage}[t]{0.20\linewidth}
	\centering
\includegraphics[scale=0.06]{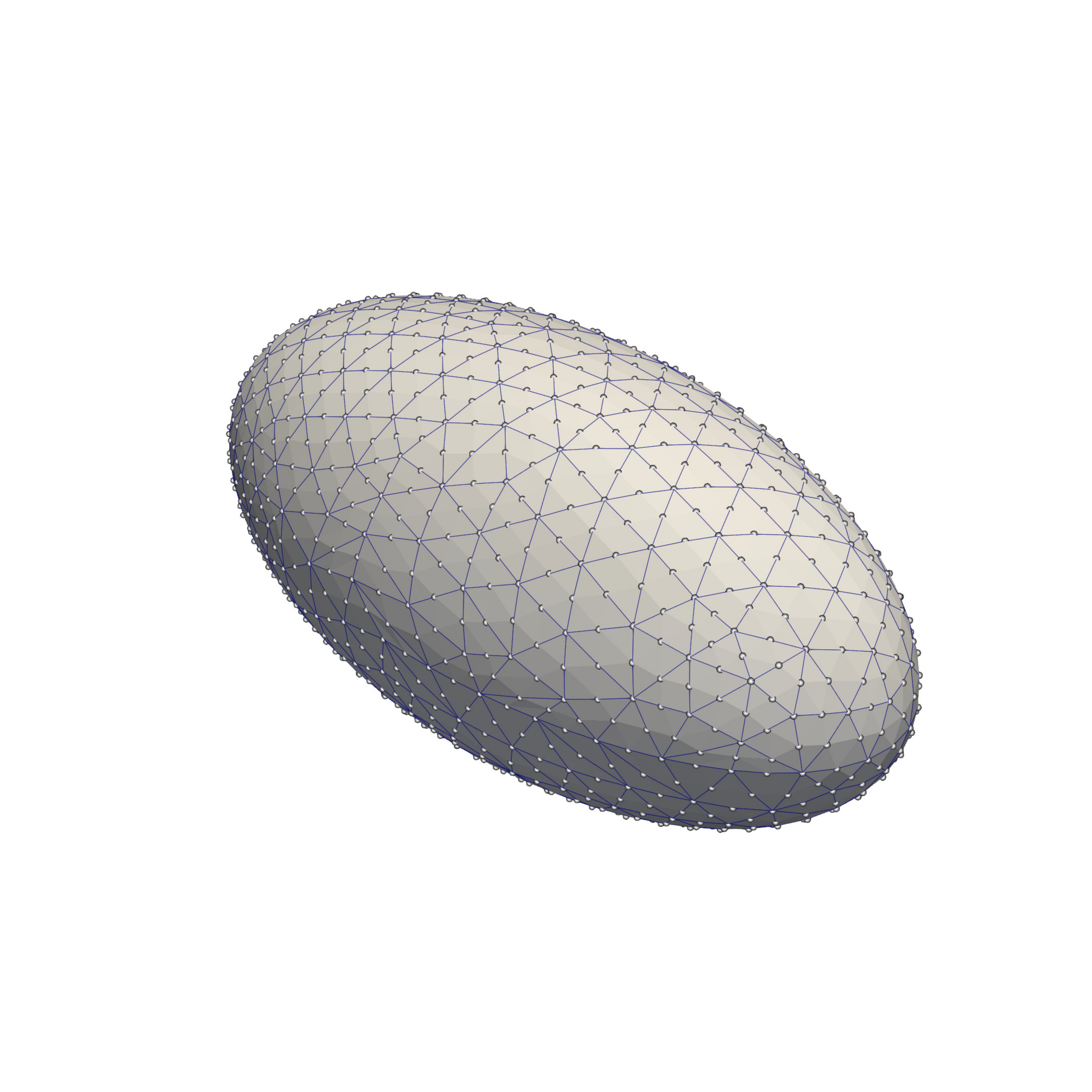}
\end{minipage}
\hspace{12pt}
\begin{minipage}[t]{0.20\linewidth}
	\centering
\includegraphics[scale=0.06]{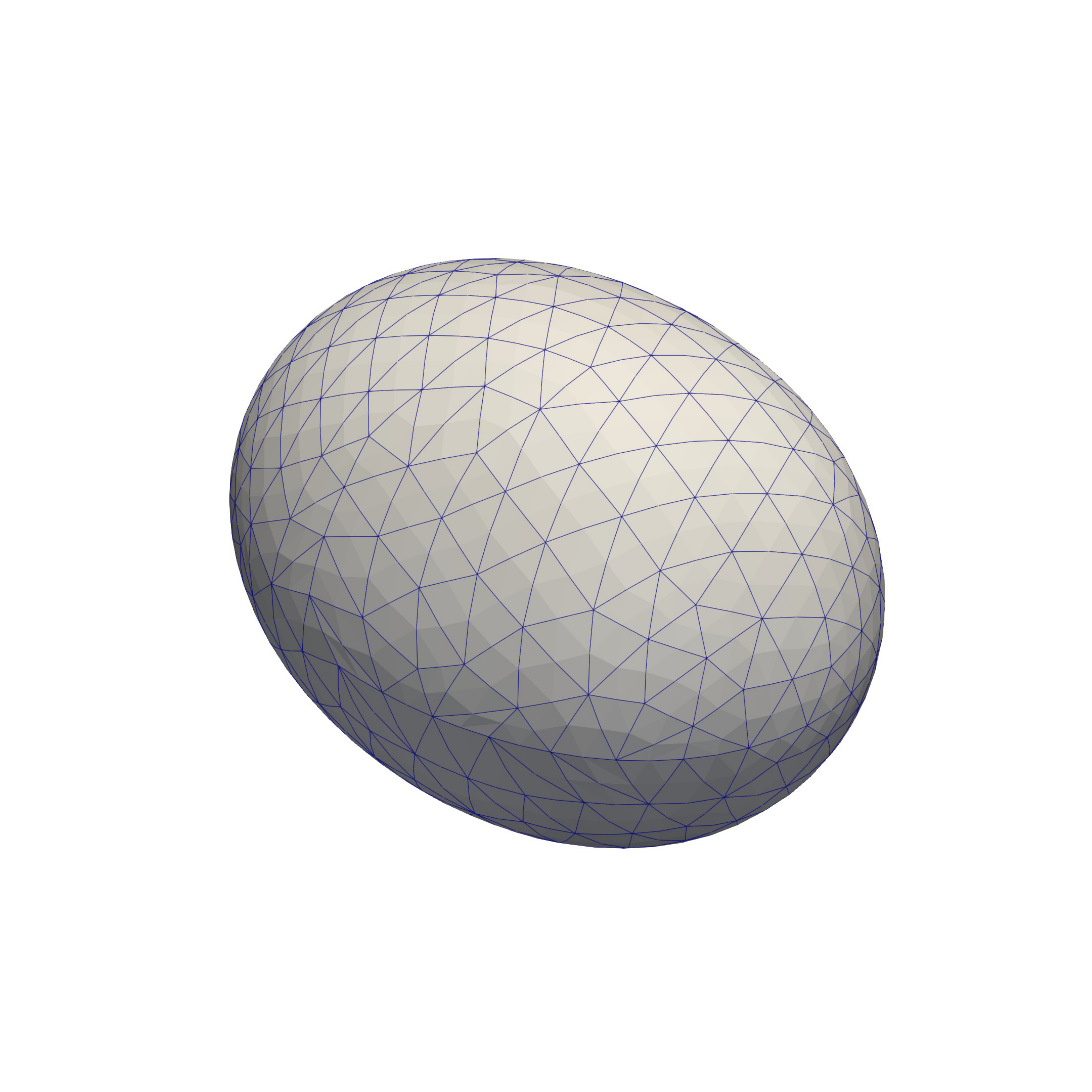}
\end{minipage}
\hspace{12pt}
\begin{minipage}[t]{0.20\linewidth}
	\centering
\includegraphics[scale=0.06]{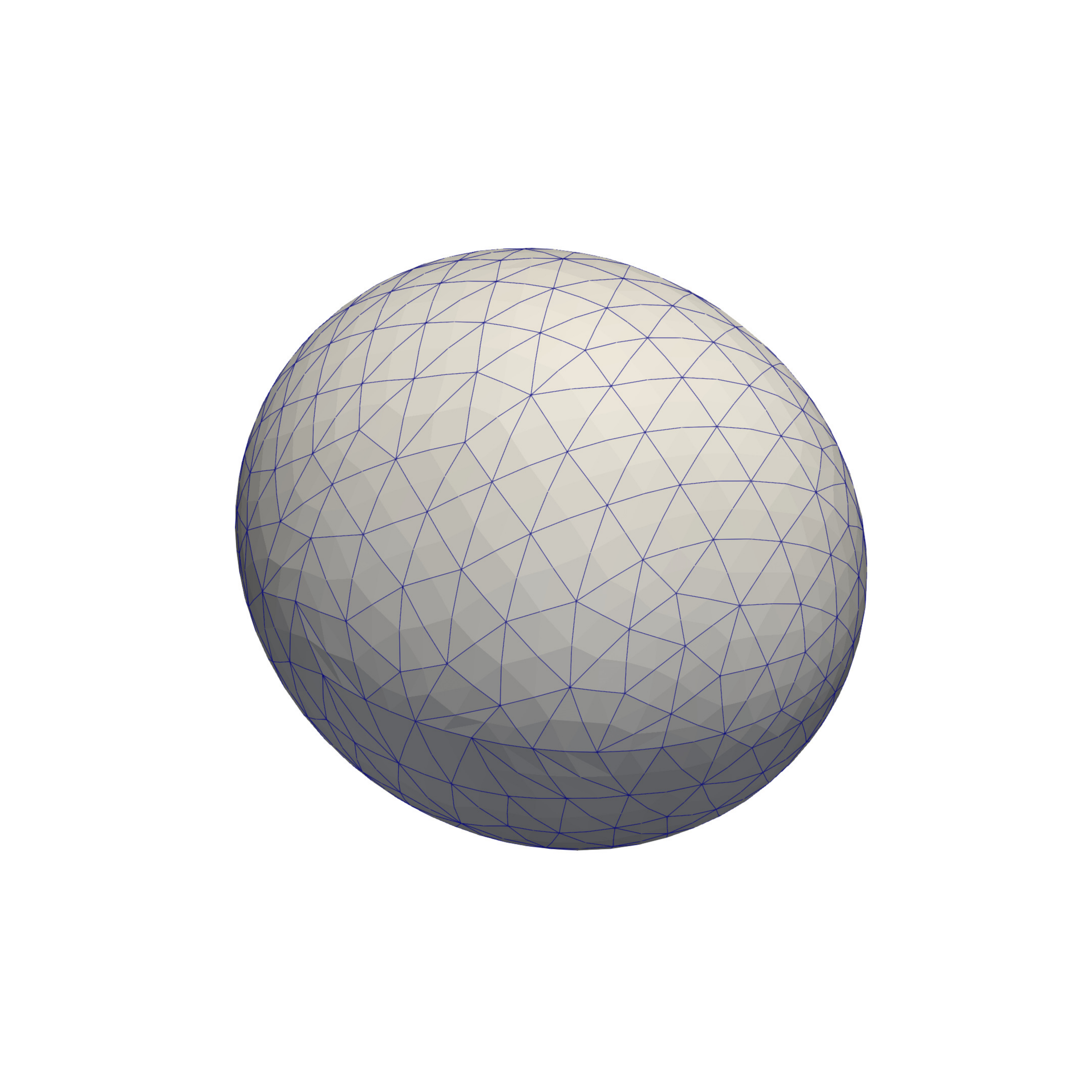}
\end{minipage}
\hspace{12pt}
\begin{minipage}[t]{0.20\linewidth}
	\centering
\includegraphics[scale=0.06]{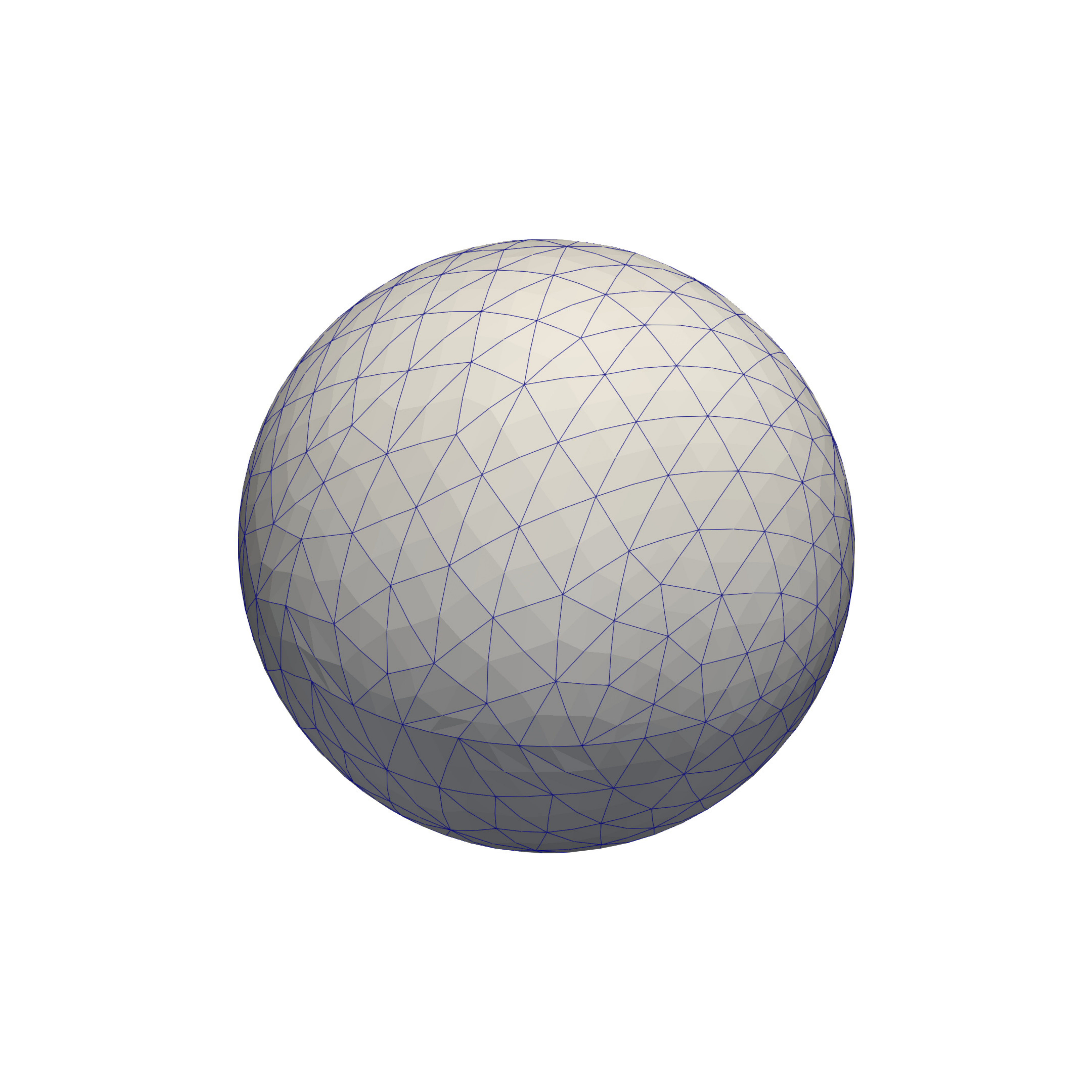}
\end{minipage}
\caption{Evolution of an ellipsoid by surface diffusion, using $\P^2$-elements. The number of total Lagrange nodes is $1354$.}
\label{sd_ellipsoid_l=2}
\end{figure}

\begin{figure}[h!]
	\centering
	\begin{subfigure}{.325\linewidth}
		\includegraphics[width=\textwidth]{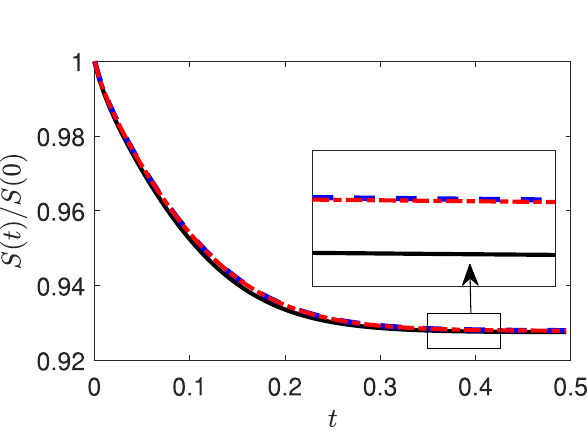}
		\caption{The normalized surface area}
	\end{subfigure}\hfill%
	\begin{subfigure}{.325\linewidth}
		\includegraphics[width=\textwidth]{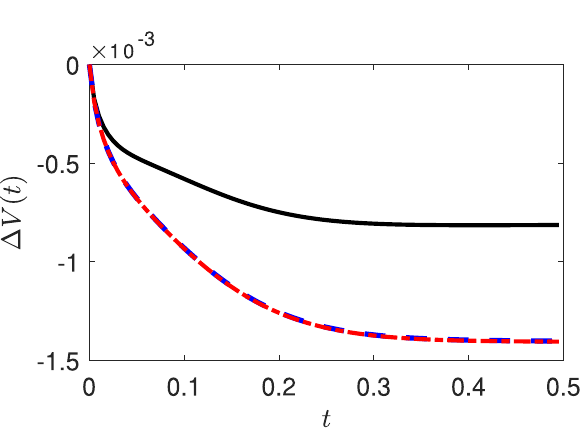}
		\caption{The relative volume loss}
	\end{subfigure}\hfill%
	\begin{subfigure}{.325\linewidth}
		\includegraphics[width=\textwidth]{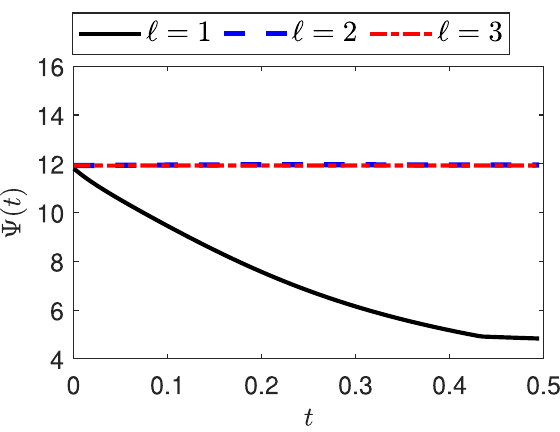}
		\caption{The mesh quality function}
	\end{subfigure}
	\caption{Evolution of geometric quantities for $\P^\ell$-elements, $\ell=1,2,3$.}
	\label{sd_ellipsoid_geo}
\end{figure}

As depicted in \Cref{sd_ellipsoid_l=1,sd_ellipsoid_l=2}, surface diffusion evolves the ellipsoid to its equilibrium, a sphere. For a fair comparison, we further show the positions of the Lagrange nodes for $\ell=1,2$, as demonstrated in the first snapshot of \Cref{sd_ellipsoid_l=1,sd_ellipsoid_l=2}. Observe that the evolution process only needs to control the motion of these Lagrange nodes. The $\P^1$-element method uses Lagrange nodes at grid vertices, resulting in a total number of Lagrange nodes of 340, while the $\P^2$-element method uses Lagrange nodes at vertices and edges, resulting in 1354 nodes in total. Moreover, in \Cref{sd_ellipsoid_geo}, we monitor the evolution of some geometric quantities. These numerical results show that our proposed isoparametric methods successfully decrease the surface area and effectively preserve the volume during the evolution while maintaining good mesh quality.

\begin{figure}[h!]
\hspace{-12pt}
\begin{minipage}[t]{0.20\linewidth}
	\centering
\includegraphics[scale=0.06]{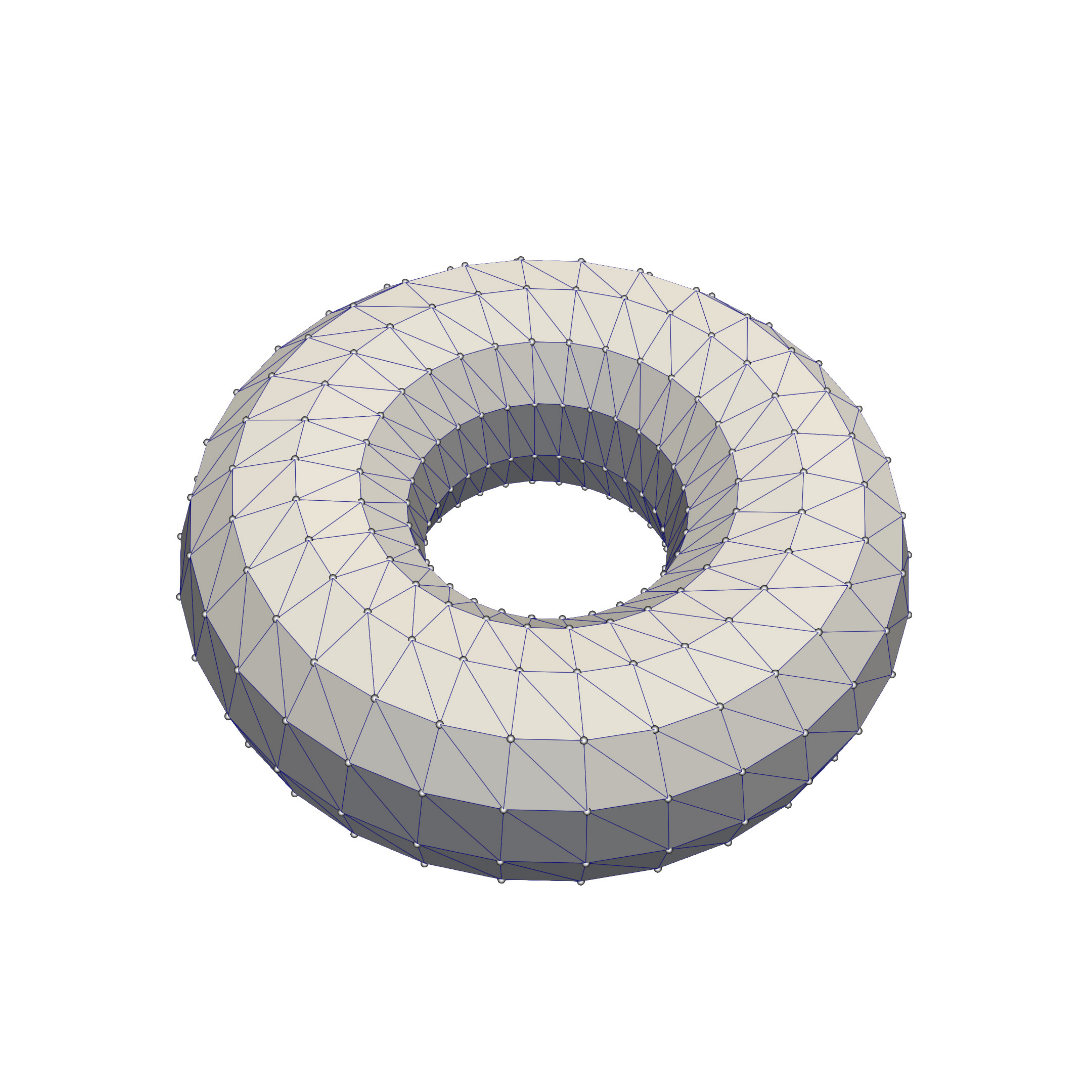}
\end{minipage}
\hspace{12pt}
\begin{minipage}[t]{0.20\linewidth}
	\centering
\includegraphics[scale=0.06]{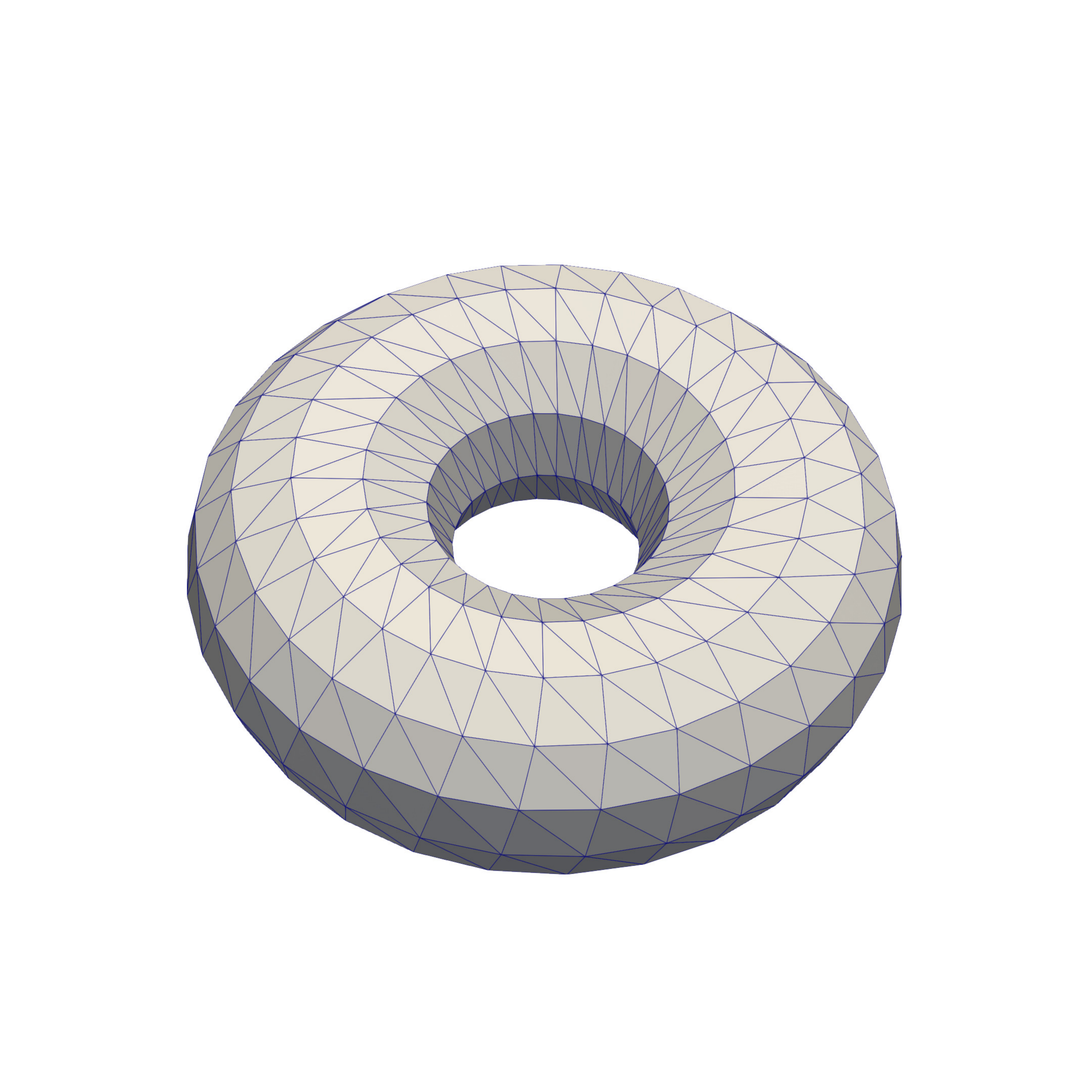}
\end{minipage}
\hspace{12pt}
\begin{minipage}[t]{0.20\linewidth}
	\centering
\includegraphics[scale=0.06]{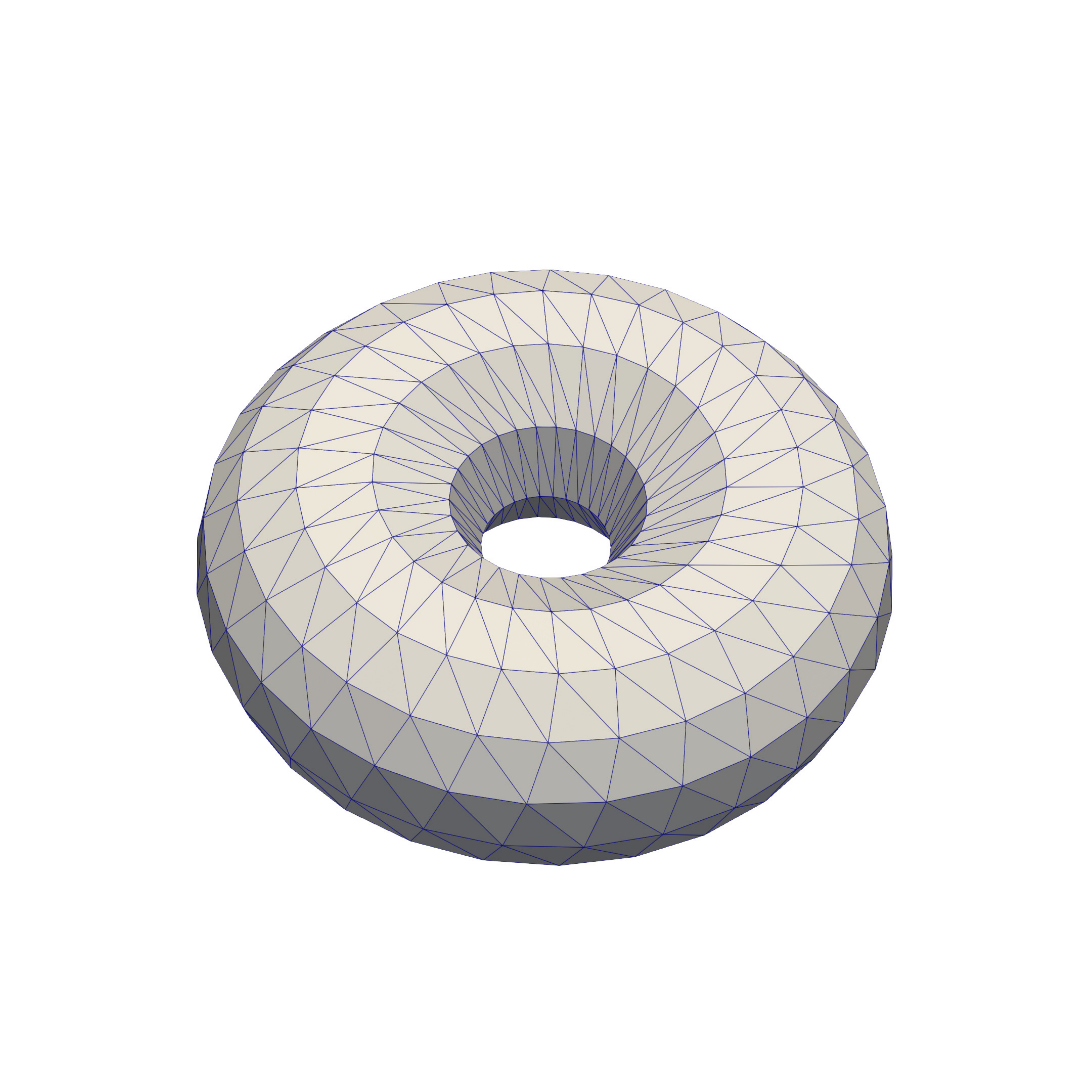}
\end{minipage}
\hspace{12pt}
\begin{minipage}[t]{0.20\linewidth}
	\centering
\includegraphics[scale=0.06]{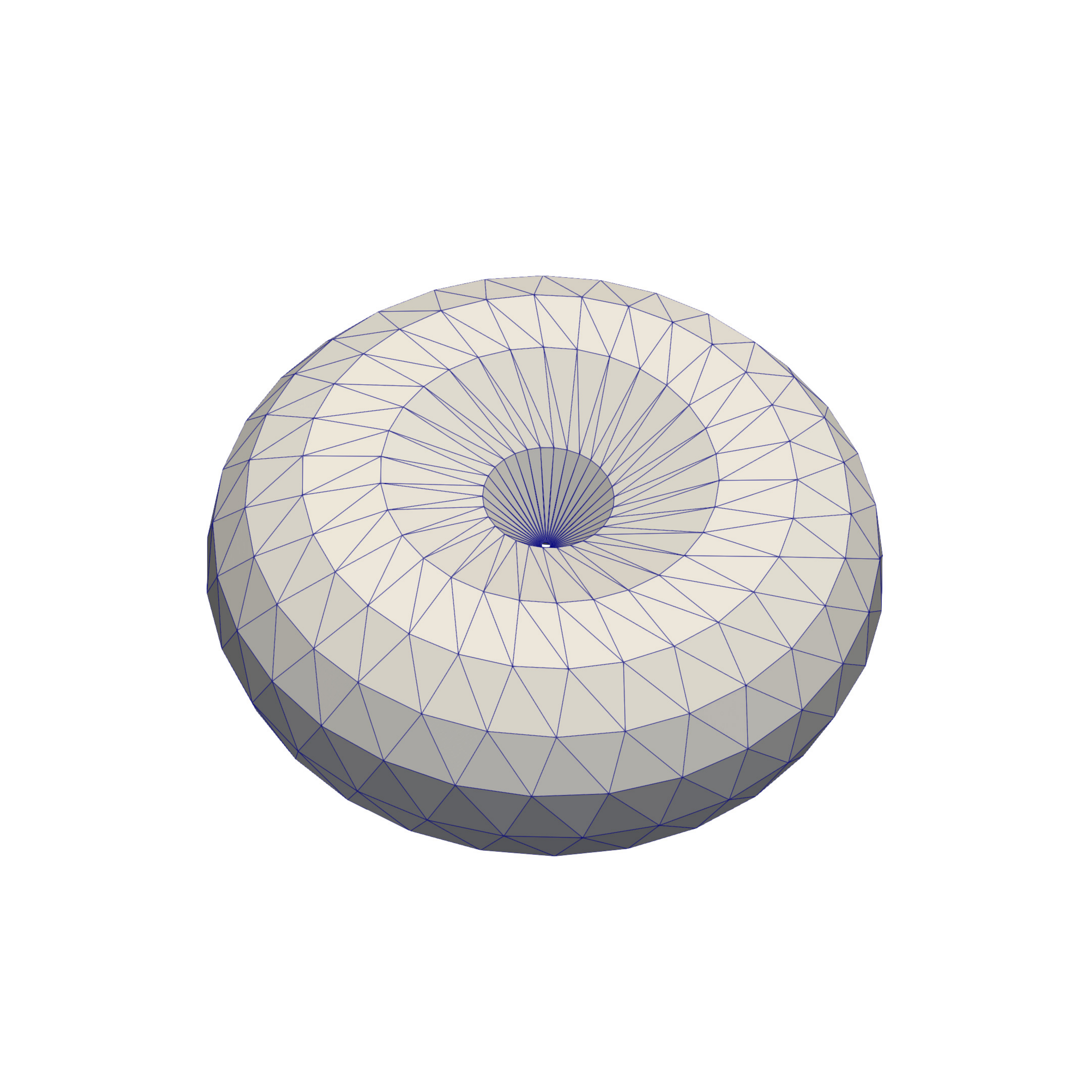}
\end{minipage}
\caption{Evolution of a torus by surface diffusion, using $\P^1$-elements. 
We used the discretization parameters $(J,K)=(720,360)$ and $\tau=0.005$. The total number of Lagrange nodes is $360$. The times are  $t=0,0.1, 0.2, 0.3$, respectively.}
\label{sd_torus_l=1}
\end{figure}

\begin{figure}[h!]
\hspace{-12pt}
\begin{minipage}[t]{0.20\linewidth}
	\centering
\includegraphics[scale=0.06]{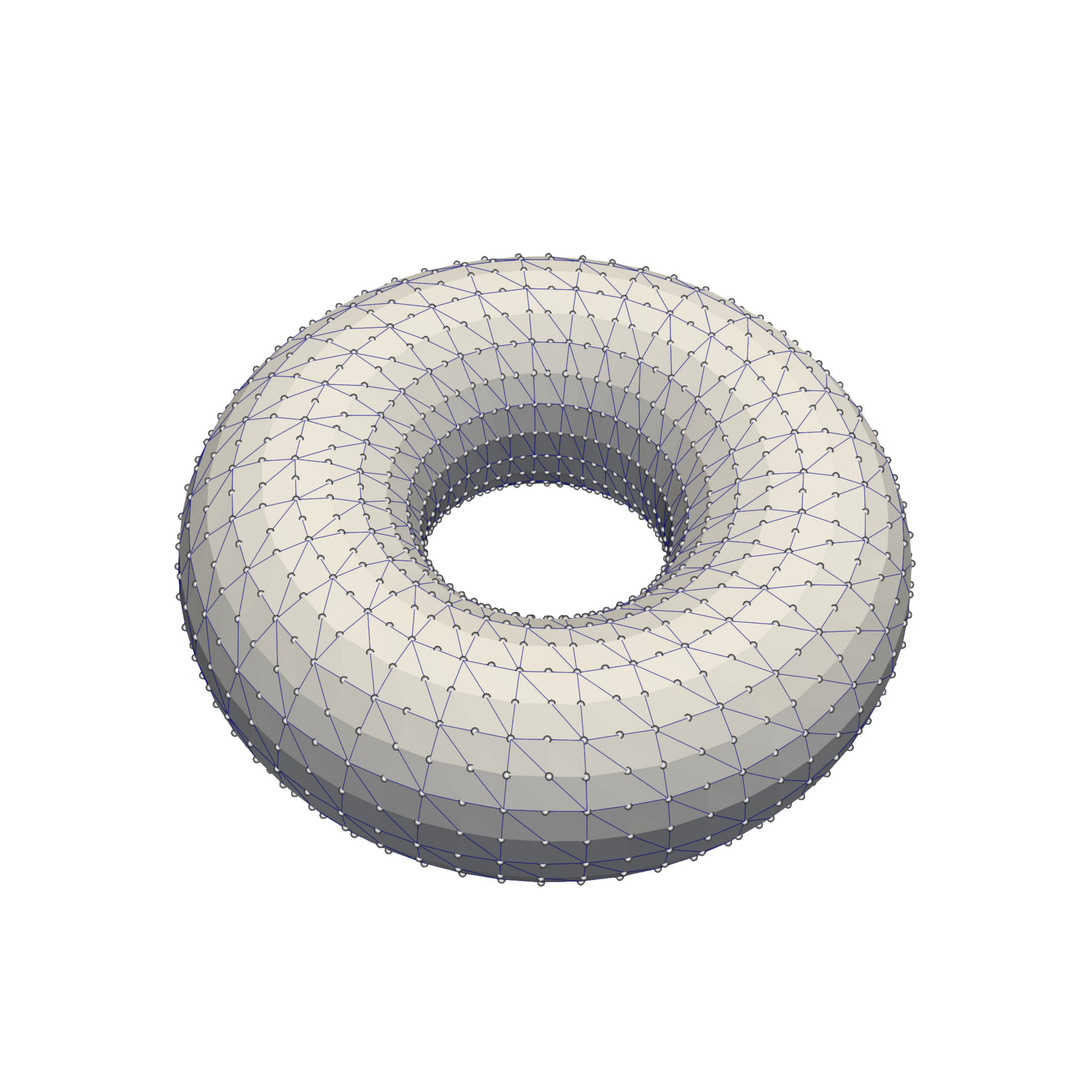}
\end{minipage}
\hspace{12pt}
\begin{minipage}[t]{0.20\linewidth}
	\centering
\includegraphics[scale=0.06]{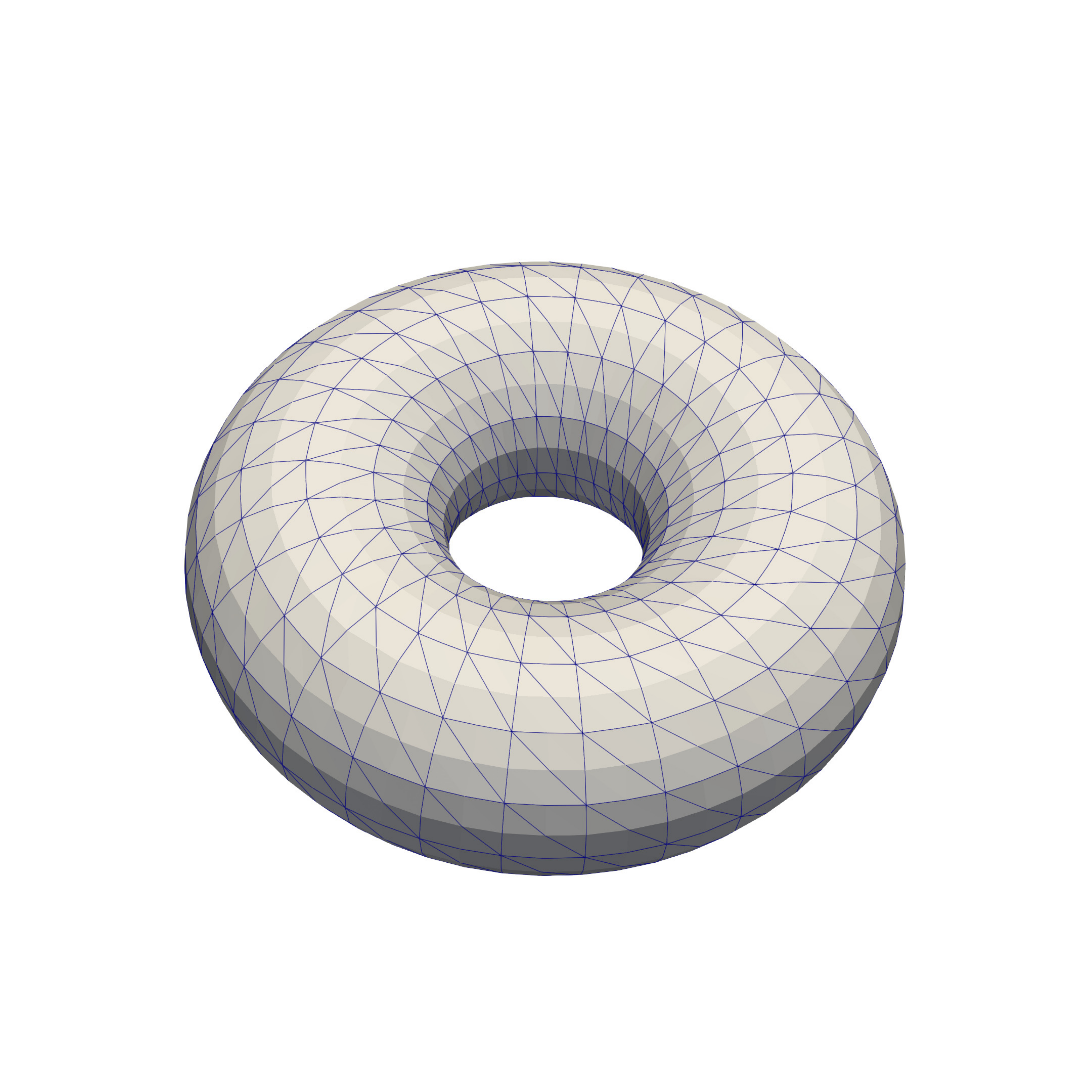}
\end{minipage}
\hspace{12pt}
\begin{minipage}[t]{0.20\linewidth}
	\centering
\includegraphics[scale=0.06]{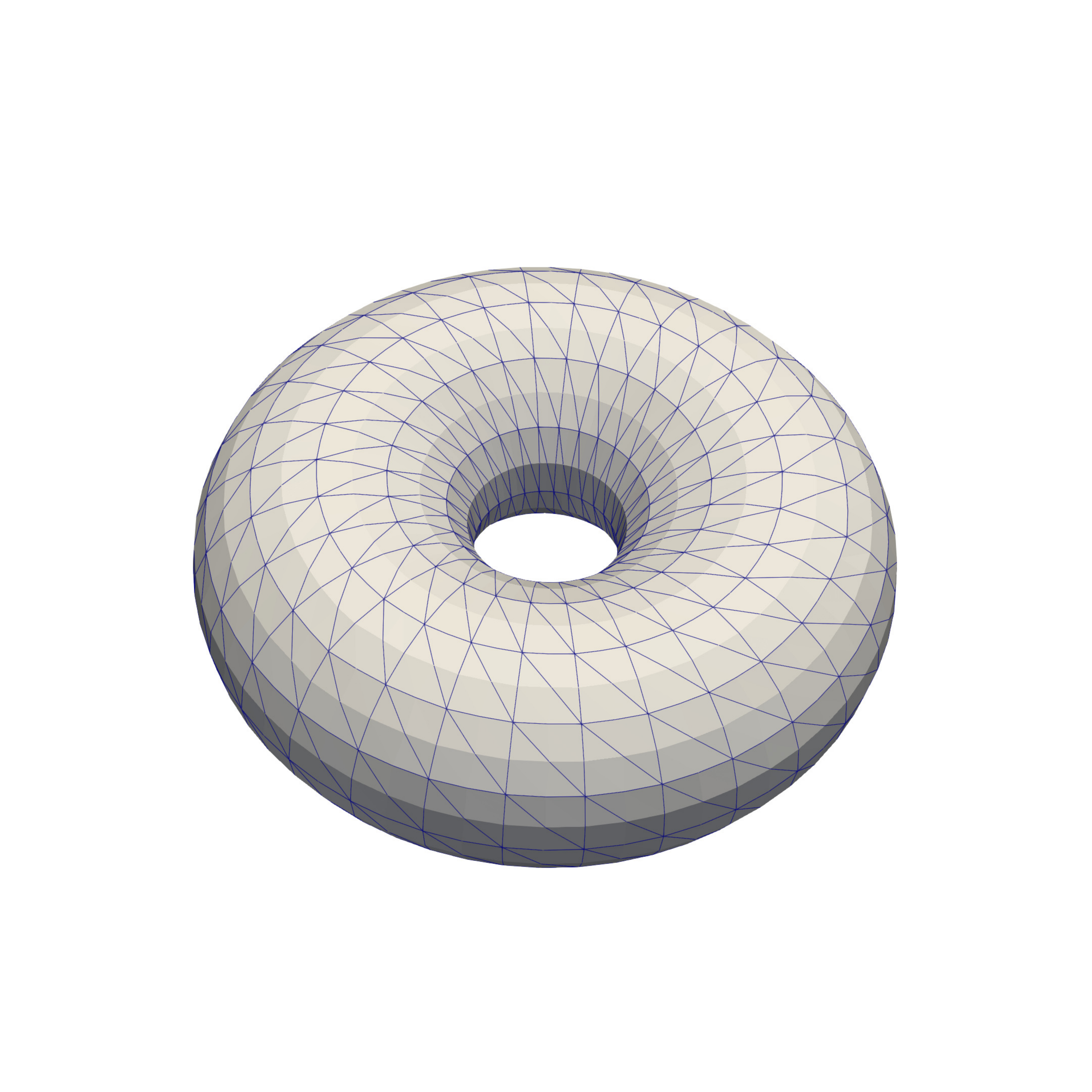}
\end{minipage}
\hspace{12pt}
\begin{minipage}[t]{0.20\linewidth}
	\centering
\includegraphics[scale=0.06]{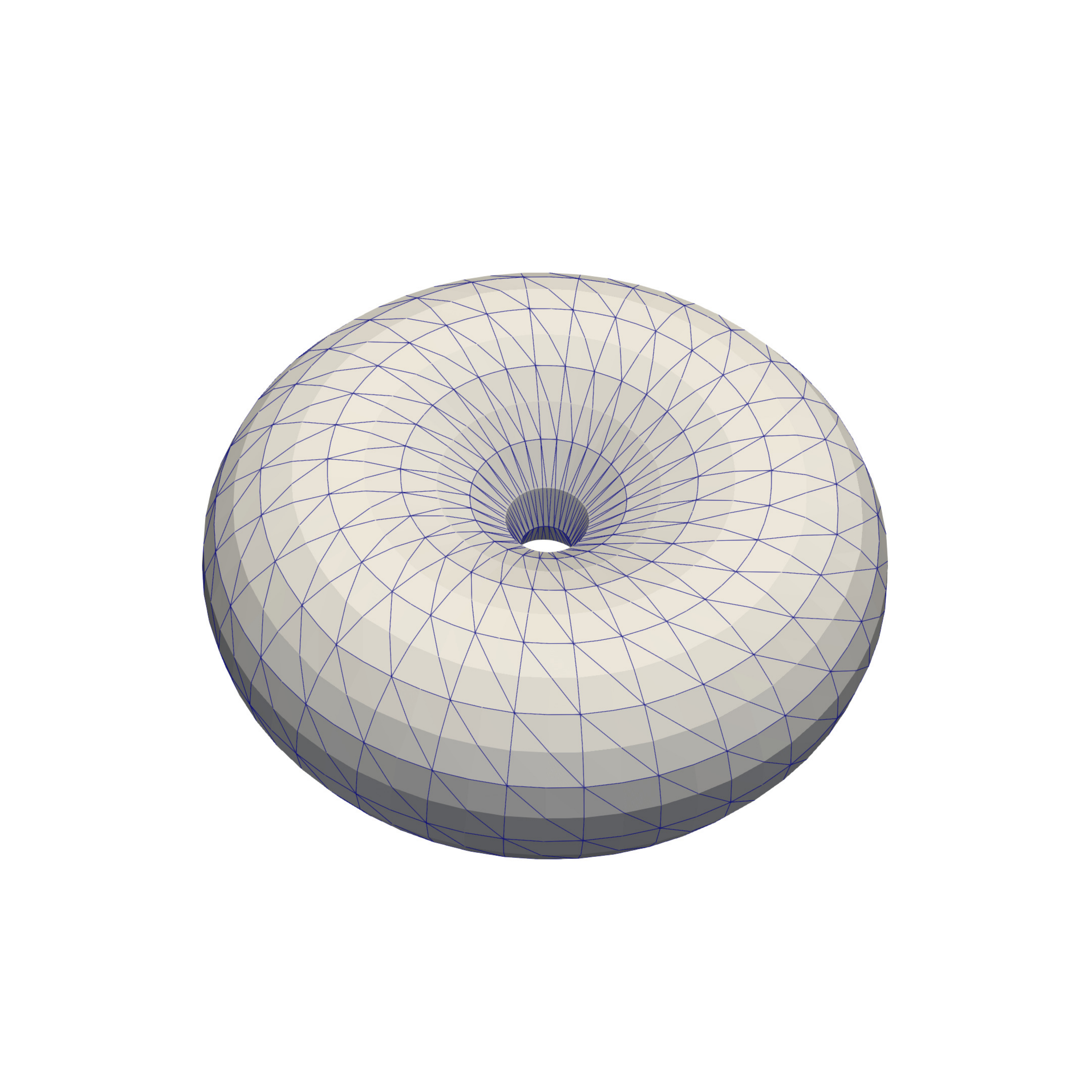}
\end{minipage}
\caption{Evolution of a torus by surface diffusion, using $\P^2$-elements. The total number of Lagrange nodes is $1440$.}
\label{sd_torus_l=2}
\end{figure}

\begin{figure}[h!]
	\centering
	\begin{subfigure}{.325\linewidth}
		\includegraphics[width=\textwidth]{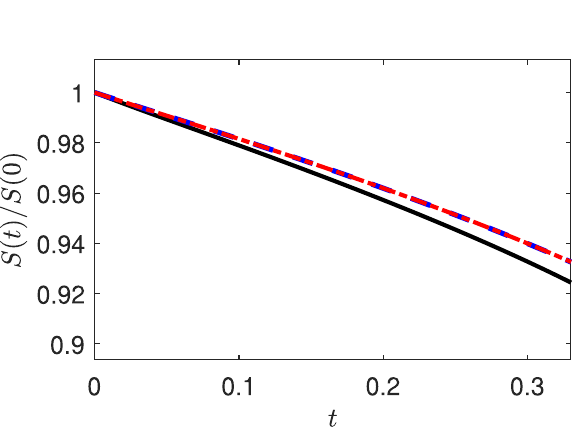}
		\caption{The normalized surface area}
	\end{subfigure}\hfill%
	\begin{subfigure}{.325\linewidth}
		\includegraphics[width=\textwidth]{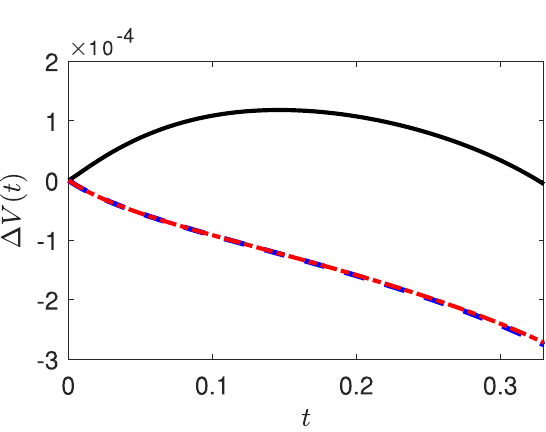}
		\caption{The relative volume loss}
	\end{subfigure}\hfill%
	\begin{subfigure}{.325\linewidth}
		\includegraphics[width=\textwidth]{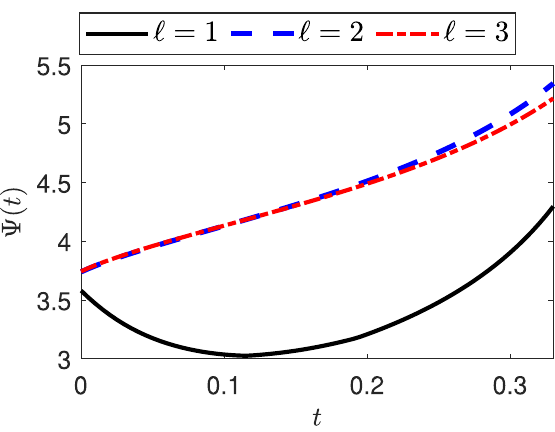}
		\caption{The mesh quality function}
	\end{subfigure}
	\caption{Evolution of geometric quantities for $\P^\ell$-elements, $\ell=1,2,3$.}
	\label{sd_torus_geo}
\end{figure}

To provide a more comprehensive comparison, we use the torus $(\sqrt{x^2+y^2}-2)^2+z^2=1$ as an initial surface and apply $\P^\ell$-elements to show the morphology evolution and monitor the evolution of some geometric quantities. \Cref{sd_torus_l=1,sd_torus_l=2} shows that the torus will develop a singularity in the center under surface diffusion. We stop the simulation before such a singularity would occur. The total number of Lagrange nodes of the $\P^1$-element method and the $\P^2$-element method is 360 and 1440, respectively. We notice that the surface area keeps decreasing, as shown in \Cref{sd_torus_geo}.

\subsection{Structure-preserving isoparametric finite element methods for surface diffusion}

In this subsection, we present some numerical results for structure-preserving isoparametric finite element methods for surface diffusion of both curves \eqref{SP:curve} and surfaces \eqref{SP:surface}. For these nonlinear methods, we use the Picard iteration method. Take the curves case \eqref{SP:curve} as an example. We first set
\[
\bX^{m+1,0}_h=\bX^{m}_h,\qquad \kappa^{m+1,0}_h=\kappa^{m}_h,
\]
and for each iteration step $i\ge 1$, we solve the following system of equations
	\begin{align*}
		\l(\frac{\bX^{m+1,i+1}_h-\id}{\tau},\varphi_h\bn^{m+\frac{1}{2},i}_h \r)^h_{\Gamma^m_h}-\l(\nabla_s\kappa^{m+1,i+1}_h,\nabla_s\varphi_h\r)^h_{\Gamma^m_h}&=0, \\
			\l(\kappa^{m+1,i+1}_h\bn^{m+\frac{1}{2},i}_h,\bm{\eta}_h \r)^h_{\Gamma^m_h} +\l(\nabla_s\bX^{m+1,i+1}_h,\nabla_s \bm{\eta}_h \r)_{\Gamma^m_h}^{h}&=0.
	\end{align*}
Here, the intermediate normal vector is given by
\begin{equation*}
	\bn^{m+\frac{1}{2},i}_h\circ\bY_h^m=\frac{(\p_{\widehat{\bm{\tau}}} \bY_h^{m+1,i}+\p_{\widehat{\bm{\tau}}} \bY_h^{m})^\perp}{2|\p_{\widehat{\bm{\tau}}} \bY_h^m|},\qquad \text{on}\quad \widehat{\sigma}_j.
\end{equation*}
Note that the above system has a very similar structure to the linear scheme \eqref{eq:2dnum}. 
We end the iteration step when the adjacent iteration error is smaller than a fixed number
\[
\|\bX^{m+1,i+1}_h-\bX^{m+1,i}_h\|_{L^{\infty}}+\|\kappa_h^{m+1,i+1}-\kappa_h^{m+1,i}\|_{L^{\infty}}\le \mathrm{tol}.
\]
 A similar Picard iteration method can be applied for the surfaces case \eqref{SP:surface}, which we omit here for simplicity. In practice, we always set $\mathrm{tol}=10^{-12}$.

\begin{figure}[h!]
	\centering
	\begin{subfigure}{.325\linewidth}
		\includegraphics[width=\textwidth]{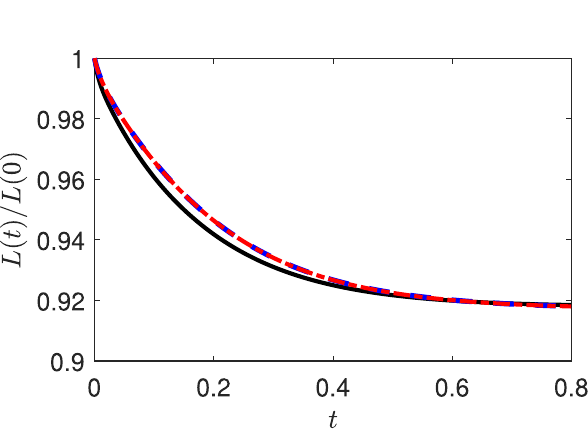}
		\caption{The normalized perimeter}
	\end{subfigure}\hfill%
	\begin{subfigure}{.325\linewidth}
		\includegraphics[width=\textwidth]{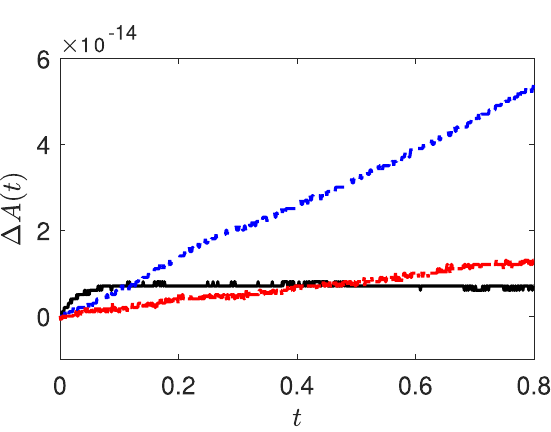}
		\caption{The relative area loss}
	\end{subfigure}\hfill%
	\begin{subfigure}{.325\linewidth}
		\includegraphics[width=\textwidth]{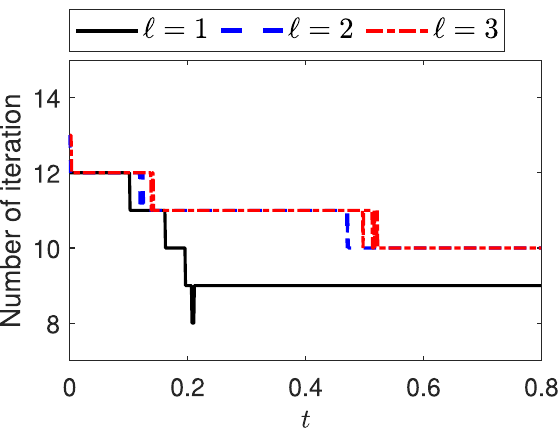}
		\caption{The number of iterations}
	\end{subfigure}
	\caption{Evolution of some geometric quantities for the structure-preserving $\P^\ell$-element scheme, $\ell=1,2,3$, starting with an ellipse. The discretization parameters are $N=128$ and $\tau=0.005$.}
	\label{curve_SP_geo}
\end{figure}

\begin{figure}[h!]
	\centering
	\begin{subfigure}{.325\linewidth}
		\includegraphics[width=\textwidth]{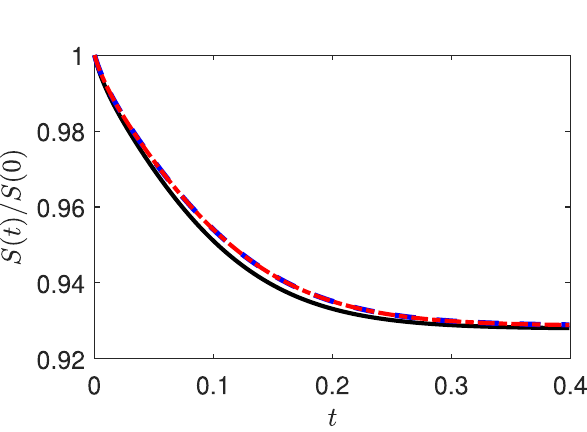}
		\caption{The normalized surface area}
	\end{subfigure}\hfill%
	\begin{subfigure}{.325\linewidth}
		\includegraphics[width=\textwidth]{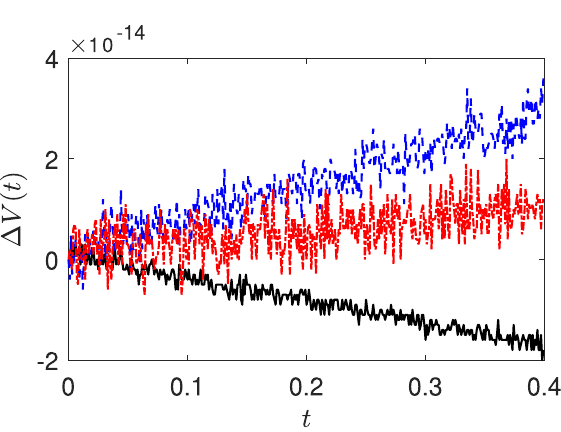}
		\caption{The relative volume loss}
	\end{subfigure}\hfill%
	\begin{subfigure}{.325\linewidth}
		\includegraphics[width=\textwidth]{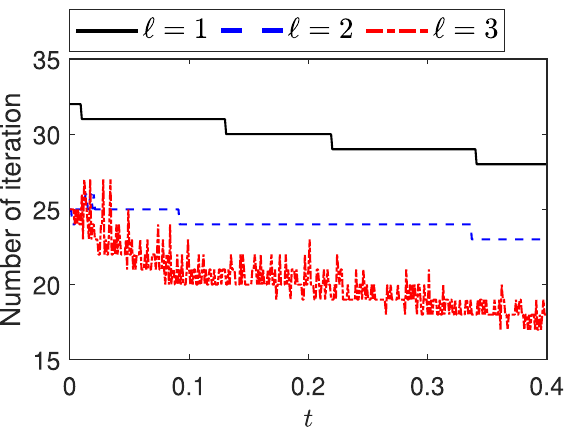}
		\caption{The number of iterations}
	\end{subfigure}
	\caption{Evolution of some geometric quantities for the structure-preserving $\P^\ell$-element scheme, $\ell=1,2,3$, starting with as an ellipsoid. The discretizaion parameters are $(J,K)=(676,340)$ and $\tau=0.001$.}
	\label{surface_SP_geo}
\end{figure}

Since the evolution of shapes and convergence order test by structure-preserving schemes makes not much difference with the schemes \eqref{eq:2dnum} and \eqref{eq:3dnum}, we only present the evolution of some geometric quantities. \Cref{curve_SP_geo,surface_SP_geo} depict the evolution of the normalized perimeter/surface area and the relative area loss/volume loss for an initial ellipse and an initial ellipsoid, respectively.  We observe that the structure-preserving schemes preserve two geometric structures for both the curves and the surfaces case, while the iteration numbers of the Picard iteration are relatively small (lower than 35 per time step).

 \begin{figure}[h!]
        \centering
        \begin{subfigure}[b]{0.4\linewidth}
            \includegraphics[width=7.5cm,height=4.5cm]{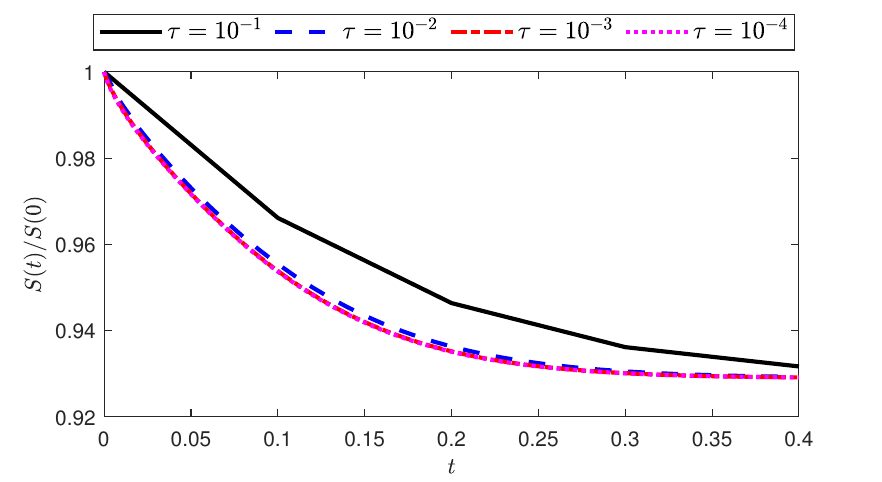}
             \caption[]%
            {{\small The normalized surface area}}    
        \end{subfigure}
        \hspace{1cm}
        \begin{subfigure}[b]{0.4\linewidth}  
            \includegraphics[width=7.5cm,height=4.5cm]{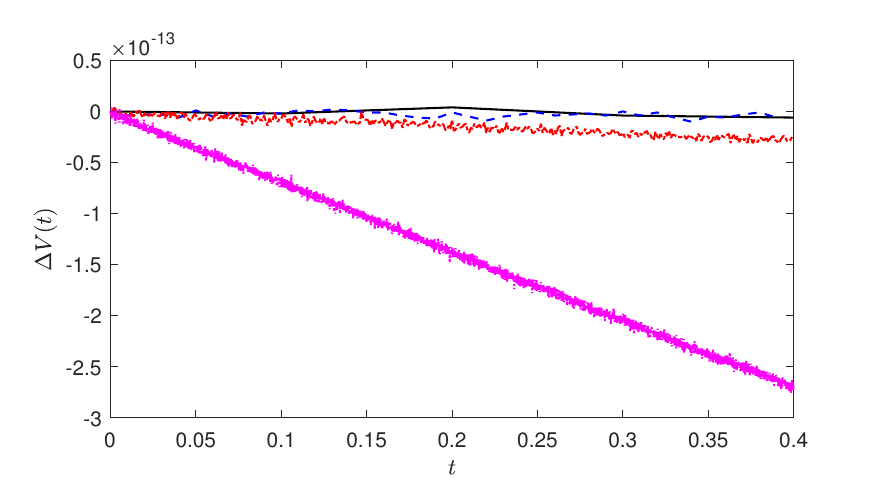}
            \caption[]%
            {{\small The relative volume loss}}    
        \end{subfigure}
        \begin{subfigure}[b]{0.4\linewidth}   
            \includegraphics[width=7.5cm,height=4.5cm]{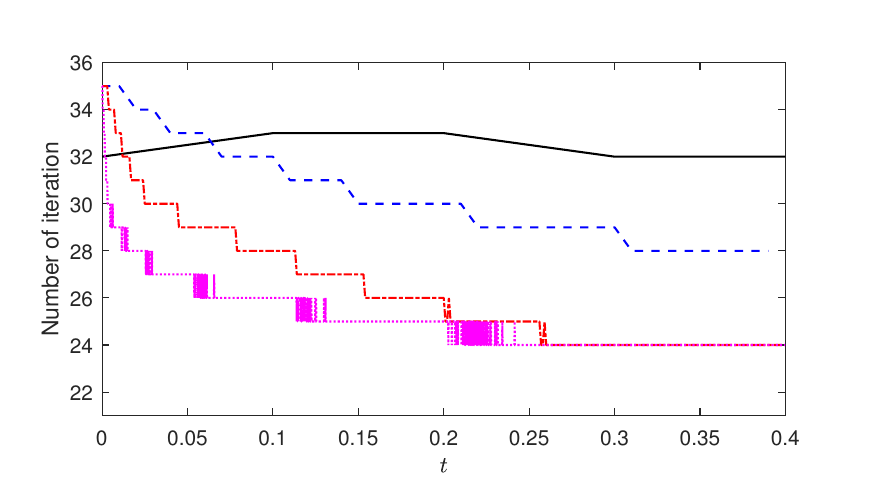}
            \caption[]%
            {{\small The number of iterations}}    
        \end{subfigure}
         \hspace{1cm}
        \begin{subfigure}[b]{0.4\linewidth}   
            \includegraphics[width=7.5cm,height=4.5cm]{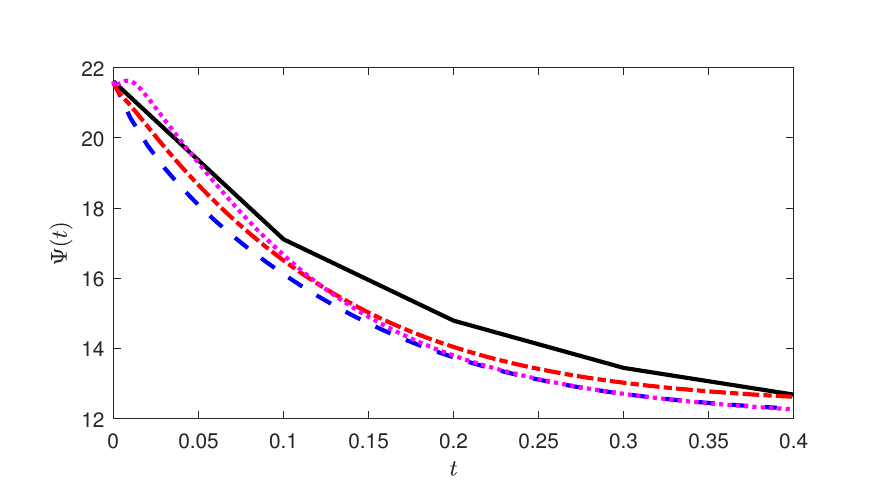}
          \caption[]%
            {{\small  The mesh ratio function}}    
        \end{subfigure}
        \caption[  ]
        {Evolution of some quantities for the structure-preserving $\mathbb{P}^2$-element scheme, starting with an ellipsoid. The discretization parameters are $(J,K)=(676,340)$ and $\tau=0.1,0.01,0.001, 0.0001$.} 
        \label{surface_SP_geo_dif_tau}
    \end{figure}

Finally, we show the robustness of our numerical method, and of the Picard iteration solver, with respect to large and small time step sizes. To this end, we compare different time step sizes ranging from $\tau=0.1$ to $\tau=0.0001$ for a fixed mesh size for the $\P^2$-finite element method. As shown in \Cref{surface_SP_geo_dif_tau}, the Picard iteration numbers remain low for all time step sizes, and the structure-preserving properties are maintained.

%

\section{Conclusions}

We have proposed unconditionally energy-stable isoparametric finite element methods for mean curvature flow and surface diffusion. Our methods can be seen as an extension of the classical BGN method \cite{BGN07A,BGN08A,BGN20}, which uses linear finite elements, to arbitrarily high order finite elements. Our approach is to fix a reference grid and prove a pointwise estimate over this grid. We implement our methods within \textsc{Dune} and demonstrate the unconditional energy stability for the discrete energy based on numerical integration. A similar idea can also be extended to structure-preserving methods for surface diffusion by choosing the intermediate normal \cite{Bao-Zhao}. Extensive numerical experiments demonstrate the expected higher spatial order accuracy, unconditional energy stability, and structure-preserving properties.

\section*{Acknowledgment}
S.\ Praetorius was partially supported by the German Research Foundation (DFG) through the research unit FOR3013, ``Vector- and Tensor-Valued Surface PDEs,'' within the project TP06, ``Symmetry, length, and tangential constraints'' (project number 417223351). G.\ Zhang was supported by the Tsinghua University Doctoral Student Short-Term Overseas Visiting Scholarship Fund.

\appendix
\renewcommand*{\appendixname}{}

\renewcommand{\theequation}{\Alph{section}.\arabic{equation}}
\section{Transformation formula for surface integrals}\label{sec:appendix:transform}

Let $\widehat{\sigma},\sigma,\widetilde{\sigma}\subset \R^3$ be three $C^1$-surfaces, possibly with a boundary. Let $\sigma$ and $\widetilde{\sigma}$ be parametrized by $\bY\colon\widehat{\sigma} \to \R^3$ and $\widetilde{\bY}=\bX\circ \bY\colon\widehat{\sigma} \to \R^3$, respectively, where $\bX\in C^1(\sigma,\R^3)$.
The surface elements $|\mathcal{J}(\widetilde{\bY})|$ and $|\mathcal{J}(\bY)|$ of $\widetilde{\sigma}$ and $\sigma$,
respectively, are defined via $\cJ(\widetilde\bY)=\p_{\widehat{\bm{\tau}}_1}\widetilde\bY\times \p_{\widehat{\bm{\tau}}_2}\widetilde\bY$ and
$\cJ(\bY)=\p_{\widehat{\bm{\tau}}_1}\bY\times \p_{\widehat{\bm{\tau}}_2}\bm Y$, with
$\{\widehat{\bm{\tau}}_1,\widehat{\bm{\tau}}_2\}$ being an orthonormal basis of the tangent space of $\widehat{\sigma}$.

\begin{lemma}\label{lemma:appendix}
 It holds that
\begin{equation} \label{eq:lemmaappendix0}
|\mathcal{J}(\widetilde{\bY})|= \sqrt{\det\l(\p_{\widehat{\bm{\tau}}_i}\widetilde{\bY}\cdot \p_{\widehat{\bm{\tau}}_j}\widetilde{\bY} \r)_{1\le i,j\le 2} }, 
\end{equation}
and
\begin{equation} \label{eq:lemmaappendix}
|\mathcal{J}(\widetilde{\bY})|= |\cJ(\bX)|\circ \bY\, |\mathcal{J}(\bY)| = \sqrt{\det\l(\p_{\bm{\tau}_i}\bX\cdot \p_{\bm{\tau}_j}\bX \r)_{1\le i,j\le 2} }\circ \bY\, |\mathcal{J}(\bY)|,
\end{equation}
where $\{\bm{\tau}_1,\bm{\tau}_2\}$ is an orthonormal basis of the tangent space of $\sigma$.
\end{lemma}

\begin{proof}
	For vectors $\bm a, \bm b, \bm c, \bm d \in \mathbb R^3$ it holds that
\begin{equation}
    (\bm a \times \bm b) \cdot (\bm c \times \bm d) = \det \begin{pmatrix} \bm a \cdot \bm c & \bm a \cdot \bm d \\ \bm b \cdot \bm c & \bm b \cdot \bm d \end{pmatrix} = (\bm a \cdot \bm d)\, \bm b \cdot \bm d - (\bm a \cdot \bm d)\, \bm b \cdot \bm c,
\end{equation}
which immediately implies \eqref{eq:lemmaappendix0}.

Locally we can find an open set $U \subset \mathbb R^2$ and a map $\widehat{\bm G} : U \to \widehat\sigma$ that parameterizes a part of $\widehat\sigma$. Let $\bm G = \bm Y \circ \widehat{\bm G}$ and
$\widetilde{\bm G} = \widetilde{\bm Y} \circ \widehat{\bm G} = \bm X \circ \bm Y \circ \widehat{\bm G}$ be the corresponding parameterizations of parts of $\sigma$ and $\widetilde\sigma$, respectively.
It is then not difficult to prove that 
\begin{equation} \label{eq:appflat}
|\cJ(\widetilde{\bm G})| = |\cJ(\bm X)| \circ \bm G\, |\cJ(\bm G)|,    
\end{equation}
where 
$|\cJ(\widetilde{\bm G})|^2 = \det ((\nabla \widetilde{\bm G})^\top \nabla \widetilde{\bm G})$ and
$|\cJ(\bm G)|^2 = \det ((\nabla \bm G)^\top \nabla \bm G)$. 
Indeed, we know that $\nabla \bm G$ maps $\R^2$ to the tangent space of
$\sigma$. Hence we can choose $\bm b_i \in \R^2$ with $(\nabla \bm G) \bm b_i =
\bm\tau_i$, $i=1,2$. It follows that $(\nabla\widetilde{\bm G}) \bm b_i =
\p_{\bm\tau_i} \bm X$, $i=1,2$.  We then compute
\begin{align*}
(\p_{\bm{\tau}_i}\bX\cdot \p_{\bm{\tau}_j}\bX ) \circ \bm G
= ((\nabla\widetilde{\bm G}) \bm b_i) \cdot ((\nabla\widetilde{\bm G}) \bm b_j)
= \bm b_i \cdot ((\nabla\widetilde{\bm G})^\top(\nabla\widetilde{\bm G}) \bm b_j)
= \bm e_i \cdot (B^\top (\nabla\widetilde{\bm G})^\top(\nabla\widetilde{\bm G}) B)
\bm e_j,
\end{align*}
where $B \in \R^{2\times2}$ is the matrix with columns $\bm b_1$ and $\bm b_2$.
Hence, we have that
\begin{equation} \label{eq:appJX2}
|\cJ(\bm X)|^2 \circ \bm G = \det\l(\p_{\bm{\tau}_i}\bX\cdot \p_{\bm{\tau}_j}\bX \r)_{1\le i,j\le 2} = \det(B^2) \det ( (\nabla\widetilde{\bm G})^\top\nabla\widetilde{\bm G})=\det(B^2) |\cJ(\widetilde{\bm G})|^2.
\end{equation}
On the other hand, we have
\[
\bm \tau_i \cdot \bm \tau_j = ((\nabla \bm G) \bm b_i) \cdot 
(\nabla \bm G) \bm b_j = \bm e_i \cdot (B^\top (\nabla{\bm G})^\top(\nabla{\bm G}) B) \bm e_j,
\]
which implies that $1 = \det(B^2) \det ( (\nabla{\bm G})^\top\nabla{\bm G})
= \det(B^2) |\cJ(\bm G)|^2$.
Combining this with \eqref{eq:appJX2} yields the claim \eqref{eq:appflat}. 

Now applying \eqref{eq:appflat} two times, to different function compositions, yields that
\[
|\cJ(\widetilde{\bm G})| = |\cJ(\widetilde{\bm Y})| \circ \widehat{\bm G}\, |\cJ(\widehat{\bm G})|
\quad\text{and} \quad
|\cJ(\bm G)| = |\cJ(\bm Y)| \circ \widehat{\bm G}\, |\cJ(\widehat{\bm G})|.
\]
Combining these gives our desired claim \eqref{eq:lemmaappendix} first on a patch, and with the help of a complete chart for all of $\widehat\sigma$.\end{proof}

\renewcommand{\theequation}{\Alph{section}.\arabic{equation}}
\section{Numerical computation of the Hausdorff distance}\label{sec:appendix:hausdorff-distance}

In this appendix, we give precise details on how we compute the distance used in \eqref{L2error}. We implemented the presented algorithm within the numerical framework \textsc{Dune} \cite{Bastian2021}.

Let $\Gamma_1\subset\mathbb{R}^{n+1}$ and $\Gamma_2\subset\mathbb{R}^{n+1}$ denote two parametrized $n$-dimensional grids, $n=1,2$, with $\boldsymbol{Y}_1$ and $\boldsymbol{Y}_2$ their parametrization over corresponding reference grids $\widehat{\Gamma}_1=\bigcup_{j=1}^{J_1}\widehat{\sigma}_{1,j}=\bigcup_{j=1}^{J_1} A_{1,j}(\widehat{\mathbb{K}})$ and $\widehat{\Gamma}_2=\bigcup_{j=1}^{J_2}\widehat{\sigma}_{2,j}=A_{2,j}(\widehat{\mathbb{K}})$.
We consider the approximate $L^2$-distance $d_{L^2}(\Gamma_1, \Gamma_2)$ in terms of the pointwise closest-point distance function $\operatorname{dist}(\cdot, \Gamma)$,
\begin{align} \label{eq:dL2}
    d_{L^2}^2(\Gamma_1, \Gamma_2) &\colonequals \sum_{j=1}^{J_1}\sum_{i=1}^Q\operatorname{dist}^2(\bF_{1,j}(\qbx_i), \Gamma_2)\, |\nabla_{\bxHat}\bF_{1,j}|(\qbx_i)\, \qw_i\,,
\end{align}
with $J_1$ the number of elements in the reference grid $\widehat{\Gamma}_1$ and $Q$ the number of quadrature points used to approximate the integral, where $\qbx_i$ are the quadrature points and $\qw_i$ the quadrature weights of a quadrature rules to use for the approximation, and $\bF_{1,j}=\bY_1\circ A_{1,j}$.

While this problem has been investigated for piecewise flat triangular grids,  \cite{Cignoni1998Metro,Straub2007Exact,Bartin2010Precise}, point clouds \cite{Chen2017Local,Zhang2018Efficient}, or freeform surfaces \cite{Elber2008Hausdorff}, little can be found about the computation of the distance between parametrized grids.
The challenge lies in the efficient computation of the closest-point distance to the parametrized surface. It can be written in its squared form as
\begin{align}
    \operatorname{dist}^2(\boldsymbol{p},\Gamma_2) &=  \inf_{j=1}^{J_2}\inf_{\bxHat\in\widehat{\mathbb{K}}} \|p - \bF_{2,j}(\bxHat)\|^2
\end{align}
for $\boldsymbol{p}\in\Gamma_1$, with $\bF_{2,j}=\bY_2\circ A_{2,j}$. The inner $\inf$-problem is a constrained least-squares problem, with linear constraints restricting the point $\bxHat$ to the domain of $\widehat{\mathbb{K}}$. Since $\widehat{\mathbb{K}}$ is the standard $n$-simplex, we have the box-constraints $0\leq \widehat{x}_i\leq 1$, $i=1\ldots n$, and additionally the barycentric coordinate condition $\sum_{i=1}^n \widehat{x}_i \leq 1$. 
The actual least-squares problem is then solved using a Levenberg--Marquardt algorithm, implemented in the ALGLIB library \cite{alglib}.

The outer $\inf$-problem is simplified by searching only in elements that are close to the target point $\boldsymbol{p}$. We thereby represent each grid element of $\Gamma_2$ by its center in the sense of the mapping of the barycenter in the reference element $\bxHat_c\in\widehat{\mathbb{K}}$, to build a point cloud:
\begin{align}
    \chi_2 \colonequals\left\{\boldsymbol{q}=\bF_{2,j}(\bxHat_c)\mid j=1\ldots J_2\right\}\,.
\end{align}
The set $\operatorname{kNN}(\boldsymbol{p},\chi_2)$ of $k$ nearest neighbors to the point $\boldsymbol{p}$, represented in terms of their indices $j\in\{1\ldots J_2\}$, is computed using a kd-tree data structure and the nanoflann library \cite{nanoflann}. We thus obtain the procedure
\begin{align}
    \operatorname{dist}^2(\boldsymbol{p},\Gamma_2) &= \inf_{j\in\operatorname{kNN}(\boldsymbol{p},\chi_2)} \inf_{\bxHat\in\widehat{\mathbb{K}}} \|\boldsymbol{p} - \bF_{2,j}(\bxHat)\|^2\,.
\end{align}

Denote by $L$ the maximal number of Levenberg--Marquardt iterations. The total complexity of the Hausdorff-distance computation can thus be estimated as $\mathcal{O}(J_1\cdot Q\cdot k\log(J_2)\cdot L)$ plus an additional initialization phase of the kd-tree data structures, which is independent of $J_1$ and $Q$.


\vspace{-8pt}

\bibliographystyle{abbrv}
\bibliography{bibliography}

\begin{thebibliography}{10}

\bibitem{Bai-Hu-Li}
G.~Bai, J.~Hu, and B.~Li.
\newblock A convergent evolving finite element method with artificial tangential motion for surface evolution under a prescribed velocity field.
\newblock {\em SIAM Journal on Numerical Analysis}, 62(5):2172--2195, 2024.

\bibitem{Bai2024}
G.~Bai and B.~Li.
\newblock Convergence of a stabilized parametric finite element method of the {Barrett--Garcke--N{\"u}rnberg} type for curve shortening flow.
\newblock {\em Mathematics of Computation}, 2024.

\bibitem{Bai2023}
G.~Bai and B.~Li.
\newblock A new approach to the analysis of parametric finite element approximations to mean curvature flow.
\newblock {\em Foundations of Computational Mathematics}, 24(5):1673--1737, 2024.

\bibitem{Bansch}
E.~B{\"a}nsch, K.~Deckelnick, H.~Garcke, and P.~Pozzi.
\newblock {\em Interfaces: modeling, analysis, numerics}.
\newblock Springer, 2023.

\bibitem{Bao-Garcke-Nurnberg-Zhao2022}
W.~Bao, H.~Garcke, R.~N{\"u}rnberg, and Q.~Zhao.
\newblock Volume-preserving parametric finite element methods for axisymmetric geometric evolution equations.
\newblock {\em Journal of Computational Physics}, 460:111180, 2022.

\bibitem{Bao-Jiang-Li}
W.~Bao, W.~Jiang, and Y.~Li.
\newblock A symmetrized parametric finite element method for anisotropic surface diffusion of closed curves.
\newblock {\em SIAM Journal on Numerical Analysis}, 61(2):617--641, 2023.

\bibitem{Bao-Li2023}
W.~Bao and Y.~Li.
\newblock A symmetrized parametric finite element method for anisotropic surface diffusion in three dimensions.
\newblock {\em SIAM Journal on Scientific Computing}, 45(4):A1438--A1461, 2023.

\bibitem{Bao-Li2024}
W.~Bao and Y.~Li.
\newblock A structure-preserving parametric finite element method for geometric flows with anisotropic surface energy.
\newblock {\em Numerische Mathematik}, 156(2):609--639, 2024.

\bibitem{Bao-Zhao}
W.~Bao and Q.~Zhao.
\newblock A structure-preserving parametric finite element method for surface diffusion.
\newblock {\em SIAM Journal on Numerical Analysis}, 59(5):2775--2799, 2021.

\bibitem{BGN07B}
J.~W. Barrett, H.~Garcke, and R.~N{\"u}rnberg.
\newblock On the variational approximation of combined second and fourth order geometric evolution equations.
\newblock {\em SIAM Journal on Scientific Computing}, 29(3):1006--1041, 2007.

\bibitem{BGN07A}
J.~W. Barrett, H.~Garcke, and R.~N{\"u}rnberg.
\newblock A parametric finite element method for fourth order geometric evolution equations.
\newblock {\em Journal of Computational Physics}, 222(1):441--467, 2007.

\bibitem{BGN08A}
J.~W. Barrett, H.~Garcke, and R.~N{\"u}rnberg.
\newblock On the parametric finite element approximation of evolving hypersurfaces in {$\Bbb R^3$}.
\newblock {\em Journal of Computational Physics}, 227(9):4281--4307, 2008.

\bibitem{BGN20}
J.~W. Barrett, H.~Garcke, and R.~N{\"u}rnberg.
\newblock Parametric finite element approximations of curvature-driven interface evolutions.
\newblock In {\em Handbook of numerical analysis}, volume~21, pages 275--423. Elsevier, 2020.

\bibitem{Bartin2010Precise}
M.~Barto{\v{n}}, I.~Hanniel, G.~Elber, and M.-S. Kim.
\newblock Precise {Hausdorff} distance computation between polygonal meshes.
\newblock {\em Computer Aided Geometric Design}, 27(8):580--591, 2010.

\bibitem{dune}
P.~Bastian, M.~Blatt, S.~Burballa, A.~Burchardt, A.~Dedner, C.~Engwer, J.~Fahlke, C.~Gr{\"a}ser, C.~Gr{\"u}ninger, D.~Kempf, R.~Kl{\"o}fkorn, T.~Koch, S.~Ospina~De~Los~R{\'i}os, S.~Praetorius, and O.~Sander.
\newblock \textsc{Dune} -- the distributed and unified numerics environment.
\newblock \url{https://dune-project.org}, 2024.
\newblock v2.10.

\bibitem{Bastian2021}
P.~Bastian, M.~Blatt, A.~Dedner, N.-A. Dreier, C.~Engwer, R.~Fritze, C.~Gr{\"a}ser, C.~Gr{\"u}ninger, D.~Kempf, R.~Kl{\"o}fkorn, et~al.
\newblock The {Dune} framework: Basic concepts and recent developments.
\newblock {\em Computers \& Mathematics with Applications}, 81:75--112, 2021.

\bibitem{Binz}
T.~Binz and B.~Kov{\'a}cs.
\newblock A convergent finite element algorithm for generalized mean curvature flows of closed surfaces.
\newblock {\em IMA Journal of Numerical Analysis}, 42(3):2545--2588, 2022.

\bibitem{nanoflann}
J.~L. Blanco and P.~K. Rai.
\newblock nanoflann: a {C}++ header-only fork of {FLANN}, a library for nearest neighbor ({NN}) with kd-trees.
\newblock \url{https://github.com/jlblancoc/nanoflann}, 2014.
\newblock Version v1.6.2.

\bibitem{alglib}
S.~Bochkanov.
\newblock {ALGLIB}.
\newblock \url{https://www.alglib.net}.
\newblock Version v4.03.0.

\bibitem{Brenner-Scott}
S.~C. Brenner and L.~R. Scott.
\newblock {\em The mathematical theory of finite element methods}.
\newblock Springer, 2008.

\bibitem{Chen2017Local}
Y.~Chen, F.~He, Y.~Wu, and N.~Hou.
\newblock A local start search algorithm to compute exact {Hausdorff} distance for arbitrary point sets.
\newblock {\em Pattern Recognition}, 67:139--148, 2017.

\bibitem{Cignoni1998Metro}
P.~Cignoni, C.~Rocchini, and R.~Scopigno.
\newblock {Metro}: measuring error on simplified surfaces.
\newblock {\em Computer Graphics Forum}, 17(2):167--174, 1998.

\bibitem{Deckelnick-Dziuk-Elliott}
K.~Deckelnick, G.~Dziuk, and C.~M. Elliott.
\newblock Computation of geometric partial differential equations and mean curvature flow.
\newblock {\em Acta numerica}, 14:139--232, 2005.

\bibitem{csftime}
K.~Deckelnick and R.~N\"urnberg.
\newblock Second order in time finite element schemes for curve shortening flow and curve diffusion.
\newblock arXiv:2502.19277, 2025.

\bibitem{Demlow}
A.~Demlow.
\newblock Higher-order finite element methods and pointwise error estimates for elliptic problems on surfaces.
\newblock {\em SIAM Journal on Numerical Analysis}, 47(2):805--827, 2009.

\bibitem{Duan2021}
B.~Duan, B.~Li, and Z.~Zhang.
\newblock High-order fully discrete energy diminishing evolving surface finite element methods for a class of geometric curvature flows.
\newblock {\em Annals of Applied Mathematics}, 37(4):405--436, 2021.

\bibitem{Dziuk1991}
G.~Dziuk.
\newblock An algorithm for evolutionary surfaces.
\newblock {\em Numerische Mathematik}, 58(1):603--611, 1990.

\bibitem{Elber2008Hausdorff}
G.~Elber, T.~Grandine, F.~Chen, and B.~J{\"u}ttler.
\newblock {Hausdorff} and minimal distances between parametric freeforms in $\mathbb{R}^2$ and $\mathbb{R}^3$.
\newblock In {\em Advances in Geometric Modeling and Processing}, pages 191--204, Berlin, Heidelberg, 2008. Springer Berlin Heidelberg.

\bibitem{Elliott2013}
C.~M. Elliott and T.~Ranner.
\newblock Finite element analysis for a coupled bulk--surface partial differential equation.
\newblock {\em IMA Journal of Numerical Analysis}, 33(2):377--402, 2013.

\bibitem{Elliott}
C.~M. Elliott and T.~Ranner.
\newblock A unified theory for continuous-in-time evolving finite element space approximations to partial differential equations in evolving domains.
\newblock {\em IMA Journal of Numerical Analysis}, 41(3):1696--1845, 2021.

\bibitem{dune-foamgrid}
B.~Flemisch, T.~Koch, N.~Schr{\"o}der, and O.~Sander.
\newblock {Dune-FoamGrid} -- implements one- and two-dimensional simplex grids, which can be embedded in {Euclidean} spaces of arbitrary dimension.
\newblock \url{https://gitlab.dune-project.org/extentsions/dune-foamgrid}, 2024.
\newblock v2.10.

\bibitem{Garcke2024}
H.~Garcke, W.~Jiang, C.~Su, and G.~Zhang.
\newblock Structure-preserving parametric finite element method for curve diffusion based on {Lagrange} multiplier approaches.
\newblock {\em arXiv preprint arXiv:2408.13443}, 2024.

\bibitem{higher-order-bgn}
H.~Garcke, R.~N{\"u}rnberg, S.~Praetorius, and G.~Zhang.
\newblock Code and data: Isoparametric finite element methods for mean curvature flow and surface diffusion.
\newblock \url{https://github.com/Ganghui-Zhang/higher-order-bgn}, 2025.

\bibitem{Garcke-Nurnberg-Zhao2023}
H.~Garcke, R.~N{\"u}rnberg, and Q.~Zhao.
\newblock Structure-preserving discretizations of two-phase {Navier--Stokes} flow using fitted and unfitted approaches.
\newblock {\em Journal of Computational Physics}, 489:112276, 2023.

\bibitem{Heine03PhD}
C.-J. Heine.
\newblock {\em Computations of Form and Stability of Rotating Drops with Finite Elements}.
\newblock PhD thesis, University of Aachen, Aachen, 2003.

\bibitem{Heine04}
C.-J. Heine.
\newblock Isoparametric finite element approximation of curvature on hypersurfaces.
\newblock Technical Report~26, Fak. f. Math. Phys. Univ. Freiburg, 2004.

\bibitem{Hu-Li}
J.~Hu and B.~Li.
\newblock Evolving finite element methods with an artificial tangential velocity for mean curvature flow and {Willmore} flow.
\newblock {\em Numerische Mathematik}, 152(1):127--181, 2022.

\bibitem{Jiang21}
W.~Jiang and B.~Li.
\newblock A perimeter-decreasing and area-conserving algorithm for surface diffusion flow of curves.
\newblock {\em Journal of Computational Physics}, 443:110531, 2021.

\bibitem{Jiang23}
W.~Jiang, C.~Su, and G.~Zhang.
\newblock A second-order in time, {BGN}-based parametric finite element method for geometric flows of curves.
\newblock {\em Journal of Computational Physics}, 514:113220, 2024.

\bibitem{Kovacs}
B.~Kov{\'a}cs.
\newblock High-order evolving surface finite element method for parabolic problems on evolving surfaces.
\newblock {\em IMA Journal of Numerical Analysis}, 38(1):430--459, 2018.

\bibitem{Kovacs-Li-Lubich2019}
B.~Kov{\'a}cs, B.~Li, and C.~Lubich.
\newblock A convergent evolving finite element algorithm for mean curvature flow of closed surfaces.
\newblock {\em Numerische Mathematik}, 143:797--853, 2019.

\bibitem{Kovacs-Li-Lubich2020}
B.~Kov{\'a}cs, B.~Li, and C.~Lubich.
\newblock A convergent algorithm for forced mean curvature flow driven by diffusion on the surface.
\newblock {\em Interfaces and Free Boundaries}, 22(4):443--464, 2020.

\bibitem{Kovacs-Li-Lubich2021}
B.~Kov{\'a}cs, B.~Li, and C.~Lubich.
\newblock A convergent evolving finite element algorithm for {Willmore} flow of closed surfaces.
\newblock {\em Numerische Mathematik}, 149(3):595--643, 2021.

\bibitem{Li1}
B.~Li.
\newblock Convergence of {Dziuk}'s linearly implicit parametric finite element method for curve shortening flow.
\newblock {\em SIAM Journal on Numerical Analysis}, 58(4):2315--2333, 2020.

\bibitem{Li2}
B.~Li.
\newblock Convergence of {Dziuk}'s semidiscrete finite element method for mean curvature flow of closed surfaces with high-order finite elements.
\newblock {\em SIAM Journal on Numerical Analysis}, 59(3):1592--1617, 2021.

\bibitem{Nurnberg}
R.~N{\"u}rnberg.
\newblock A structure preserving front tracking finite element method for the {Mullins--Sekerka} problem.
\newblock {\em Journal of Numerical Mathematics}, 31(2):137--155, 2023.

\bibitem{dune-curvedgeometry}
S.~Praetorius.
\newblock {Dune-CurvedGeometry} -- parametrizations of curved geometries.
\newblock \url{https://gitlab.dune-project.org/extentsions/dune-curvedgeometry}, 2024.
\newblock v2.10.

\bibitem{dune-curvedgrid}
S.~Praetorius.
\newblock {Dune-CurvedGrid} -- parametrizations of curved grids.
\newblock \url{https://gitlab.dune-project.org/extentsions/dune-curvedgrid}, 2024.
\newblock v2.10.

\bibitem{Praetorius2022}
S.~Praetorius and F.~Stenger.
\newblock {Dune-CurvedGrid} -- a dune module for surface parametrization.
\newblock {\em Archive of Numerical Software}, 6:1--27, 2022.

\bibitem{Sander2020}
O.~Sander.
\newblock {\em \textsc{Dune} --- The distributed and unified numerics environment}, volume 140.
\newblock Springer Nature, 2020.

\bibitem{Sander2017}
O.~Sander, T.~Koch, N.~Schr{\"o}der, and B.~Flemisch.
\newblock The {Dune FoamGrid} implementation for surface and network grids.
\newblock {\em Archive of Numerical Software}, 5:217--244, 2017.

\bibitem{Schmidt93}
A.~Schmidt.
\newblock {\em {D}ie {B}erechnung dreidimensionaler {D}endriten mit {F}initen {E}lementen}.
\newblock PhD thesis, University Freiburg, Freiburg, 1993.

\bibitem{Schmidt96}
A.~Schmidt.
\newblock Computation of three dimensional dendrites with finite elements.
\newblock {\em J. Comput. Phys.}, 195(2):293--312, 1996.

\bibitem{Straub2007Exact}
R.~Straub.
\newblock Exact computation of the {Hausdorff} distance between triangular meshes.
\newblock In {\em Eurographics (Short Papers)}, pages 17--20, 2007.

\bibitem{Zhang2018Efficient}
D.~Zhang, L.~Zou, Y.~Chen, and F.~He.
\newblock Efficient and accurate {Hausdorff} distance computation based on diffusion search.
\newblock {\em IEEE Access}, 6:1350--1361, 2017.

\bibitem{Zhao2021}
Q.~Zhao, W.~Jiang, and W.~Bao.
\newblock An energy-stable parametric finite element method for simulating solid-state dewetting.
\newblock {\em IMA Journal of Numerical Analysis}, 41(3):2026--2055, 2021.

\end{thebibliography}

\end{document}